\documentclass[11pt]{amsart}
\usepackage{amsmath, amscd, amssymb}
\usepackage{times}
\usepackage[enableskew,vcentermath]{youngtab}
\usepackage{ytableau}
\usepackage{xcolor}
\ytableausetup{centertableaux, boxsize=1em}
\usepackage[mathscr]{eucal}
\usepackage[normalem]{ulem} 
\usepackage{soul} 
\usepackage[all,cmtip]{xy}
\usepackage{tikz-cd}   
\usepackage{mathabx}

\usepackage{color} 
\setulcolor{red}
\sethlcolor{yellow}
\setstcolor{blue}
\usepackage{amscd}
\usepackage{xfrac} 
\usepackage{yfonts}
\usepackage{ytableau}
\usepackage{geometry}
\usepackage{amsfonts}
\usepackage[colorlinks,linkcolor=blue,citecolor=blue, pdfstartview=FitH]{hyperref}
\usepackage{graphicx}
\usepackage{cite}

\allowdisplaybreaks
\numberwithin{equation}{section}

\theoremstyle{definition}
\newtheorem{thm}{Theorem}[section]
\newtheorem{theorem}[thm]{Theorem}

\newtheorem{lemma}[thm]{Lemma}
\newtheorem{corollary}[thm]{Corollary}
\newtheorem{proposition}[thm]{Proposition}

\newtheorem{remark}[thm]{Remark}

\newtheorem{definition}[thm]{Definition}

\newtheorem{example}[thm]{Example}

\newtheorem{defn-thm}[thm]{Definition-Theorem}
\newtheorem{conjecture}[thm]{Conjecture}

\newenvironment{observe}{\noindent\textcolor{blue}{\textit{Observation}}.}{\hfill \textcolor{blue}{$\blacktriangleleft$}\par}
\newtheorem*{theorem*}{Theorem}
\newtheorem*{proposition*}{Proposition}
\usepackage{todonotes}
\definecolor{pistachio}{rgb}{0.58, 0.77, 0.45}
\definecolor{eggshell}{rgb}{0.94, 0.92, 0.84}

\newcommand{\CC}{{\mathbb C}}

\newcommand{\Ext}{\operatorname{Ext}}

\newcommand{\Hom}{\operatorname{Hom}}

\newcommand{\cont}{{\operatorname{cont}}}

\newcommand{\Res}{\operatorname{Res}}

\newcommand{\btheorem}{\begin{theorem}}
\newcommand{\etheorem}{\end{theorem}}
\newcommand{\bproposition}{\begin{proposition}}
\newcommand{\eproposition}{\end{proposition}}
\newcommand{\bdefinition}{\begin{definition}}
\newcommand{\edefinition}{\end{definition}}
\newcommand{\bcorollary}{\begin{corollary}}
\newcommand{\ecorollary}{\end{corollary}}
\newcommand{\bproof}{\begin{proof}}
\newcommand{\eproof}{\end{proof}}
\newcommand{\bremark}{\begin{remark}}
\newcommand{\eremark}{\end{remark}}
\newcommand{\eexample}{\end{example}}
\newcommand{\bexample}{\begin{example}}
\newcommand{\elemma}{\end{lemma}}
\newcommand{\blemma}{\begin{lemma}}
\newcommand{\bobserve}{\begin{observe}}
\newcommand{\eobserve}{\end{observe}}

\renewcommand{\phi}{\varphi}
\newcommand{\ee}{\end{eqnarray*}}
\newcommand{\be}{\begin{eqnarray*}}

\newcommand{\beq}{\begin{equation}}
\newcommand{\eeq}{\end{equation}}
\newcommand{\bd}{\begin{enumerate}}
\newcommand{\ed}{\end{enumerate}}
\newcommand{\bti}{\begin{tikzcd}}
\newcommand{\eti}{\end{tikzcd}}

\renewcommand{\tilde}{\widetilde}

\title{On Ext-groups between Arthur-type modules of the graded Hecke algebra of type $A$}

\author{Kei Yuen Chan and Peng Zhou}
\address{Department of Mathematics, The University of Hong Kong,  Hong Kong
}
\email{kychan1@hku.hk}
\address{Department of Mathematics, The University of Hong Kong,  Hong Kong
}
\email{zhpmath@connect.hku.hk}

\date{}

\begin{document}

\begin{abstract}
Let $\mathbb H_n$ be the graded Hecke algebra of type $A_{n-1}$. This article establishes a formula for Ext-groups between a Speh $\mathbb H_n$-module and a $W$-type $\mathbb H_n$-module in terms of Littlewood-Richardson coefficients. As a consequence, we obtain an explicit and efficient method to compute Ext-groups between a Speh $\mathbb H_n$-module and an Arthur-type module of $\mathbb H_n$. In particular, this proves new instances of a conjecture of Chan and Prasad on Ext non-vanishing between Arthur-type representations. We also discuss some other non-Ext vanishing results for calibrated modules.
\end{abstract}

%\begin{abstract}
%    In this paper, we compute the extension groups $\textnormal{Ext}^{i}_{\mathbb{H}}(X, Y)$ of graded Hecke algebra $\mathbb{H}$, where $X$ is a Speh representation and $Y$ is a $W$-type representation. The primary methodology relies on decomposing the relevant Hom spaces and differentials using combinatorial tools, particularly Littlewood-Richardson coefficients and  Gelfand-Tsetlin basis. \\

%    mention a conjectural generalization of generalized Steinberg case
%\end{abstract}

\maketitle
\tableofcontents

\part{Preliminaries}
\section{Introduction}

\subsection{Background}
 The graded Hecke algebra was introduced by Lusztig \cite{lusztig1989affine} in the study of the representation theory of reductive $p$-adic groups. The Iwahori component, namely the indecomposable block containing the trivial representation, can be much reduced to the module categories of corresponding graded Hecke algebras.

The Ext-group is a fundamental invariant for studying the homological properties of modules.
The Ext-groups for $p$-adic reductive groups and related settings have fruitful studies and applications in the literature:
\begin{itemize}
\item \cite{casselmanvanishing} shows the Ext-vanishing between the Steinberg representation and the trivial representation;
\item \cite{orlik2005extensions} computes the Ext-groups between two generalized Steinberg representations of $p$-adic reductive groups;
\item \cite{opdamsolleveldhomo09} shows Ext-vanishing between discrete series for affine Hecke algebras. This has applications to a classification of discrete series for these algebras \cite{opdamsolleveld10}.
\item \cite{opdamsolleveldExttemp} computed extensions between tempered representations in terms of $R$-groups.
\item \cite{chan2016duality} shows the Ext-vanishing between discrete series for graded Hecke algebras with arbitrary parameters, and has an application to the classification of discrete series of the graded Hecke algebra of type $H_4$ \cite{chanhuang}.
\item \cite{chanqarciinductiveExt} showed an inductive theorem on deducing some non-vanishing Ext-groups.
\end{itemize}
Other related developments include Dirac cohomology \cite{barbaschciubotrapa} and homological branching law \cite{PrasadExt}.

This article aims to compute Ext-groups between certain Arthur type modules for the graded Hecke algebra of type A, extending some of the Ext-groups discussed above. Arthur-type \(\mathbb H_n\)-modules form a large class of unitary modules, generalizing discrete series and tempered modules. A main tool is the Koszul resolution constructed in \cite{chan2016duality}. We notice that there are some relevant computations in a recent preprint \cite{chen2026extension}.

\subsection{Main results}

To explain our results, we need to introduce notation. Let $\mathbb H_n$ be the graded Hecke algebra of type $A_{n-1}$ in Definition \ref{def graded hecke type A} specialized at $c=-1$. Let $\lambda$ be a partition of $n$, and let $\pi_{\lambda}$ be its associated $S_n$-representation. For each $z \in \mathbb C^{\times}$, we lift $\pi_{\lambda}$ to an $\mathbb H_n$-module, denoted by $\tau^*_z(\lambda)$, with the action given in Proposition \ref{prop w module jucys murphy element}. We shall refer to these modules as $W$-type modules, following  similar terminology in \cite{babasical1999Ktype}. We emphasize that the modules in \cite{babasical1999Ktype} are unitary, but the modules we consider here are not necessarily unitary. This modification is also needed for the study of Ext-groups of some Arthur type representations.

A $W$-type module $\tau^*_z(\lambda)$ is said to be Speh if $\lambda$ is rectangular, that is all parts in the partition are the same. These modules correspond to essentially Speh representations of $p$-adic general linear groups, and this explains the notion.

When $\lambda$ is rectangular, we further say that $\tau^*_z(\lambda)$ is a unitary Speh module if $\tau^*_z(\lambda)$ is unitary in the sense of \cite[Section 2.5]{barbaschciubotrapa}. We shall not need the precise notion of unitary modules in our context, and we shall rephrase a convenient and equivalent combinatorial condition in Section \ref{ss arthur type modules}. In particular, the Arthur type modules are unitary.

For an $\mathbb H_{n_1}$-module $\pi_1$ and an $\mathbb H_{n_2}$-module $\pi_2$, define $\pi_1 \times \pi_2$ to be the induced module defined in (\ref{eqn parabolic induction}). An irreducible $\mathbb H_n$-module $X$ is said to be of Arthur-type if $X\cong X_1\times \ldots X_k$ for some unitary Speh modules $X_1, \ldots, X_k$.

A central conjecture underlying this study is the following Ext non-vanishing prediction \cite[Conjecture 5.1]{qadri2025extbranch}:

\begin{conjecture}[Chan-Prasad] \label{conj non vanishing Ext}
Let $X$ and $Y$ be Arthur type $\mathbb H_n$-modules. Then $\mathrm{Ext}^i_{\mathbb H_n}(X, Y)\neq 0$ if and only if $X$ and $Y$ have the same central character.
\end{conjecture}

Toward computating Ext-groups between some Arthur type representations, we first have the following explicit formula:

\begin{theorem} \label{thm ext formula young diagram} (=Theorems \ref{thm ext in term hom}+ \ref{dim of Hom}) 
Let $X$ be a Speh  $\mathbb H_n$-module and let $X'$ be a $W$-type  $\mathbb H_n$-module. Then
\[ \dim~\mathrm{Ext}^i_{\mathbb H_n}(X, X') =\dim~\mathrm{Hom}_{S_n}(X \otimes \bigwedge^i\mathbb C^n , X'), \]
and
\[ \dim~\mathrm{Ext}^i_{\mathbb H_n}(X', X) = \dim~\mathrm{Hom}_{S_n}(X'\otimes \bigwedge^i \mathbb C^n, X) .
\]
\end{theorem}

As a consequence, we verify Conjecture \ref{conj non vanishing Ext} in some cases:

\begin{theorem} (=Theorem \ref{thm conj speh}) \label{thm speh module hn}
Conjecture \ref{conj non vanishing Ext} holds when at least one of $X$ and $Y$ is a Speh module of $\mathbb H_n$.
\end{theorem}

The general case of Conjecture \ref{conj non vanishing Ext} is more complicated and seems to require new techniques. One complication is that Frobenius reciprocity and the geometric lemma usually lead to computations of Ext-groups between another family of $\mathbb H_n$-modules called calibrated modules. The precise definition of calibrated modules is provided in Definition \ref{def calibrate}. To this end, we propose another conjecture on the smallest degree for non-vanishing Ext-groups between calibrated modules, see Conjecture \ref{conj smallest ladder non vanishing}. We verify this conjecture for several examples, given in Propositions \ref{prop conj smallest ladder 1}, \ref{prp conj speh w type} and \ref{prop conj two row speh}.

One important computation is the following result on non-vanishing Ext-groups, which is proved by a method similar to that used in the proof of Theorem \ref{thm ext formula young diagram}:

\begin{theorem}  \label{thm ext smallest}
Let $X$ and $Y$ be $W$-type $\mathbb H_n$-modules whose underlying Young diagrams are quasi-rectangular. Then there exists an integer $i^*$ such that 
\[   \mathrm{Ext}^{i^*}(X, Y) \neq 0,
\]
and, furthermore, the smallest such integer can be determined.
\end{theorem}

One useful technique in the general case is the following transfer result between Ext-groups:

\begin{theorem} (=Corollary \ref{cor transfer thm}) \label{thm ext exchange}
Let $X$ be a Speh $\mathbb H_n$-module and let $Y$ be a calibrated $\mathbb H_n$-module. Then there exists a $W$-type module $Z$ such that 
\[  \mathrm{Ext}^i_{\mathbb H_n}(X, Y)\cong \mathrm{Ext}^i_{\mathbb H_n}(X, Z)
\]
for all $i$.
\end{theorem}

Combining Theorems \ref{thm ext smallest} and \ref{thm ext exchange}, we find the smallest degree for non-vanishing Ext-groups between a \(W\)-type module and a calibrated module.

\subsection{Strategy for proving Theorems \ref{thm ext formula young diagram} and \ref{thm ext smallest}}

We first explain the proof for Theorem \ref{thm ext formula young diagram}. Let \(X\) be a Speh \(\mathbb H_n\)-module and let \(X'\) be a \(W\)-type \(\mathbb H_n\)-module. A necessary condition for \(\mathrm{Ext}^i_{\mathbb H_n}(X,X')\neq 0 \) for some \(i\) is that \(X\) and \(X'\) have the same central character. By some combinatorial arguments, the central character requirement is equivalent to the condition that the underlying Young diagram of $X'$ is of quasi-rectangular shape.  In other words, the partition of $X'$ has to be of the form $(l^{d}, d^{k}) $ for  some positive integers $l$, $d$ and $k$. This is shown in  Corollary \ref{cor: shape for same central character}.

To compute the Ext-groups in Theorem \ref{thm ext formula young diagram}, we shall exploit the Koszul complex in Section \ref{ss koszul resol}. Via Littlewood-Richardson coefficients, Theorem \ref{dim of Hom} and Corollary \ref{cor: dim of Hom for l<d case}  provide combinatorial descriptions of the complexes in terms of Young diagrams. 

The final step is to prove that all differentials in this complex vanish, that is Theorem \ref{d_i=o}. A proof for this step requires some structural information of \(W\)-type modules in Section \ref{sss orthogonal gt basis}. A more elaborate exposition on this step involving some other technical components will be given in Section \ref{ss proof strategy}.

The proof for Theorem \ref{thm ext smallest} goes along a similar line. See Remark \ref{rm differentials non vanish ex} for further comments on the statement of  Theorem \ref{thm ext smallest}.

\subsection{Strategy for proving Theorem \ref{thm speh module hn}}

By applying a standard duality, it suffices to consider the situation in which \(X\) is an Arthur-type module and \(X'\) is a Speh \(\mathbb H_n\)-module. We write an Arthur-type module \(X\) as:
\[  X \cong X_1 \times \ldots \times X_k
\]
for some Speh \(\mathbb H_{n_i}\)-modules \(X_i\) (\(i=1,\ldots, k\)). By Frobenius reciprocity,
\[  \mathrm{Ext}^j_{\mathbb H_n}(X_1\times \ldots \times X_k, X')\cong \mathrm{Ext}^j_{\mathbb H_{n_1}\otimes \ldots \otimes \mathbb H_{n_k}}(X_1\boxtimes \ldots \boxtimes X_k, X') .
\]

A special product feature of an Arthur type module is that 
\[  X_1 \times \ldots \times X_k \cong X_{\sigma(1)} \times \ldots \times X_{\sigma(k)}
\]
for some permutation $\sigma$. Under a suitable arrangement of the Speh modules, the direct summand in \(X'|_{\mathbb H_{n_1}\otimes \ldots \otimes \mathbb H_{n_k}}\) that has the same central character as \(X_1 \boxtimes \ldots \boxtimes X_k\) is isomorphic to the tensor product of some \(W\)-type modules in quasi-rectangular shape. This is mainly shown in Lemma \ref{lem decompose }. Now, Theorem  \ref{thm speh module hn} will follow from Theorem \ref{thm ext formula young diagram}.

\subsection{Strategy for proving Theorem \ref{thm ext exchange}}

We first state the main trick in proving Theorem \ref{thm ext exchange}. Let \(X\) be a Speh module and let \( \chi\) be a character of \(\mathbb H_1\). If \( \chi \) lies in the set of numbers that determines the central character of \(X\), then 
\( X \times \chi \cong \chi \times X \).

The first step in proving Theorem \ref{thm ext exchange} is to find a \(W\)-type $\mathbb H_{n+k}$ module $\widetilde{Y}$ (for some \(k\)) and some characters $\chi_1, \ldots, \chi_k$ such that the restriction of $\widetilde{Y}$ to $\mathbb H_1^{\otimes k}\otimes \mathbb H_n$ has a direct summand isomorphic to:
\[   \chi_1 \boxtimes \ldots \boxtimes \chi_k \boxtimes Y .
\]
Then we show that \(  \mathrm{Ext}^i_{\mathbb H_n}(\chi_1\times \ldots \times \chi_k\times X, \widetilde{Y}) 
\)
can be used to recover the Ext-group \( \mathrm{Ext}^i_{\mathbb H_n}(X, Y) \). 

The next step is to apply the trick above and Frobenius reciprocity,
\[  \mathrm{Ext}^i_{\mathbb H_{n+k}}(\chi_1\times \ldots \times \chi_k \times X, \widetilde{Y})\cong \mathrm{Ext}^i_{\mathbb H_{n+k}}(X \times \chi_1\times \ldots \times \chi_k, \widetilde{Y}) \cong \mathrm{Ext}^i_{\mathbb H_n\otimes \mathbb H_1^{\otimes k}}(X\boxtimes \chi_1\boxtimes \ldots \boxtimes \chi_k, \widetilde{Y}). 
\]
Then, one shows that a direct summand in the restriction of $\widetilde{Y}$ to \( \mathbb H_n\otimes \mathbb H_1^{\otimes k}\), with the same central character as \(\mathbb X\boxtimes \chi_1\boxtimes \ldots \boxtimes \chi_k\), takes the form \( Z\boxtimes \chi_1\boxtimes \ldots \boxtimes \chi_k \) for some $W$-type module \(Z\).

Finally, we show that \(\mathrm{Ext}^i_{\mathbb H_n}(X, Z)\) can be computed from \(\mathrm{Ext}^i_{\mathbb H_n}(X, Y)\) and coincides with \( \mathrm{Ext}^i_{\mathbb H_n}(X, Y) \).

\subsection{Organization}
Section \ref{Sec: Hecke alg} reviews some basic notations of graded Hecke algebras of type A and their modules, including central characters, Koszul-type resolutions and  $\bullet$-dual.
Section \ref{Sec preliminaries for sym gp} reviews two important tools in symmetric groups: Littlewood-Richardson coefficients and the Okounkov-Vershik approach. Section \ref{Sec:W-type} introduces the notation of $W$-type modules of graded Hecke algebras and shows when a $W$-type module has same central character with a Speh module. In Section \ref{Sec:comb of Hom}, we use Littlewood-Richardson coefficients to give a combinatorial formula for the dimension of $\mathrm{Hom}_{S_n}(X\otimes\bigwedge^{i}\mathbb{C}^{n}, Y)$, where $X$ is a Speh module and $Y$ is a \(W\)-type module. In Section \ref{sec vanish of differential}, we prove that the differentials in the Koszul resolution vanish when $X$ is a Speh module and $Y$ is a \(W\)-type module. Section \ref{sec first non-vanish Ext for two quasi-rec} investigate the first non-vanishing Ext-groups for two quasi-ractangular type modules. Section \ref{sec calibrated mod and skew young diagram} studies calibrated modules and their relation to skew Young daigrams. Then we introduce the notation of uniform skew Young diagrams in Section \ref{sec uniform skew diagram}. In Section \ref{sec ext non-vanishing conj for Arthur}, we check Conjecture \ref{conj non vanishing Ext} in the case of one of the modules is unitary Speh. In Section 
\ref{sec elementary transform}, we introduce an important class operations called elementary transforms. Then in Section \ref{sec uniform calibrated modules} we introduce the notation of uniform calibrated modules and show the relation between the uniformity and the central character of a calibrated modules. In Section \ref{sec ext between speh and calibrated}, we formulate and verify Conjecture \ref{conj smallest ladder non vanishing} in severtal cases. 

\subsection{Acknowledgment}
 This project is supported in part by the Research Grants Council of the Hong Kong Special Administrative Region, China (Project No: 17305223, 17308324) and the National Natural Science Foundation of China (Project No. 12322120). 

The present article used AI for correcting the grammar, searching for some references, and simplifying calculations in Section \ref{sec vanish of differential}. 

 After the main body of the article had been written up, Quanyong Chen showed an AI-generated proof of Conjecture \ref{conj non vanishing Ext} to the first-named author. The approach of the AI-generated proof seems to be quite different from the treatment in this article. We heavily use the Koszul resolution of graded Hecke algebras, and the article also contains some other results not presented in the AI-generated article. The present article does not use any form of results or ideas from the AI-generated article.

\section{Preliminaries for graded Hecke algebras}\label{Sec: Hecke alg}
\subsection{Graded Hecke Algebras}

\begin{definition} \label{def graded hecke type A}
   The graded Hecke algebra $\mathbb{H}_{n}(c)$ with equal parameter $c \in \mathbb C$ of type $A_{n-1}$ is generated by $\{s_{i}|1\leq i\leq n-1\}$ and $\{e_{1},...,e_{n}\}$,  with the relations given as follows:
   \begin{enumerate}
       \item[(1)] $s_{i}^{2}=1;\quad s_{i}s_{i+1}s_{i}=s_{i+1}s_{i}s_{i+1};\quad s_{i}s_{j}=s_{j}s_{i}$ for $j\neq i,i+1$.
       \item[(2)] $e_{i}e_{j}=e_{j}e_{i}$ for all $i$ and $j$.
       \item[(3)] $s_{i}e_{j}=e_{j}s_{i}$ for $j\neq i,i+1$;
       \item[(4)] $s_{i}e_{i}-e_{i+1}s_{i}=c$ for all $i$. 
   \end{enumerate}
\end{definition}

Note that \(\langle s_1,\ldots, s_{n-1} \rangle \) generate the group algebra \(\mathbb C[S_n]\) of \(S_n\). We shall identify \(s_i\) with the element permuting \(i\) and \(i+1\) in  \(S_n\). For \(1\leq i \neq j\leq n\), define \(s_{i,j}\) to be the involution in \(S_n\) that swaps the numbers \(i\) and \(j\). In particular, we have \(s_i=s_{i,i+1}\) for \(i=1,\ldots, n-1\).

\begin{lemma}
Let $c_1, c_2\in \mathbb C\setminus \left\{ 0 \right\}$. Then $\mathbb{H}_{n}(c_{1})$ and $\mathbb{H}_{n}(c_{2})$ are isomorphic algebras.
\end{lemma} 
\begin{proof}
Define $F: \mathbb H(c_1)\rightarrow \mathbb H(c_2)$  by:
\[   \frac{1}{c_1}e_i \mapsto \frac{1}{c_2}e_i, \quad \mbox{for $i\in \left\{1, \ldots, n\right\}$} ,\quad    s_i \mapsto s_i, \quad \mbox{for $i \in \left\{1, \ldots, n-1\right\}$}.
\]
It is straightforward to check that $F$ determines an isomorphism.
\end{proof}

Due to the above lemma, we will use $\mathbb{H}_{n}(-1)$ (or, for simplicity, $\mathbb{H}_{n}$) to denote graded Hecke algebras of type $A_{n-1}$ with equal parameters in this paper.

%\subsection{Extension Groups and Frobenius reciprocity}

%Let $\mathbb H$ be a graded Hecke algebra. Let $\pi$ be an $\mathbb H$-module. Let 
%\[  \ldots P_2 \stackrel{d_1}{\rightarrow} \rightarrow P_1   \stackrel{d_0}{\rightarrow} P_0 \rightarrow \pi \rightarrow 0 .
%\]
%be a projective resolution of $\mathbb H$-modules on $\pi$. For another $\mathbb H$-module $\pi'$, denote by $d_i^*: \mathrm{Hom}_{\mathbb H}(P_i, \pi') \rightarrow \mathrm{Hom}_{\mathbb H}(P_{i_1}, \pi')$ the map induced from $d_i$. Define the $i$-th extension group $\mathrm{Ext}^i_{\mathbb H}(\pi, \pi')$ between $\pi$ and $\pi'$ to be the $i$-th homology of the complex 
%\[  \left\{ \mathrm{Hom}_{\mathbb H}( P_i, \pi') , d_i^* \right\} .\]

\subsection{Central Characters}

For a finite-dimensional $\mathbb H$-module $X$, $\chi \in \mathbb{C}^{n,\vee}$ is said to be a \emph{weight} of $X$ if there exists a non-zero vector $x$ in $X$ such that $(v-\chi(v)).x=0$ for all $v \in \mathbb{C}^{n}$.

\begin{proposition}\cite[Proposition 4.5]{lusztig1989affine}
    The center of the graded Hecke algebra $\mathbb{H}_n$ is $Z(\mathbb H_n)=\mathbb C[e_1,\ldots, e_n]^{S_n}$, the ring of $S_n$-invariant polynomials in $e_1, \ldots, e_n$.
\end{proposition}

The group $S_{n}$ has a natural action on $\mathbb{C}^{n,\vee}$:
\begin{equation*}
    (w.\chi)(v)=\chi(w^{-1}.v)
\end{equation*}
for $w\in S_{n}$, $\chi\in \mathbb{C}^{n,\vee}$ and $v\in \mathbb{C}^{n}$. Let $W\chi$ be the orbit of $\chi$ in $\mathbb{C}^{n,\vee}$. An $\mathbb{H}_{n}$-module $X$ is
said to have a $\emph{central character}$ $W\chi$ if there exists a function
$\chi\colon Z(\mathbb H_{n})\to\mathbb C$ such that every $z\in Z(\mathbb H)$ acts
by the scalar $\chi(z)$ on $X$. For simplicity, we also call $\chi$ the central character of $X$. 

For $\mathbb H_{n}$-modules $X$ and $Y$, let $\operatorname{Ext}_{\mathbb H_{n}}^i(X,Y)$
be the $i$-th extension group of $X$ and $Y$ in the category of $\mathbb H_{n}$-modules. We have the following standard and useful result on a relation between Ext-groups and central characters: 

\begin{proposition} \label{prop cc vanish}
    Let $X$ and $Y$ be $\mathbb H_{n}$-modules. If $X$ and $Y$ have different central characters, then 
    \[ \operatorname{Ext}_{\mathbb H_{n}}^i(X,Y)=0 \] for all $i$.
\end{proposition}

\subsection{Koszul-type Resolution for $\mathbb{H}_{n}$-modules} \label{ss koszul resol}
In this section, we introduce a Koszul-type resolution for $\mathbb{H}_{n}$-modules, and use it to introduce a derived complex, which is the main tool for computing the Ext-groups in this paper.

For $v\in \mathbb{C}^{n}$, we define an element in $\mathbb{H}_{n}$ as follows:
\begin{equation}\label{eq:v-tilde}
  \widetilde{v}
  =v-\frac12\sum_{1\leq i<j\leq n}
  (-1)\langle v,e_{i}-e_{j}\rangle s_{i,j}.
\end{equation}

Let $X$ be an $\mathbb{H}_{n}$-module. Define a sequence of $\mathbb{H}_{n}$-maps $\widetilde d_i$ as follows:
\begin{equation}\label{eq:resolution}
\begin{aligned}
0 &\longrightarrow
\mathbb{H}_{n}\otimes_{\CC[S_n]}\bigl(\Res_{S_{n}} X\otimes\textstyle\bigwedge^n \mathbb{C}^{n}\bigr)
\xrightarrow{\widetilde d_{n-1}}\cdots\xrightarrow{\widetilde d_i}
\mathbb{H}_{n}\otimes_{\CC[S_n]}\bigl(\Res_{S_n} X\otimes\textstyle\bigwedge^i \mathbb{C}^{n}\bigr)
\\
&\xrightarrow{\widetilde d_{i-1}}\cdots\xrightarrow{\widetilde d_0}
\mathbb{H}_{n}\otimes_{\CC[S_n]}\Res_{S_n} X
\xrightarrow{\epsilon}X\longrightarrow0,
\end{aligned}
\end{equation}
where $\epsilon\colon\mathbb{H}_{n}\otimes_{\CC[S_n]}\Res_{S_n} X\to X$ is given by
\begin{equation*}
   \epsilon(h\otimes x)=h.x, 
\end{equation*}
and, for $i\geq1$,
\begin{equation*}
  \widetilde d_i\colon
  \mathbb{H}_{n}\otimes_{\CC[S_n]}\bigl(\Res_{S_n} X\otimes\textstyle\bigwedge^{i+1}\mathbb{C}^{n}\bigr)
  \longrightarrow
  \mathbb{H}_{n}\otimes_{\CC[S_n]}\bigl(\Res_{S_n} X\otimes\textstyle\bigwedge^{i}\mathbb{C}^{n}\bigr)
\end{equation*}
is given by
\begin{align*}
&\widetilde d_i\bigl(h\otimes(x\otimes v_1\wedge\cdots\wedge v_{i+1})\bigr)
\\
&\quad=
\sum_{j=1}^{i+1}(-1)^{j+1}
\Bigl(
  h\widetilde v_j\otimes x\otimes
    v_1\wedge\cdots\wedge\widehat{v_j}\wedge\cdots\wedge v_{i+1}
\\[-2pt]
&\hspace{47mm}{}
  -h\otimes \widetilde v_j. x\otimes
    v_1\wedge\cdots\wedge\widehat{v_j}\wedge\cdots\wedge v_{i+1}
\Bigr).
\end{align*}

\begin{theorem}\cite[Theorem 3.4]{chan2016duality}
    For any $\mathbb{H}_{n}$-module $X$, the sequence \eqref{eq:resolution} is a complex and forms a projective resolution for $X$.
\end{theorem}

\begin{remark}
    In \cite{chan2016duality}, Chan introduced two equivalent Koszul type resolutions for $\mathbb{H}_{n}$-modules. The complex \eqref{eq:resolution} is called the alternate form of the Koszul-type resolution and has advantages for computations.
\end{remark}

Now we use the resolution \eqref{eq:resolution} to construct a complex for computing 
extension groups. Let $X$ and $Y$ be $\mathbb{H}_{n}$-modules. Applying the functor $\Hom_{\mathbb{H}_{n}}(\,\cdot\,,Y)$  to the resolution of $X$ as in \eqref{eq:resolution}, we get the induced complex:
\begin{equation}\label{eq:ext-complex}
\begin{aligned}
0&\longleftarrow
\Hom_{S_n}(\Res_{S_n} X\otimes\textstyle\bigwedge^n\mathbb{C}^{n},\Res_{S_{n}} Y)
\xleftarrow{d_{n-1}^\vee}\cdots
\xleftarrow{d_{i}^\vee}
\Hom_{S_n}(\Res_{S_n} X\otimes\textstyle\bigwedge^i\mathbb{C}^{n},\Res_{S_{n}} Y)
\\
&\xleftarrow{d_{i-1}^\vee}\cdots
\xleftarrow{d_{0}^\vee}
\Hom_{S_n}(\Res_{S_n} X,\Res_{S_n} Y)\longleftarrow0,
\end{aligned}
\end{equation}
where the differential $d_{i}^\vee$ can be explicitly written as
\begin{align*}
d_{i}^\vee(\eta)
\bigl(x\otimes v_1\wedge\cdots\wedge v_{i+1}\bigr)
&=
\sum_{j=1}^{i+1}(-1)^{j+1}\widetilde v_j\cdot
\eta\bigl(x\otimes
v_1\wedge\cdots\wedge\widehat{v_j}\wedge\cdots\wedge v_{i+1}\bigr)
\\
&\quad-
\sum_{j=1}^{i+1}(-1)^{j+1}
\eta\bigl(\widetilde v_j\cdot x\otimes
v_1\wedge\cdots\wedge\widehat{v_j}\wedge\cdots\wedge v_{i+1}\bigr).
\end{align*}
Here the action of $\widetilde v_j$ on the term
$\eta(x\otimes v_1\wedge\cdots\wedge\widehat{v_j}\wedge\cdots\wedge
v_{i+1})$ is via the action of $\widetilde v_j$ on $Y$, and the action of $\widetilde v_j$ on the term $x$ is via the action of $\widetilde v_j$ on $X$.

\begin{proposition}\cite[Proposition 3.7]{chan2016duality}
Let $X$ and $Y$ be $\mathbb{H}_{n}$-modules. Then $\Ext_{\mathbb{H}_{n}}^i(X,Y)$ is
isomorphic to the $i$-th homology of the complex in
\eqref{eq:ext-complex}.
\end{proposition}

\subsection{Hecke subalgebras of $\mathbb H_n$}

For a sequence of positive integers $\mathbf n=(n_1, \ldots, n_k)$ with $n_1+\ldots + n_k=n$, we shall usually regard $\mathbb H_{n_1}\otimes \ldots \otimes \mathbb H_{n_k}$ as the subalgebra of $\mathbb H_n$ via the following embedding: for $t \in \left\{ 1, \ldots k \right\}$, and for $j \in \left\{ 1, \ldots, n_t \right\}$,
\[  1^{\otimes (t-1)} \otimes e_j \otimes 1^{\otimes (k-t)} \mapsto e_{n_1+\ldots+n_{t-1}+j}
\]
and, for $j \in \left\{ 1, \ldots, n_t \right\}$,
\[ 1^{\otimes (t-1)} \otimes s_j \otimes 1^{\otimes (k-t)} \mapsto s_{n_1+\ldots +n_{t-1}+j} .
\]

For $i=1, \ldots, k$, let $X_i$ be an $\mathbb H_{n_i}$-module. Define an \(\mathbb H_n\)-module as follows,
\begin{equation} \label{eqn parabolic induction}
X_1\times \ldots \times X_k := \mathbb H_{n} \otimes_{\mathbb H_{n_1}\otimes \ldots \otimes \mathbb H_{n_k}}(X_1\boxtimes \ldots \boxtimes X_k) .
\end{equation}

\begin{lemma} \label{lem frobenius}
Let $\sigma$ be an $\mathbb H_{n_1}\otimes \ldots \otimes \mathbb H_{n_k}$-module and let $\pi$ be an $\mathbb H_n$-module. Then, for any non-negative integer $i$,
\[  \mathrm{Ext}^i_{\mathbb H_n}(\mathbb H_n\otimes_{\mathbb H_{n_1}\otimes \ldots \otimes \mathbb H_{n_k}}\sigma, \pi ) \cong \mathrm{Ext}^i_{\mathbb H_{n_1}\otimes \ldots \otimes \mathbb H_{n_k}}(\sigma, \pi|_{\mathbb H_{n_1}\otimes \ldots \otimes \mathbb H_{n_k}}) .
\]
\end{lemma}

\begin{proof}
This follows from that $\mathbb H_n\otimes_{\mathbb H_{n_1}\otimes\ldots \otimes \mathbb H_{n_k}}$ is left adjoint to the restriction functor, and the restriction functor is exact.
\end{proof}

\subsection{$\bullet$-dual}

Define $\bullet: \mathbb H_n \rightarrow \mathbb H_n$ to be the \(\mathbb C\)-linear anti-involution determined by:
\[   v^{\bullet} = v , \quad v \in V, \quad \mbox{and} \quad  t_w^{\bullet} =t_{w^{-1}}=t_w^{-1} . \]
For any $\mathbb H_n$-module $\pi$, the $\bullet$-dual, denoted by $X^{\bullet}$, of $X$ is the $\mathbb H_n$-module whose underlying space is the space of linear functions from $X$ to $\mathbb C$ with the $\mathbb H_n$-action given by:
\[  (h.f)(x)=f(h^{\bullet}.x) .\]
One important property of $\bullet$ is that $\bullet$ sends any irreducible $\mathbb H_n$-module to any irreducible $\mathbb H_n$-module.

\begin{lemma} \label{lem irr bullet self dual}
Let \(X\) be an irreducible \(\mathbb H_{n_1}\otimes \cdots \otimes \mathbb H_{n_k}\)-module. Then \(X^{\bullet} \cong X\).
\end{lemma}

\begin{proof}
Since the weight spaces of \(X\) and \(X'\) are the same, the lemma follows from \cite[Propsoiton 4.3.1]{barbasch2013hermitian}.
\end{proof}

By abuse of notation, we shall similarly write $\bullet$ to be the analogous anit-involution for $\mathbb H_{n-i}\otimes \mathbb H_i$ to $\mathbb H_{n-i}\otimes \mathbb H_i$. 

For an $\mathbb H_{n-i}\otimes \mathbb H_i$-module $\pi$, we also similarly define $X^{\bullet}$.

\begin{lemma} \label{lem dual of bullet}
Let $X$ be a finite-dimensional $\mathbb H_n$-module. Then 
\[   (X^{\bullet})|_{\mathbb H_{n-i}\otimes \mathbb H_i} \cong (X|_{\mathbb H_{n-i}\otimes \mathbb H_i})^{\bullet} .
\]

\end{lemma}

\begin{proof}
Note that for any $h \in \mathbb H_{n-i}\otimes \mathbb H_i$, the $\bullet_{\mathbb H_{n-i}\otimes \mathbb H_i}$-operation on $h$ coincides with the $\bullet_{\mathbb H_n}$-operation on $h$. Then the isomorphism follows from the definition of the $\bullet$-dual.
\end{proof}

\begin{proposition} \label{prop dual ext} \cite[Proposition 4.10]{chan2016duality}
Let $X$ and $Y$ be finite dimensional $\mathbb H_n$-modules. Then
\[   \mathrm{Ext}^i_{\mathbb H_n}(X, Y) \cong \mathrm{Ext}^i_{\mathbb H_n}(Y^{\bullet}, X^{\bullet}) .
\]
\end{proposition}

\section{Preliminaries for symmetric groups}\label{Sec preliminaries for sym gp}
\subsection{Some notations for symmetric groups}
For a partition $\lambda=(\lambda_{1}\geq...\geq\lambda_{t}>0)\vdash n$ of $n$, there is an associated Young diagram (also denoted by $\lambda$), and an associated irreducible representation, dentoed by $\pi_{\lambda}$, of $S_{n}$.  If we flip $\lambda$ along the left-top to right-bottom diagonal line, then it is also a Young diagram, denoted by $\lambda^{t}$. 
We denote by $P_{n}$ the set of partitions of $n$, and also the set of Young diagrams with $n$ boxes. If we can get a Young diagram after deleting a rightmost box in a row of a Young diagram $\lambda$, then we call this box to be a \emph{corner} of $\lambda$. A Young diagram $\lambda$ is said to be a hook if $\lambda=(\ell, 1^k)$ for some $\ell$ and $k$.

If two Young diagrams $\mu=(\mu_{1},...,\mu_{t})$ and $\lambda=(\lambda_{1},...,\lambda_{t})$ satisfy $\mu_{i}\leq\lambda_{i}, \forall 1\leq i\leq t$, then we call $\mu$ to be a \emph{subdiagram} of $\lambda$ and denote by $\mu\subseteq\lambda$. In this case, we can obtain a \emph{skew diagram} by removing the boxes of $\mu$ in $\lambda$, which is denoted by $\lambda-\mu$ (In this paper, we always assume that $\mu_1<\lambda_1$). We can put numbers into a skew diagram $\lambda-\mu$, if the entries consist of  $\xi_{1}$'s $1$, \ldots , $\xi_{k}$'s $k$, then we say that we put numbers via $\xi=(\xi_{1},...,\xi_{k})$.
The sequence $\textnormal{word}(\lambda-\mu)$ of positive integers obtained by reading the entries of the skew tableau from top row to bottom row and right column to left column is called the \emph{word} of the skew tableau. 

For a rectangular Young diagram $\lambda=((\lambda_1)^t)$ and a Young subdiagram $\mu=(\mu_1,\ldots, \mu_t)$ of $\lambda$, define $(\lambda-\mu)^{\circlearrowleft}$ to be the skew Young diagram obtained by rotating $\lambda-\mu$ counterclockwisely by $180^{\circ}$, i.e. 
\begin{equation} \label{eqn clockwise}
(\lambda-\phi)^{\circlearrowleft}=(\lambda_{1}-\mu_{t},\lambda_{1}-\mu_{t-1},...,\lambda_{1}-\mu_{1}) . %-(\lambda_{1}-\lambda_{t},\lambda_{1}-\lambda_{t-1},...,\lambda_{1}-\lambda_{2},0) .
\end{equation}

We sometimes label each box of a Young diagram $\lambda$ with a unique number. We shall write $\boxed{a}$ to represent the box labeled with the unique number $a$. 

Given a Young diagram $\lambda$, we can label the boxes of $\lambda$ with numbers as follows:
\begin{enumerate}
    \item Fill the top-left box with the complex number $z$.
    \item The numbers in the adjacent boxes of the same row increase by one from left to right.
    \item The numbers in the adjacent boxes of same column decrease from up to down.
\end{enumerate}
Given a box $\Box$ in $\lambda$, the $z$-content of $\Box$ is the number filled in that box by the above process.
The multi-set of numbers in $\lambda$ is called the {\it $z$-content} of $\lambda$.  For a skew Young diagram $\lambda-\mu$, we define the \(z\)-content of \(\lambda-\mu\) to be the multiset difference between the \(z\)-content of \(\lambda\) and the \(z\)-content of \(\mu \).

When $z=0$, we sometimes simply call the $z$-content of \(\lambda-\mu\) to be the content of \(\lambda -\mu\), i.e.
if $\Box$ is in the $c$-th column and $r$-th row of the Young diagram (we often denote it by $(r,c)$),then the content of this box is:
\begin{equation}
    \textnormal{cont}_{\lambda}(\Box)=c-r .
\end{equation}
In other words, the content of a box in a Young diagram $\lambda$ is the same as the $0$-content of the box and the $z$-content of $\Box$ is
\[  z+\mathrm{cont}_{\lambda}(\Box).
\]

For a Young diagram $\lambda$, denote by $\mathrm{SYT}(\lambda)$, the set of standard Young tableaux of shape $\lambda$. Let $T \in \mathrm{SYT}(\lambda)$. We also define $\mathrm{cont}_T(\Box)=\mathrm{cont}_{\lambda}(\Box)$.

\begin{example}
The set of numbers filled in the following Young diagram \begin{ytableau}
*(white)  2& *(white)3 & *(white) 4& *(white)5\\
*(white)  1& *(white) 2& *(white)3\\
*(white)  0& *(white) 1& *(white) 2\\
*(white) -1& *(white)0 
\end{ytableau}
represent the $2$-content of the corresponding boxes.
\end{example}

\begin{example}
In the Young tableau
$T=$\begin{ytableau}
*(white)  1& *(white)2 & *(white) 4& *(white)7\\
*(white)  3& *(white) 5& *(white)6\\
*(white)  8& *(white) 10\\
*(white) 9& *(white)11 
\end{ytableau}
, $\textnormal{cont}_{T}(\boxed{6})=3-2=1$, $\textnormal{cont}_{T}(\boxed{9})=1-4=-3$.
\end{example}

\begin{lemma}\label{lem multipicity of content determine a diagram}
Let $\lambda=(\lambda_0\geq\lambda_1\geq...\geq\lambda_{s-1}>0)$ be a Young diagram. Define
\begin{equation}
    m_i(\lambda)=\#\{(r,c)\in\lambda|c-r=i\}.
\end{equation}
 Then the set $\{m_{i}(\lambda)|i\in\mathbb{\mathbb{C}}\}$ determines the shape of $\lambda$.
\end{lemma}
\begin{proof}
    It's not difficult to check that
    \[
    m_i(\lambda)-m_{i-1}(\lambda)=\textbf{1}_{\{-s+1,...,0\}}(i)-\#\{r|\lambda_{r}-r=i\},
    \]
    where $\textbf{1}_{\{-s+1,...,0\}}$ is the characteristic function of the set $\{-s+1,...,0\}$. Since $r\rightarrow\lambda_r-r$ is strictly decreasing, each value appears at most once. Then we can get the minimal content, $s$ and the set $\{\lambda_r-r|r=0,...,s-1\}$ from the set $\{m_{i}(\lambda)|i\in\mathbb{C}\}$. If we rearrange these values from largest to smallest:
    \[
    i_0>i_1\>...>i_{s-1},
    \]
    then we get 
    \[
    \lambda_r=i_r+r.
    \]
\end{proof}

\subsection{Littlewood-Richardson coefficients}\label{Sec:LR coefficients}

\begin{definition}
A \emph{Littlewood-Richardson (LR) tableau} is a skew tableau $\lambda-\mu$ with the following
properties:
\begin{itemize}
    \item[(\romannumeral1)]  It is \emph{semistandard}, i.e., the numbers in each row of $\lambda-\mu$ weakly increase from left-to-right, and the numbers in each column of $\lambda-\mu$ strictly increase from top-to-bottom.
    \item[(\romannumeral2)] It satisfies the \emph{Yamanouchi word condition} (YWC), i.e., for each positive integer $j$, starting from the first entry of $\textnormal{word}(\lambda-\mu)$ to any place in $\textnormal{word}(\lambda-\mu)$, there are at least as many $j$ as $(j+1)$.
\end{itemize}
\end{definition}

\begin{example}
The number-filled part of $T_{1}=$\begin{ytableau}
*(white)  & *(white) & *(white) 1& *(white)1\\
*(white)  & *(white) 1& *(white)2\\
*(white)  2& *(white) 3\\
*(white) 3& *(white)4 
\end{ytableau}
is a LR tableau in the shape of $(4,3,2,2)-(2,1)$ and is number-filled via $(3,2,2,1)$, its word is $\{1,1,2,1,3,2,4,3\}$.
The number-filled part of $T_{2}=$\begin{ytableau}
*(white)  & *(white) & *(white) 1& *(white)2\\
*(white)  & *(white) 1& *(white)2\\
*(white)  2& *(white) 3\\
*(white) 3& *(white)4 
\end{ytableau} is not a LR tableau because its word is $\{2,1,2,1,3,2,4,3\}$, which does not satisfy the Yamanouchi word condition.
\end{example}

\begin{definition}\label{def of LR}
    For Young diagrams $\lambda$, $\nu$ and $\xi$, the \emph{Littlewood-Richardson (LR) coefficient} $\textnormal{LR}^{\lambda}_{\nu,\ \xi}$ is  the number of LR tableaux in the shape of $\lambda-\nu$ which are numbering via $\xi$.
\end{definition}

\begin{example}
Let $\lambda=(4,4,4,4)\vdash 16$, $\nu=(3,2)\vdash 5$, then the unique partition $\xi\vdash11$ which makes $LR^{\lambda}_{\nu,\ \xi}\neq0$ is $(4,4,2,1)$. In this case, $\textnormal{LR}^{\lambda}_{\nu,\xi}=1$ and the unique method to put numbers $\xi$ is as follow: 
\[
\ytableausetup{boxsize=1.2em}
\ytableaushort{
   {} {} {} 1,
   {} {} 1 2,
   1  1 2 3, 
   2 2 3 4
}
\]

\end{example}

It is well-known (see e.g. \cite{fulton1997young} or \cite[Section 7]{vershik2005new}) that  the Littlewood-Richardson coefficient can be used to describe the restriction of irreducible \(S_n\)-representations to \(S_k\times S_{n-k}\)-representations.

\begin{lemma}
Let \(\lambda \vdash n\). 
$$\text{Res}_{S_k \times S_{n-k}}^{S_n}\pi_{\lambda} \cong \bigoplus_{\mu \vdash k, \nu \vdash n-k} (\pi_{\mu} \boxtimes \pi_{\nu})^{\oplus \textnormal{LR}_{\mu, \nu}^\lambda}.$$
\end{lemma}

\begin{lemma}\label{lem: sym of LR}
Let \( \lambda, \nu, \zeta\) be Young diagrams. Then
$\textnormal{LR}^{\lambda}_{\nu,\xi}=\textnormal{LR}^{\lambda}_{\xi,\nu}$,\quad $\textnormal{LR}^{\lambda}_{\nu,\xi}=\textnormal{LR}^{\lambda^{t}}_{\nu^{t},\xi^{t}}$.
\end{lemma}

\begin{proof}
This result is quite standard, and we only breifly outline a proof. Let \( w_k \) be the longest element in the set of minimal representatives of \( S_n/(S_{n-k}\times S_k)\). Then, for \(i=1, \ldots, k\)
\[  w_k(i)=n-k+i , 
\]
and, for \(i=k+1, \ldots, n\),
\[  w_k(i) = i-k . \]
This induces an isomorphism from \(S_k \times S_{n-k}\) to \(S_{n-k}\times S_k\), and induces a natural isomorphism
\[  \mathrm{Res}^{S_n}_{S_k \times S_{n-k}}\pi_{\lambda} \rightarrow \mathrm{Res}^{S_n}_{S_{n-k}\times S_k}\pi_{\lambda} .
\]
From this, one deduces that $\textnormal{LR}^{\lambda}_{\nu,\xi}=\textnormal{LR}^{\lambda}_{\xi,\nu}$.
\end{proof}

\subsection{The Okounkov–Vershik approach} \label{ss ov approach}

The Okounkov--Vershik approach studies the irreducible complex
representations of the symmetric groups through the inductive chain:
\begin{equation*}
    S_{0}\subset S_{1}\subset\cdots\subset S_{n}.
\end{equation*}

The main idea is spectral: the Gelfand–Zetlin algebra is generated by the
Jucys--Murphy elements, whose joint eigenvectors form the
Gelfand--Tsetlin basis.  Their eigenvalues are the contents of the boxes of
standard Young tableaux.  This description leads directly to explicit
formulas for the action of the Coxeter generators; see
\cite{okounkov1996new,vershik2005new,ceccherini2010representation}.

\subsubsection{Branching and the Gelfand--Tsetlin basis} \label{sss gt basis}

Recall that $\pi_\lambda$ is the irreducible $S_{n}$-module indexed by
$\lambda\vdash n$. The notation $\mu\nearrow\lambda$ means that $\lambda$ is
obtained from $\mu$ by adding one box. The branching rule can be described as
\begin{equation}
  \Res^{S_{n}}_{S_{n-1}}\pi_\lambda
  \cong
  \bigoplus_{\substack{\mu\vdash n-1\\ \mu\nearrow\lambda}}\pi_\mu.
  \label{eq:branching}
\end{equation}
Iterating \eqref{eq:branching}
therefore gives a canonical decomposition
\begin{equation}
  \pi_\lambda=\bigoplus_{T\in\textnormal{SYT}(\lambda)}V_T,
  \qquad \dim V_T=1,
  \label{eq:gt-decomposition}
\end{equation}
where a standard Young tableau $T$ is identified with the path
\begin{equation*}
\varnothing=\lambda^{(0)}\nearrow\lambda^{(1)}
\nearrow\cdots\nearrow\lambda^{(n)}=\lambda  
\end{equation*}
obtained by successively adding the boxes labelled with the numbers $1,2\ldots,n$.

Given a standard Young tableau $T$ and $1 \leq i \leq n$, let $T_{\leq i}$ be the Young subtableaux of $T$ containing the first $i$ boxes, i.e. the subdiagram of $T$ consisting of boxes labeled with $1, \ldots, i$. For example, let $T=$
\begin{ytableau}
    1 & 2 & 5 & 7\\
    3 & 4 & 6\\
    8
\end{ytableau}
, then $T_{\leq 1}=$
\begin{ytableau}
    1
\end{ytableau}
, $T_{\leq 5}=$
\begin{ytableau}
    1 & 2 & 5\\
    3 & 4
\end{ytableau}.

The following lemma provides a characterization of the space $V^T$:

\begin{lemma}\label{lem: branching for gt basis}
Let $T \in \mathrm{SYT}(\lambda)$. Let $v \in \pi_{\lambda}$ be in $V_T$.  Let $\mu \vdash i$. Then $v$ is in the $\pi_{\mu}$-isotypic component of $\pi_{\lambda}|_{S_i}$ if and only if $\mathrm{shape}(T_{\leq i})=\mu$.
\end{lemma}

For each $T \in \mathrm{SYT}(\lambda)$, we fix a nonzero vector $v_{T}\in V_{T}$.  The set $\left\{ v_T: T \in \mathrm{SYT}(\lambda) \right\}$ is called a Gelfand--Tsetlin basis for $\pi_{\lambda}$.
The lines $V_{T}$ are canonical, while the vectors are determined only up to
scalars. With an invariant Hermitian inner product, they may be chosen
orthonormally.

\subsubsection{Jucys--Murphy elements} \label{ss JM elements}

Let $Z_{{\color{teal} n}}=Z(\mathbb{C}[S_{n}])$, viewed as a subalgebra of $\mathbb{C}[S_{n}]$. The
Gelfand--Tsetlin algebra is defined as
\begin{equation*}
     \mathcal{GZ}_{n}=\langle Z_{1},\ldots,Z_{n}\rangle.
\end{equation*}
Define the Jucys--Murphy elements by
\begin{equation}
  X_1=0,
  \qquad
  X_k=\sum_{j=1}^{k-1}s_{j,k},
  \quad 2\leq k\leq n.
  \label{eq:yjm}
\end{equation}
They commute pairwise.

\begin{proposition}\cite[Proposition 2.1]{okounkov1996new}
    The Jucys-Murphy elements $X_k$, $1\leq k\leq n$ generate the Gelfand–Zetlin algebra:
    \begin{equation}
\mathcal{GZ}_n=\mathbb{C}[X_1,\ldots,X_n].
  \label{eq:gz-yjm}
\end{equation}
\end{proposition}

For $T\in\textnormal{SYT}(\lambda)$, the defining spectral formula is

\begin{equation} 
  X_{k}.v_{T}=\textnormal{cont}_{T}(\boxed{k })v_{T},
  \qquad 1\leq k\leq n.
  \label{eq:content-eigenvalue}
\end{equation}

Thus the Gelfand--Tsetlin basis is the simultaneous eigenbasis of Jucys–Murphy elements
$X_{1},\ldots,X_{n}$.  Moreover, the content vector
$\bigl(\textnormal{cont}_{T}(\boxed{1}),\ldots,\textnormal{cont}_{T}(\boxed{n})\bigr)$ determines $T$ uniquely.  In the
Okounkov--Vershik construction, the Coxeter relations first determine the possible joint spectra; these spectra are then identified with the content vectors of standard tableaux.  Consequently, the Young graph and the branching rule arise from the internal structure of the group algebras.

Let
\begin{equation*}
    \Omega_{k}=\sum_{1\leq p<q\leq k}s_{p,q}\in Z_{k},
\end{equation*}  
then $X_{k}=\Omega_{k}-\Omega_{k-1}$.  On the irreducible
$S_{k}$-module $V^\nu$, the central element $\Omega_{k}$ acts by the scalar
\begin{equation*}   \sum_{\Box\in\nu}\cont_{T}(\Box)
  =\frac12\sum_{j\geq1}\nu_j(\nu_j-2j+1).
\end{equation*}

Along the path associated with $T$, subtracting the scalars for
$\lambda^{(k)}$ and $\lambda^{(k-1)}$ leaves precisely the content of
the box added at the $k$-th step, proving
\eqref{eq:content-eigenvalue}.  This also makes the spectral and branching
descriptions compatible: the prefix $1,\ldots,k$ records the
$S_{k}$-constituent, while $X_{k}$ records the newly added box.

\subsubsection{Action on the orthogonal Gelfand--Tsetlin basis} \label{sss orthogonal gt basis}

From \eqref{eq:yjm}, one obtains
\begin{equation} 
  X_{i+1}=s_{i}X_{i} s_{i}+s_{i},
  \qquad
  X_{i+1}s_{i}-s_{i}X_{i}=1,
  \qquad
  s_{i}X_{j}=X_{j}s_{i}\quad(j\neq i,i+1).
  \label{eq:local-relations}
\end{equation}
For a standard tableau \(T\), define the \emph{axial distance}:
\begin{equation}
  r_{i}(T)=\textnormal{cont}_{T}(\boxed{i+1})-\textnormal{cont}_{T}(\boxed{i}).
  \label{eq:axial-distance}
\end{equation}
Write $s_{i}.T$ to be the Young tableau obtained by replacing the filling in the box \(\boxed{i} \) with \(\boxed{i+1}\), and replacing the filling in the \(\boxed{i+1}\) with \(\boxed{i}\).

The following theorem is standard, for example, in \cite{vershik2005new}. 

\begin{theorem}\cite[Proposition 4.1]{vershik2005new}\label{thm: Young's orthogonal form}
The unit vectors $v_{T}$ can be chosen compatibly so that
\begin{equation}
  s_{i}v_{T}=\frac{1}{r_i(T)}v_T
  +\sqrt{1-\frac{1}{r_i(T)^2}}\,v_{(s_{i,i+1}T)},
  \label{eq:orthogonal-action}
\end{equation}
where $v_{(s_{i}.T)}=0$ if $s_{i}T$ is not standard.  More explicitly:
\begin{itemize}
  \item if $i$ and $i+1$ are in the same row, then
        $r_{i}(T)=1$ and $s_{i}v_{T}=v_{T}$;
  \item if they are in the same column, then
        $r_{i}(T)=-1$ and $s_{i}v_{T}=-v_{T}$;
  \item otherwise $T'=s_{i}.T$ is standard and, 
        \begin{equation}
          [s_{i}]_{(v_{T},v_{T'})}=
          \begin{pmatrix}
            r_{i}(T)^{-1}&\sqrt{1-r_{i}(T)^{-2}}\\
            \sqrt{1-r_{i}(T)^{-2}}&-r_{i}(T)^{-1}
          \end{pmatrix}.
          \label{eq:orthogonal-block}
        \end{equation}
\end{itemize}
\end{theorem}

\noindent

We provide a brief explanation of Theorem \ref{thm: Young's orthogonal form} for the convenience of the reader. We first write
\begin{equation} \label{eqn w type module}
s_{i}v_{T}=av_{T}+bv_{s_{i}T} 
\end{equation}
for some \(a, b \in \mathbb C^*\). Here we use the last relation in 
\eqref{eq:local-relations} to see that \(v_{T'}\) does not appear in (\ref{eqn w type module}) if \(T'\neq T\) or \(s_i.T\).

Applying $X_{i+1}s_{i}-s_{i}X_{i}=1$ to $v_{T}$ and using
\eqref{eq:content-eigenvalue} gives
\begin{equation*}
    \bigl(\textnormal{cont}_{T}(\boxed{i+1})-\textnormal{cont}_{T}(\boxed{i})\bigr)a=1,
\end{equation*}
so $a=r_{i}(T)^{-1}$. The involution condition for $s_{i}$ then forces
$b^2=1-r_{i}(T)^{-2}$.  This explains \eqref{eq:orthogonal-action}.

\begin{remark}
There is also a rationally normalized seminormal basis $f_{T}$, described in \cite[Page 23]{vershik2005new}. If
$T'=s_{i}T$ is standard and the pair is oriented in increasing Coxeter
length, then
\begin{align}
  s_{i}f_{T}&=r_{i}(T)^{-1}f_{T}+f_{T'},
  &s_{i}f_{T'}&=(1-r_{i}(T)^{-2})f_{T}-r_{i}(T)^{-1}f_{T'}.
  \label{eq:seminormal-action}
\end{align}
Normalizing these vectors gives \eqref{eq:orthogonal-block}.  Hence the
seminormal and orthogonal bases use the same Gelfand--Tsetlin lines, only
their scalings differ.
\end{remark}

\part{Ext-groups between a Speh module and a \(W\)-type module}
\section{$W$-type Modules}\label{Sec:W-type}

In this section, we introduce $W$-type modules of graded Hecke algebras.

\begin{definition}(cf.  \cite{babasical1999Ktype})
    An $\mathbb{H}_{n}$-module $X$ is called a \emph{$W$-type module} if its restriction to $S_n$ is irreducible. Specifically, if the restriction of $X$ to $S_{n}$ is an irreducible representation associated with a rectangular Young diagram, then $X$ is called a {\it Speh module}.
\end{definition}

We rephrase the relation (\ref{eq:local-relations}) as an algebra morphism between \(\mathbb H_n\) and \(\mathbb C[S_n]\).

\begin{proposition}  \label{prop w module jucys murphy element}  (cf.\cite[Lemma 3.2.1]{barbasch2014unitary}\cite[Section 7.1]{chan2019bernstein})
    For every \( z \in \mathbb{C} \), there exists a surjective algebra homomorphism
\[
\tau_z : \mathbb{H}_{n} \to \mathbb{C}[S_n]
\]
given by
\begin{eqnarray*}
    w &\mapsto& w, \quad w \in S_n,\\
    e_k &\mapsto& s_{1,k}+s_{2,k}+...+s_{k-1,k}+z,\\
    e_1 &\mapsto& z.
\end{eqnarray*}
\end{proposition}

Given an irreducible representation of $S_{n}$, we can associate with a $W$-type module of the graded Hecke algebra  $\mathbb{H}_{n}$ of type $A$ as follows.  For every partition $\lambda\vdash n$ and $z\in\mathbb{C}$, define a \(W\)-type module $\tau_{z}^{*}(\lambda)$ to be the \(\mathbb H_n\)-module obtained by pulling back the \(S_n\)-representation \(\pi_{\lambda}\)  via $\tau_{z}$. We shall use the notations in Section \ref{ss ov approach}.

\begin{proposition}\label{prop:central character}\cite{barbasch2014unitary}
    Let \(\lambda\) be a Young diagram. The central character of $\tau_{z}^{*}(\lambda)$ is determined by the $z$-content of the Young diagram $\lambda$.
\end{proposition}
\begin{proof}
     By Proposition \ref{prop w module jucys murphy element} and (\ref{eq:content-eigenvalue}),  for each $T\in\textnormal{SYT}(\lambda)$, each $e_{k}$ acts on $v_{T}$ of $\tau_{z}^{*}(\lambda)$ by multiplying the scalar $\textnormal{cont}_{T}(\boxed{k})+z$. Hence, the central character of \(\tau^*_{z}(\lambda)\) is the $S_{n}$-orbit of 
     \[  (\textnormal{cont}_T(\boxed{1})+z, \ldots, \textnormal{cont}_T(\boxed{n})+z) ,
     \]
     which is identified with the $z$-content of $\lambda$.
\end{proof}

\begin{proposition} \label{prop multiplicity one weight}
Let $X$ be a $W$-type  $\mathbb H_n$-module. Then any weight of $X$ appears with multiplicity one. 
\end{proposition}

\begin{proof}
As in the proof of Proposition \ref{prop:central character}, $v_{T}$ is a weight vector with weight 
 \begin{equation*}
     \gamma_{T}=(\textnormal{cont}_{T}(\boxed{1})+z,...,\textnormal{cont}_{T}(\boxed{n})+z) .
 \end{equation*}
Hence, if \(T\) and \(T'\) are two different standard Young tableaux in $\mathrm{SYT}(\lambda)$, then \(v_T\) and \(v_{T'}\) have different weights.

As
\begin{equation*}
\pi=\bigoplus\limits_{T\in\textnormal{SYT}(\lambda)}\mathbb{C}v_{T},
\end{equation*}
the proposition then follows. 
%If $\gamma_{T}=\gamma_{T'}$, then $\textnormal{cont}_{T}(\boxed{k})=\textnormal{cont}_{T'}(\boxed{k})$ for any $k$, hence $T=T'$ by induction for $k$.
\end{proof}

\begin{proposition}\label{shape for same central character}
    Let $\lambda=(l^{d})$ be the rectangular Young diagram with  $l$ columns and $d$ rows. Let $\mu$ be a skew Young diagram. We label the boxes of $\lambda$ with its $z$-content.
    \begin{enumerate}
        \item If $l\geq d$, then the $z'$-content of $\mu$ is equal to the $z$-content of $\lambda$ as a multiset for some \(z' \in \mathbb C\) if and only if there exists a unique integer $k \geq 0$ and a partition $\alpha=(\alpha_{1},...,\alpha_{k})$ such that $\alpha_{1}\leq l-k-d$  and $\mu=((l-k)^{d}, d+\alpha_{1},...,d+\alpha_{k})\backslash(\alpha_{1},...,\alpha_{k})$.
        \item If $l<d$, the $z'$-content of $\mu$ is equal to the $z$-content of $\lambda$ as a multiset for some \(z' \in \mathbb C\) if and only if there exists a unique integer $k \geq 0$ and a partition $\alpha=(\alpha_{1},...,\alpha_{d-l-k})$ such that  $\alpha_{1}\leq k$ and $\mu=((l+k)^{l}, l+\alpha_{1},...,l+\alpha_{d-l-k})\backslash(\alpha_{1},...,\alpha_{k})$.
    \end{enumerate}
   
\end{proposition}
\begin{proof}
    Firstly, we prove the statement in $l\geq d$ case.  For convenience, we can assume $\lambda$ is number-filled with $0$-content.
     
    The "if" direction is easy to check: The $k$-content of $\mu=((l-k)^{d}, d+\alpha_{1},...,d+\alpha_{k})/(\alpha_{1},...,\alpha_{k})$ is same as the $0$-content of $\lambda$.

    Then we explain for the "only if" direction. We assume taht the $k$-content of skew Young diagram $\mu=\nu/\alpha$  is same as the $0$-contents of $\lambda$.
    
    For convenience, we denote
    \[
     N_j=\#\{(r,c)\in\lambda|c-r=j\},
    \]
    and
    \[
    C_{k}=((l-k)^d,d^k).
    \]
    Then it's easy to check that 
    \[
    m_i(C_k)=N_{i+k},
    \]
    where $m_i(C_k)$ is defined as in Lemma \ref{lem multipicity of content determine a diagram}. Hence a skew Young diagram $\mu=\nu/\alpha$ has $k$-content such that the multiplicity set of the content is same as the $0$-content of the rectangular $\lambda$ if and only if
    \[
    m_i(\nu/\alpha)=m_{i}(C_k).
    \]
    Notice that 
    \[
    m_i(\nu)-m_i(\alpha)=m_i(\nu/\alpha),
    \]
    then the above condition becomes
    \[
    m_i(\nu)=m_{i}(C_k)+m_i(\alpha).
    \]
    Now let us perform a row-by-row calculation. If we add $\alpha_{j}$ boxes to the $d+j$ row of $C_k$, then the set of the contents of these newly added boxes is exactly equal to the set of contents of the $j$-th row of $C_k$, hence
    \[
    m_i((l-k)^d,d+\alpha_1,...,d+\alpha_k)=m_{i}(C_k)+m_i(\alpha).
    \]
    By Lemma \ref{lem multipicity of content determine a diagram}, the shape of $\nu$ is uniquely determined by $\{m_{i}(C_k)+m_i(\beta)|i\in\mathbb{C}\}$, therefore
    \[
    \nu=((l-k)^{d}, d+\alpha_{1},...,d+\alpha_{k}).
    \]
     Then 
     \[
     \mu=\nu/\alpha=((l-k)^{d}, d+\alpha_{1},...,d+\alpha_{k})\backslash(\alpha_{1},...,\alpha_{k}).
     \]

We now prove the second statement. There is a natural one-to-one correspondence between the $z$-content of a (skew) Young diagram $\nu$ and the $(-z)$-content of $\nu^{t}$ by multiplying the scalar $-1$ to the contents of each box. For example:
\[
\begin{ytableau}
    1 & 2 & 3\\
    0 & 1 & 2\\
    -1 & 0 & 1\\
    -2 & -1 & 0\\
    -3 & -2 & -1\\
    -4 & -3 & -2
\end{ytableau}
\]
has $1$-content $\{1,2,3,0,1,2,-1,0,1,-2,-1,0,-3,-2,-1,-4,-3,-2\}$, 
\[ \begin{ytableau}
    -1 & 0 & 1 & 2 & 3 & 4\\
    -2 & -1 & 0 & 1 & 2 & 3\\
    -3 & -2 & -1 & 0 & 1 & 2
\end{ytableau}
\]
has $-1$-content $\{-1,-2,-3,0,-1,-2,1,0,-1,2,1,0,3,2,1,4,3,2\}$, which is natural one-to-one correspondence to the previous one by multiplying the scalar $-1$ to the content of each box. The second statement then follows from the first statement with the above correspondence.
\end{proof}

\begin{corollary}\label{cor: shape for same central character}
    Let $X=\tau_{0}^{*}(l^{d})$ be a Speh $\mathbb{H}_{n}$-module. Suppose that the $W$-type $\mathbb H_n$-module $Y$ has the same central character as $X$. Then the following hold:
    \begin{enumerate}
     \item If $l\geq d$, then  $Y=\tau_{k}^{*}\big((l-k)^{d}),d^{k}\big)$ for some $k\geq 0$. 
     \item If $l<d$, then $Y=\tau_{-k}^*((l+k)^l, l^{d-l-k})$ for some $k \geq 0$.
     \end{enumerate}
\end{corollary}

\begin{example}
    \begin{ytableau}
        *(cyan)0 & *(cyan)1 & *(cyan)2 & *(cyan)3 & *(cyan)4 & *(cyan)5\\
        *(cyan)-1 & *(red!40)0 & *(red!40)1 & *(red!40)2 & *(red!40)3 & *(red!40)4\\
        *(cyan)-2 & *(red!40)-1 & *(yellow!40)0 & *(yellow!40)1 & *(yellow!40)2 & *(yellow!40)3
    \end{ytableau}
    ,
    \begin{ytableau}
        *(cyan)1 & *(cyan)2 & *(cyan)3 & *(cyan)4 & *(cyan)5\\
        *(cyan)0 & *(red!40)1 & *(red!40)2 & *(red!40)3 & *(red!40)4 \\
        *(cyan)-1 & *(red!40)0 & *(yellow!40)1 & *(yellow!40)2 & *(yellow!40)3\\
        *(cyan)-2 & *(red!40)-1 & *(yellow!40)0
    \end{ytableau}
    ,
    \begin{ytableau}
         \none & *(cyan)2 & *(cyan)3 & *(cyan)4 & *(cyan)5\\
        *(cyan)0 & *(cyan)1 & *(red!40)2 & *(red!40)3 & *(red!40)4 \\
        *(cyan)-1 & *(red!40)0 & *(red!40)1 & *(yellow!40)2 & *(yellow!40)3\\
        *(cyan)-2 & *(red!40)-1 & *(yellow!40)0 & *(yellow!40)1
    \end{ytableau}
    ,
     \begin{ytableau}
        *(cyan)2 & *(cyan)3 & *(cyan)4 & *(cyan)5\\
        *(cyan)1 & *(red!40)2 & *(red!40)3 & *(red!40)4 \\
        *(cyan)0 & *(red!40)1 & *(yellow!40)2 & *(yellow!40)3\\
        *(cyan)-1 & *(red!40)0 & *(yellow!40)1\\
        *(cyan)-2 & *(red!40)-1 & *(yellow!40)0
    \end{ytableau}
    ,
    \begin{ytableau}
        \none & *(cyan)3 & *(cyan)4 & *(cyan)5\\
        *(cyan)1 & *(cyan)2 & *(red!40)3 & *(red!40)4 \\
        *(cyan)0 & *(red!40)1 & *(red!40)2 & *(yellow!40)3\\
        *(cyan)-1 & *(red!40)0 & *(yellow!40)1 & *(yellow!40)2\\
        *(cyan)-2 & *(red!40)-1 & *(yellow!40)0
    \end{ytableau}
    ,
    \begin{ytableau}
        \none & *(cyan)3 & *(cyan)4 & *(cyan)5\\
        \none & *(cyan)2 & *(red!40)3 & *(red!40)4 \\
        *(cyan)0 & *(cyan)1 & *(red!40)2 & *(yellow!40)3\\
        *(cyan)-1 & *(red!40)0 & *(red!40)1 & *(yellow!40)2\\
        *(cyan)-2 & *(red!40)-1 & *(yellow!40)0 & *(yellow!40)1
    \end{ytableau}
    and
    \begin{ytableau}
          *(cyan)3 & *(cyan)4 & *(cyan)5\\
          *(cyan)2 & *(red!40)3 & *(red!40)4 \\
          *(cyan)1 & *(red!40)2 & *(yellow!40)3\\
          *(cyan)0 & *(red!40)1 & *(yellow!40)2\\
         *(cyan)-1 & *(red!40)0 & *(yellow!40)1\\
         *(cyan)-2 & *(red!40)-1 & *(yellow!40)0
    \end{ytableau}
    have same contents. We can also notice that they hame same amount of (skew) hooks, where the definition of skew hooks is in Definition \ref{def skew hook}.
\end{example}

\section{Combinatorial description for some Hom spaces} \label{Sec:comb of Hom}

In this section, we will give a combinatorial description for $\textnormal{Hom}_{S_{n}}(X\otimes\bigwedge^{i}\mathbb{C}^{n}, Y)$ in the complex \eqref{eq:ext-complex}, where $X$ is a Speh module and $Y$ is a quasi-rectangular $W$-type module. 

%\begin{definition}
%     If a Young diagram is $\big((l-k)^{d}),d^{k}\big)$ for some $k, d, l$, then $((l+k)^{l}, l^{d-l-k})$ is said to be of \emph{quasi-rectangular} shape.
%\end{definition}
 %    Note that Corollary \ref{cor: shape for same central character} is to say that $Y$ is in quasi-rectangular shape with same amounts of hooks as $X$. 

We denote $\pi_{(n-1,1)}$ by $V$ and call it the \emph{standard representation} of $S_{n}$. There are two well-
known results:
\[ \bigwedge^{i}V=\pi_{(n-i,1,...,1)}\] and 
\[\bigwedge^{i}\mathbb{C}^{n}=\bigwedge^{i-1}V\oplus\bigwedge^{i}V.\]
We will also use these results to compute $ \textnormal{dim}\ \textnormal{Hom}_{S_{n}}(\pi_{\lambda}\otimes\bigwedge^{i}V, \pi_{\mu})$, and this will be useful in the context where $V$ is used to define the graded Hecke algebras and resolution \ref{eq:ext-complex} instead of \(\mathbb C^n\).

\begin{theorem}\cite[Theorem 3.5]{hayashi2015decomposition}\label{hayashi theorem} Let \(\lambda, \mu \vdash n\). Let \(0 \leq i \leq n\). Then the following hold:
    \begin{itemize}
    \item[(1)] 
        \begin{equation*}
           \textnormal{dim}\ \textnormal{Hom}_{S_{n}}(\pi_{\lambda}\otimes\bigwedge^{i}\mathbb{C}^{n}, \pi_{\mu})= \sum\limits_{\nu\vdash n-i}\sum\limits_{\xi\vdash i} \textnormal{LR}^{\lambda}_{\nu,\ \xi}\textnormal{LR}^{\mu}_{\nu,\ \xi^{t}}
        \end{equation*}
    \item[(2)] 
        \begin{equation*}
             \textnormal{dim}\ \textnormal{Hom}_{S_{n}}(\pi_{\lambda}\otimes\bigwedge^{i}V, \pi_{\mu})= \sum_{k=0}^{i}(-1)^{i-k} \textnormal{dim}\ \textnormal{Hom}_{S_{n}}(\pi_{\lambda}\otimes\bigwedge^{i}\mathbb{C}^{n}, \pi_{\mu}).
        \end{equation*}
    \end{itemize}
\end{theorem}

\begin{proof}
    Notice that $\bigwedge^{i}\mathbb{C}^{n}\cong\textnormal{Ind}^{S_{n}}_{S_{n-i}\times S_{i}} \textnormal{triv}\boxtimes\textnormal{sgn}$. Then 
    \begin{eqnarray*}
        \textnormal{dim}\ \textnormal{Hom}_{S_{n}}(\pi_{\lambda}\otimes\bigwedge^{i}\mathbb{C}^{n}, \pi_{\mu})&=&\textnormal{dim}\ \textnormal{Hom}_{S_{n}}(\pi_{\lambda}\otimes \textnormal{Ind}^{S_{n}}_{S_{n-i}\times S_{i}} \textnormal{triv}\boxtimes\textnormal{sgn}, \pi_{\mu})\\
        &=& \textnormal{dim}\ \textnormal{Hom}_{S_{n}}(\textnormal{Ind}^{S_{n}}_{S_{n-i}\times S_{i}} \textnormal{Res}^{S_{n}}_{S_{n-i}\times S_{i}}(\pi_{\lambda})\otimes(\textnormal{triv}\boxtimes\textnormal{sgn}), \pi_{\mu})\\
        &=&\textnormal{dim}\ \textnormal{Hom}_{S_{n-i}\times S_{i}}(\textnormal{Res}^{S_{n}}_{S_{n-i}\times S_{i}}(\pi_{\lambda})\otimes(\textnormal{triv}\boxtimes\textnormal{sgn}), \textnormal{Res}^{S_{n}}_{S_{n-i}\times S_{i}}\pi_{\mu})\\
        &=& 
        \textnormal{dim}\ \textnormal{Hom}_{S_{n-i}\times S_{i}}(\bigoplus_{\nu \vdash n-i, \xi \vdash i} (\pi_{\nu} \boxtimes \pi_{{\xi}^{t}})^{\oplus \textnormal{LR}_{\nu, \xi}^\lambda}, \bigoplus_{\nu \vdash n-i, \xi \vdash i} (\pi_{\nu} \boxtimes \pi_{\xi})^{\oplus \textnormal{LR}_{\nu, \xi}^\mu})\\
        &=&\sum\limits_{\nu\vdash n-i}\sum\limits_{\xi\vdash i} \textnormal{LR}^{\lambda}_{\nu,\ \xi}\textnormal{LR}^{\mu}_{\nu,\ \xi^{t}}.
    \end{eqnarray*}
The second statement in the theorem follows from the isomorphism $\bigwedge^{i}\mathbb{C}^{n}\cong\bigwedge^{i-1}V\oplus\bigwedge^{i}V$.
\end{proof}

The following corollary is a direct application of Theorem \ref{hayashi theorem}. We can use this corollary to get the multiplicity directly in degree one and degree two cases from the Young diagram.
\begin{corollary} 
Let \(\lambda, \mu \vdash n\). We have:
\begin{itemize}
    \item[(1)]
    \begin{equation*}
    \textnormal{dim}\ \textnormal{Hom}_{S_{n}}(\pi_{\lambda}\otimes V, \pi_{\mu})= 
    \begin{cases}
        0,\quad & \textnormal{if}\ \#\{\textnormal{boxes}\ \textnormal{of}\ \lambda\cap\mu\}\leq n-2,\\
        1,\quad & \textnormal{if}\  \#\{\textnormal{boxes}\ \textnormal{of}\  \lambda\cap\mu\}=n-1,\\
        \textnormal{if}\ \#\{\textnormal{corners}\ \textnormal{of}\ \lambda\}-1,\quad & \mbox{if \(\lambda=\mu\)}.
    \end{cases}
     \end{equation*}
    \item[(2)] 
    \begin{eqnarray*}
    &&\textnormal{dim}\ \textnormal{Hom}_{S_{n}}(\pi_{\lambda}\otimes\bigwedge^{2}V, \pi_{\mu})\\&&= 
    \begin{cases}
        0,\quad & \ \textnormal{if}\ \#\{\textnormal{boxes}\ \textnormal{of}\  \lambda\cap\mu\}\leq n-3,\\ 
       \  \#\{\textnormal{corners}\ \textnormal{of}\ (\lambda-\lambda\cap\mu)\}+\#\{\textnormal{corners}\ \textnormal{of}\ (\mu-\lambda\cap\mu)\}-2,\quad &  \   \textnormal{if}\ \#\{\textnormal{boxes}\ \textnormal{of}\  \lambda\cap\mu\}=n-2,\\
        \#\{\textnormal{corners}\ \textnormal{of}\ \lambda\cap\mu\}+\#\{\textnormal{corners}\ \textnormal{of}\ \lambda\cup\mu)\}-3,\quad & \  \textnormal{if}\ \#\{\textnormal{boxes}\ \textnormal{of}\  \lambda\cap\mu\}=n-1\\
        (\#\{\textnormal{corners}\ \textnormal{of}\ \lambda\}-1)^{2},\quad & \  \textnormal{if}\  \lambda=\mu.
    \end{cases}
   \end{eqnarray*}
\end{itemize}

\end{corollary}

Now we will focus on some special cases in which $\lambda$ and $\mu$ have rectangular and quasi-rectangular shapes.

\begin{lemma}\label{lem: lr coef for rectangular}
Let $\lambda$ be a Young diagram with rectangular shape. Then
\begin{equation}
    \textnormal{LR}^{\lambda}_{\nu,\ \xi}=
    \begin{cases}
        0,\quad & \textnormal{if}\  \xi\neq(\lambda-\nu)^{\circlearrowleft}\\
        1,\quad &  \textnormal{if}\ \xi=(\lambda-\nu)^{\circlearrowleft}
    \end{cases}
\end{equation}
 for any subdiagram $\nu\subseteq\lambda$ and content $\xi$, where $(\lambda-\phi)^{\circlearrowleft}$ is defined in (\ref{eqn clockwise}). 
\end{lemma}

\begin{proof}
    This follows directly from the definition of Littlewood-Richardson coefficients \eqref{def of LR}.
\end{proof}

A Young diagram $\lambda$ is of \emph{quasi-rectangular shape} if either  \(\lambda\) is rectangular or \(\lambda=(l_{1}^{k_1},k_{1}^{k_2})\) for some positive integers $l$, $k_1$ and $k_2$.

\begin{theorem}\label{dim of Hom}
   Suppose $l\geq d$. Let $\lambda=(l^{d})$ be a rectangular Young diagram with $n$ boxes. For $l-d\geq k\geq0$  let $\mu=((l-k)^{d},d^{k})$ be a quasi-rectangular Young diagram. Then we have
    \begin{enumerate}
        \item 
        \begin{eqnarray*}
        \textnormal{dim Hom}_{S_{n}}(\pi_{\lambda}\otimes\bigwedge^{i}\mathbb{C}^{n},\pi_{\mu})=\begin{cases}
            0,\quad & \textnormal{if}\ i<kd,\\
            1,\quad& \textnormal{if}\ i=kd,\\
            \#\{\gamma\vdash t|\gamma\subseteq\lambda\cap\mu,\ \gamma=\gamma^{t}\},\quad& \textnormal{if}\ i=kd+t
        \end{cases}
        \end{eqnarray*}
        
        \item 
        \begin{eqnarray*}
        \textnormal{dim Hom}_{S_{n}}(\pi_{\lambda}\otimes\bigwedge^{i}V,\pi_{\mu})=\begin{cases}
            0,\quad & \textnormal{if}\ i<kd,\\
            1,\quad& \textnormal{if}\ i=kd,\\
            \#\{\gamma\vdash t|\gamma\subseteq\lambda\cap\mu,\ \gamma=\gamma^{t},\  \gamma\ \textnormal{has no corner}\\
            \quad\textnormal{on the top-left to bottom-right diagonal line}\},\quad& \textnormal{if}\ i=kd+t
        \end{cases}
        \end{eqnarray*}
    \end{enumerate}
\end{theorem}

\begin{proof}
    By Theorem \ref{hayashi theorem}, 
    \begin{equation*}
        \textnormal{dim}\ \textnormal{Hom}_{S_{n}}(\pi_{\lambda}\otimes\bigwedge^{i}\mathbb{C}^{n}, \pi_{\mu})= \sum\limits_{\nu\vdash n-i}\sum\limits_{\xi\vdash i} \textnormal{LR}^{\lambda}_{\nu,\ \xi}\textnormal{LR}^{\mu}_{\nu,\ \xi^{t}}.
    \end{equation*}
    When $i<kd$, $\nu\vdash n-i$ cannot be contained in $\lambda\cap\mu$, and hence $\textnormal{LR}^{\lambda}_{\nu,\ \xi}\textnormal{LR}^{\mu}_{\nu,\ \xi^{t}}=0$ for any $\nu\vdash n-i$ by definition of Littlewood-Richardson coefficients. Therefore, in this case, $ \textnormal{dim}\ \textnormal{Hom}_{S_{n}}(\pi_{\lambda}\otimes\bigwedge^{i}\mathbb{C}^{n}, \pi_{\mu})=0$.

    When $i=kd$, $\textnormal{LR}^{\lambda}_{\nu,\ \xi}$ and $\textnormal{LR}^{\mu}_{\nu,\ \xi^{t}}$ are both nonzero if and only if  $\nu=\lambda\cap\mu=((l-k)^{d})$ and $\xi=(k^{d})$ by Lemma \ref{lem: lr coef for rectangular}, and in this case, they are equal to one. Therefore, $ \textnormal{dim}\ \textnormal{Hom}_{S_{n}}(\pi_{\lambda}\otimes\bigwedge^{i}\mathbb{C}^{n}, \pi_{\mu})=1$.

    When $i=kd+t$, by Lemma \ref{lem: lr coef for rectangular}, if $\nu\subseteq\lambda\cap\mu$ and $\xi=(\lambda-\nu)^{\circlearrowleft}$, then $\textnormal{LR}^{\lambda}_{\nu,\xi}=1$ and otherwise the LR coefficient is zero. We just need to compute $\textnormal{LR}^{\mu}_{\nu,\xi^{t}}$ for such $\lambda$ and $\xi$. We can write $\xi=(k^{d})+\gamma$, which means each row of $\xi$ has $k$ more boxes than $\gamma$ i.e. $\xi_{j}=k+\gamma_{j}$ for each $j$. Then it's easy to see that $\textnormal{LR}^{\mu}_{\nu,\xi^{t}}=1$ if $\gamma=\gamma^{t}$ and $\gamma\subseteq\lambda\cap\mu$ (to make sure $\xi\subseteq\lambda\cap\mu^{t}$), otherwise the LR coefficient is zero. Therefore   $\textnormal{dim}\ \textnormal{Hom}_{S_{n}}(\pi_{\lambda}\otimes\bigwedge^{i}\mathbb{C}^{n}, \pi_{\mu})=\#\{\gamma\vdash t|\gamma\subseteq\lambda\cap\mu,\ \gamma=\gamma^{t}\}$. This proves the first formula.
    
    As for the standard representation case, notice that symmetric Young diagrams with $n$ boxes can be divided into two classes. The first class consists of Young diagrams with no corner on the top-left to bottom-right diagonal line, and the second class consists of Young diagrams which have a corner on the top-left to bottom-right diagonal line. The number of Young diagrams in the second class of $i$ is equal to the number of Young diagrams in the first class of $i-1$: the bijection is given by deleting/adding the corner on the diagonal line. Then the result is a direct corollary of this observation, the first formula and the well-known decomposition $\bigwedge^{i}\mathbb{C}^{n}=\bigwedge^{i-1}V\oplus\bigwedge^{i}V$ by induction. 
\end{proof}

\begin{corollary}\label{cor: dim of Hom for l<d case}
     Suppose \(l <d\). Let $\lambda=(l^{d})$ be a rectangular Young diagram with $n$ boxes. For $d-l\geq k\geq0$ let $\mu=((l+k)^{l},l^{d-l-k})$  be a quasi-rectangular Young diagram. Then we have
    \begin{enumerate}
        \item 
        \begin{eqnarray*}
        \textnormal{dim Hom}_{S_{n}}(\pi_{\lambda}\otimes\bigwedge^{i}\mathbb{C}^{n},\pi_{\mu})=\begin{cases}
            0,\quad & \textnormal{if}\ i<kl,\\
            1,\quad& \textnormal{if}\ i=kl,\\
            \#\{\gamma\vdash t|\gamma\subseteq\lambda\cap\mu,\ \gamma=\gamma^{t}\},\quad& \textnormal{if}\ i=kl+t
        \end{cases}
        \end{eqnarray*}
        
        \item 
        \begin{eqnarray*}
        \textnormal{dim Hom}_{S_{n}}(\pi_{\lambda}\otimes\bigwedge^{i}V,\pi_{\mu})=\begin{cases}
            0,\quad & \textnormal{if}\ i<kl,\\
            1,\quad& \textnormal{if}\ i=kl,\\
            \#\{\gamma\vdash t|\gamma\subseteq\lambda\cap\mu,\ \gamma=\gamma^{t},\  \gamma\ \textnormal{has no corner}\\
            \quad\textnormal{on the top-left to bottom-right diagonal line}\},\quad& \textnormal{if}\ i=kl+t
        \end{cases}
        \end{eqnarray*}
    \end{enumerate}
\end{corollary}
\begin{proof}
    This is a direct consequence of Theorem \ref{dim of Hom} and Lemma \ref{lem: sym of LR}.
\end{proof}

\begin{example}
 Let   $\lambda=\mu=$
$
\ytableausetup{boxsize=1.2em} % 设置格子大小
\ytableaushort{
  {} {} {} ,
  {} {} {}
}$. In this case, $l=3$, $d=2$, $k=0$. By Theorem \eqref{dim of Hom}, we have
\begin{equation*}
      \textnormal{dim}\ \textnormal{Hom}_{S_{n}}(\pi_{\lambda}\otimes\bigwedge^{i}\mathbb{C}^{n}, \pi_{\lambda})=\begin{cases}
    0,\quad & i=2,\ 5,\ 6,\\
    1,\quad & i=0,\ 1,\ 3,\ 4,
    \end{cases}
\end{equation*}
  \begin{equation*}
      \textnormal{dim}\ \textnormal{Hom}_{S_{n}}(\pi_{\lambda}\otimes\bigwedge^{i}V, \pi_{\lambda})=\begin{cases}
    0,\quad & i=1,\ 2,\ 4,\ 5,\\
    1,\quad & i=0,\ 3.
\end{cases}
\end{equation*}

\end{example}

\begin{example}\label{eg: dim of Hom of rec and quasi rec}
Let $\lambda=$
\begin{ytableau}
    *(white)  & *(white)  & *(white)  & *(white)  & *(white)  &   &  \\
    *(white)  & *(white)  & *(white)  & *(white)  & *(white)  &   &   \\
    *(white)  & *(white)  & *(white)  & *(white)  & *(white)  &   & 
\end{ytableau}
and $\mu=$
\begin{ytableau}
    *(white)  & *(white)  & *(white)  & *(white)  & *(white) \\
    *(white)  & *(white)  & *(white)  & *(white)  & *(white) \\
    *(white)  & *(white)  & *(white)  & *(white)  & *(white) \\
     & & \\
     & &
\end{ytableau}
. In this case, $l=7$, $d=3$, $k=2$. The intersection part \(\lambda\cap \mu=)5^3) \) contains the following self-conjugate subdiagrams: 
\[
\begin{ytableau}
    *(white)
\end{ytableau},
\begin{ytableau}
    *(white) & *(white)\\
    *(white)
\end{ytableau}, 
\begin{ytableau}
    *(white) & *(white)\\
    *(white) & *(white)
\end{ytableau}, 
\begin{ytableau}
    *(white) & *(white) &*(white) \\
    *(white)\\
    *(white)
\end{ytableau}, 
\begin{ytableau}
    *(white) & *(white) &*(white) \\
    *(white)&*(white)\\
    *(white)
\end{ytableau}, 
\begin{ytableau}
    *(white) & *(white) &*(white) \\
    *(white)&*(white)&*(white) \\
    *(white)& 
\end{ytableau}, 
\begin{ytableau}
    *(white) & *(white) &*(white) \\
    *(white)&*(white)&*(white) \\
    *(white)&*(white)&*(white)
\end{ytableau}. 
\]
Among them, the second, fourth, and sixth subdiagrams contain no corner on the top-left to bottom-right diagonal line, while others contain. Therefore, we have 
 \begin{equation*}
        \textnormal{dim Hom}_{S_{n}}(\pi_{\lambda}\otimes\bigwedge^{i}\mathbb{C}^{21},\pi_{\mu})=\begin{cases}
            0,\quad & \textnormal{otherwise,}\\
            1,\quad& i=6,\ 7, \ 9,\ 10,\ \ 11,\ 12,\ 14,\ 15.
        \end{cases}
        \end{equation*}
 \begin{equation*}
        \textnormal{dim Hom}_{S_{n}}(\pi_{\lambda}\otimes\bigwedge^{i}V,\pi_{\mu})=\begin{cases}
            0,\quad & \textnormal{otherwise,}\\
            1,\quad& i=6,\ 9,\ 11,\ 14.
        \end{cases}
        \end{equation*}
        
\end{example}

\section{Vanishing of Differentials}\label{sec vanish of differential}

 In this section, let $\lambda=(l^{d})$ be a rectangular Young diagram with $n$ boxes.  Let $X=\tau_{0}^{*}(\lambda)$.
 \begin{itemize}
 \item  When  $l\geq d$, let $\mu=\big((l-k)^{d}),d^{k}\big)$, for some $k\geq0$, be a quasi-rectangular Young diagram. Let $Y=\tau_{k}^{*}(\mu)$.
 \item When $l<d$, let $\mu=((l+k)^{l},l^{d-l-k})$, for some $k\geq0$, be a quasi-rectangular Young diagram. Let $Y=\tau_{-k}^{*}(\mu)$.
 \end{itemize} 
  Hence $X$ and $Y$ have same central character as in Section \ref{Sec:W-type}. 
If there is no confusion, we sometimes use $X$ and $Y$ to represent $\textnormal{Res}_{S_{n}}X$ and $\textnormal{Res}_{S_{n}}Y$ respectively.

%\subsection{Quasi-symmetric conditions}

 \begin{definition}
 A Young diagram $\nu$ is said to satisfy the \emph{$k$-th column (resp. row) type quasi-symmetric condition} if
 \begin{enumerate}
     \item the first $k$ columns (resp. rows) of \(\nu\) have the same number of boxes,
     \item after removing the first $k$ columns (resp. rows) of $\nu$, the resulting diagram is symmetric (self-conjugate).
 \end{enumerate} 
 \end{definition}

 \begin{example}
    Let $\nu$=
     \begin{ytableau}
       *(white) &  *(white) &  *(white) &  *(white)\\
        *(white) &  *(white) &  *(white)\\
         *(white) &  *(white)
     \end{ytableau}.
  The sub-diagram of $\nu$ obtained after removing the first ($k=1$) column is symmetric. Hence, $\nu$ satisfies the first column type quasi-symmetric condition. However, $\nu$ does not satisfy the first row type quasi-symmetric condition.
 \end{example}

\subsection{Proof strategy} \label{ss proof strategy}

We first explain the main steps in proving Theorem \ref{d_i=o}. We first have the following identification:
\[ \bigwedge^i \mathbb C^n \cong \mathrm{Ind}^{S_n}_{S_i \times S_{n-i}} (\mathrm{sgn} \boxtimes \mathrm{triv}) .
\]
This then gives an identification:
\begin{equation}\label{eq: decomposition of X otimes wedge}
      X\otimes \bigwedge^i \mathbb C^n  \cong \bigoplus_{\phi \vdash i} \mathrm{Ind}_{S_i \times S_{n-i}}^{S_n} (\pi_{\phi^t}\otimes \pi_{(\lambda-\phi)^{\circlearrowleft}}) ,
\end{equation} 
where $\phi$ runs over all subdiagrams of $\lambda$ with $i$ boxes.

Then, by Lemma \ref{lem: lr coef for rectangular} and the following isomorphism
\begin{equation*}
 \mathrm{Ind}^{S_n}_{S_{n-i} \times S_{i}} (\mathrm{triv} \boxtimes \mathrm{sgn}) \cong\bigwedge^i \mathbb C^n \cong  \mathrm{Ind}^{S_n}_{S_i \times S_{n-i}} (\mathrm{sgn} \boxtimes \mathrm{triv}),
\end{equation*}
we obtain
\begin{equation} \label{eqn decompose tensor wedge}
X \otimes \bigwedge^i \mathbb C^n \cong \mathrm{Ind}^{S_n}_{S_i \times S_{n-i}} (\pi_{\lambda}|_{S_i \times S_{n-i}}\otimes (\mathrm{sgn} \boxtimes \mathrm{triv})) \cong \bigoplus_{\phi \vdash i}\mathrm{Ind}_{S_i\times S_{n-i}}^{S_n}( \pi_{\phi^t} \otimes \pi_{(\lambda-\phi)^{\circlearrowleft} }) .
\end{equation}
Consequently, we have the following decomposition of $\mathrm{Hom}_{S_n}(X\otimes \bigwedge^i \mathbb C^n, Y)$.

\begin{proposition}\label{prop: decomposition of Hom in section vanish}
Suppose  $l\geq d$ (resp.  $l< d$).
\begin{equation}\label{eq: decomposition of Hom}
    \mathrm{Hom}_{S_n}(X\otimes \bigwedge^i \mathbb C^n, Y) \cong  \bigoplus_{\phi \vdash i} \mathrm{Hom}_{S_i \times S_{n-i}}(\pi_{\phi^t}\otimes \pi_{(\lambda-\phi)^{\circlearrowleft}}, Y), 
\end{equation}   
where $\phi$ runs over all $k$-th column (resp. row) type quasi-symmetric subdiagrams with $i$ boxes of \(\lambda \cap \mu\).
\end{proposition}

\begin{proof}
The decomposition follows from (\ref{eqn decompose tensor wedge}). The constraint on \(\psi\) follows from Lemma \ref{lem: lr coef for rectangular} and the proof of Theorem \ref{dim of Hom} (resp. Corollary \ref{cor: dim of Hom for l<d case}).
\end{proof}

Via (\ref{eqn decompose tensor wedge}), we obtain a natural projection map
\begin{equation*}
       \operatorname{pr}_{i,\phi}:    \pi_{\lambda}\otimes\bigwedge^{i}\mathbb{C}^{n}
\longrightarrow
\operatorname{Ind}^{S_n}_{S_i\times S_{n-i}}
(\pi_{\phi^{t}}\otimes\pi_{(\lambda-\phi)^{\circlearrowleft}}).
   \end{equation*}

%This induces the projection from $\mathrm{Hom}_{S_n}(X\otimes \bigwedge^i\mathbb C^n, Y)$ to $\mathrm{Hom}_{S_n}(\pi_{\nu^t}\otimes \pi_{(\lambda-\nu)^{\circlearrowleft}}, Y)$.

\begin{lemma}\label{lem: decomposition of f by projection}
Suppose \(l \geq d\) (resp. \(l<d\)). Let $f \in \mathrm{Hom}_{S_n}(X\otimes \bigwedge^i\mathbb C^n, Y)$. Then
    \[
    f=\sum_{\phi\vdash i}f\circ\operatorname{pr}_{i,\phi},
    \]
    where $\phi$ runs over all $k$-th column (resp. row) type quasi-symmetric subdiagrams of \(\lambda \cap \mu\).
\end{lemma}
\begin{proof}
    This follows from Proposition \ref{prop: decomposition of Hom in section vanish} and the definition of $\operatorname{pr}_{i,\phi}$.
\end{proof}

The next step is to show that the differentials $d_i^{\vee}$ in the Koszul complex in \eqref{eq:ext-complex} vanish. By Lemma \ref{lem: decomposition of f by projection}, it is equivalent to show that for each $f \in  \mathrm{Hom}_{S_n}(X\otimes \bigwedge^i \mathbb C^n, Y)$,
\begin{equation} \label{eqn decompose formula}
d_i^{\vee}(f \circ \mathrm{pr}_{i,\phi}) \circ \mathrm{pr}_{i+1, \nu} =0 
\end{equation}
for all \(k\)-th column type quasi-symmetric $\phi \vdash i$ and $\nu \vdash i+1$.  A main step is to apply Proposition \ref{prop w module jucys murphy element} to express,
\[   d_i^{\vee}(f \circ \mathrm{pr_{i,\phi}})(v\otimes e_1\wedge \ldots \wedge e_{i+1}) 
\]
into a nicer formula compatible with the decomposition in (\ref{eq: decomposition of Hom}),

Requiring both $\phi$ and $\nu$ to be quasi-symmetric is a strong condition will allow us to show the vanishing  \eqref{eqn decompose formula} for many pairs of the \(k\)-the column type quasi-symmetric subdiagrams of \(\lambda \cap \mu\) by using the following Lemma \ref{projection}. The precise meaning of these pairs is explained in Case 2 of the proof of Lemma \ref{lem vanishing on gt   basis}.

\begin{lemma}\label{projection}
   %Let $\lambda\vdash n$ and $\phi \vdash i$ with $\phi \subset \lambda$, and let $T$ be a Young tableau associated to $\lambda$. Let 
   %\begin{equation*}
   %    \operatorname{pr}_{i,\phi}:    \pi_{\lambda}\otimes\bigwedge^{i}\mathbb{C}^{n}
   %    \longrightarrow
   %    \operatorname{Ind}^{S_n}_{S_i\times S_{n-i}}
   %   (\pi_{\phi^{t}}\otimes\pi_{(\lambda-\phi)^{\circlearrowleft}})
   %\end{equation*}
   %be the natural projection map.
   Let \(T \in \mathrm{STY}(\lambda)\). Let $v_{T}$ be a vector in the Gelfand-Tsetlin basis for \(X\) defined in Section \ref{sss gt basis}. If $\operatorname{shape}(T_{\leq i})\neq\phi$, then 
   \[
   \operatorname{pr}_{i,\phi}(v_{T}\otimes e_{1}\wedge...\wedge e_{i})=0.
   \]
\end{lemma}

\begin{proof}
This lemma is an easy application of Lemma \ref{lem: branching for gt basis}.
Put
\[
\mu=\operatorname{shape}(T_{\leq i}).
\]
By the multiplicity-free branching rule for
\[
S_n\supset S_{n-1}\supset\cdots\supset S_i .
\]
the restriction of $\pi_{\lambda}$ to $S_{i}$ contains a unique copy of $\pi_{\mu}$ corresponding to the path determined by $T_{\leq i}$, and $v_{T}$ belongs to this $\pi_{\mu}$-summand. 
And $S_{i}$ acts on $e_{1}\wedge e_{2}\wedge...\wedge e_{i}$ by the sign character, hence $v_{T}\otimes e_{1}\wedge e_{2}\wedge...\wedge e_{i}$ belongs to the representation $\pi_{\mu}\otimes \textnormal{sign}\cong\pi_{\mu^{t}}$.

On the other hand, the image of $\operatorname{pr}_{i,\phi}$ is induced from a representation whose $S_{i}$-factor is $\pi_{\phi^{t}}$. Since irreducible representations of $S_{i}$ are pairwise non-isomorphic, the $S_{i}$-isotypic components corresponding to $\pi_{\mu^{t}}$ and $\pi_{\phi^{t}}$ are different whenever $\mu\neq\phi$. Consequently,
 \[
   \mathrm{pr}_{i,\phi}(v_{T}\otimes e_{1}\wedge...\wedge e_{i})=0
   \]
   whenever $\operatorname{shape}(T_{\leq i})=\mu\neq\phi$.
This proves the lemma.
\end{proof}

The remaining pairs of the \(k\)-th quasi-symmetric types are handled in Case 1 of the proof of Lemma \ref{lem vanishing on gt   basis}, for which one exploits some explicit structure of \(W\)-modules in Section \ref{ss JM elements}.

%By Theorem \ref{dim of Hom}, we can reduce to consider $\nu$ and $\phi$ to be \(k\)-th column quasi-symmetric. 

%The vanishing differentials then follow from stronger vanishing formulas in Lemma \ref{lem vanishing on gt basis}. The main ingredient in the proof of Lemma \ref{lem vanishing on gt basis} is an explicit structure of $W$-type modules.

 \subsection{Partial actions from $\widetilde{e}_j$}
Before computing the differential $d_{i}^{\vee}$ in the complex \eqref{eq:ext-complex}, we provide an explicit action formula of $ \widetilde{e_{j}}$ on $\mathbb{H}_{n}$-modules $X$ and $Y$.

\begin{lemma}\label{lem: action of e_j}
Let \(\lambda\) be a Young diagram. The action of $ \widetilde{e_{j}}$ on $\mathbb{H}_{n}$-module $\tau_{z}^{*}(\lambda)$ is as follows:
    \begin{equation}
         \widetilde{e_{j}}|_{\tau_{z}^{*}(\pi)}=z+\frac{1}{2}X_{j}+\frac{1}{2}Y_{j},
    \end{equation}
    for $1\leq j\leq n$, where $X_{j}=\sum_{t<j}s_{t,j}$ and $Y_{j}=\sum_{s>j}s_{j,s}$.
    
    Specifically,  $ \widetilde{e_{j}}$ acts 
    on $\tau_{0}^{*}(l^{d})$ as:
    \begin{equation}
        \widetilde{e_{j}}|_{\tau_{0}^{*}(l^{d})}=\frac{1}{2}X_{j}+\frac{1}{2}Y_{j}.
\end{equation}
If $l\geq d$, then $\widetilde{e_{j}}$ acts 
    on $\tau_{k}^{*}\big((l-k)^{d}),d^{k}\big)$ as: 
\begin{equation}
    \widetilde{e_{j}}|_{\tau_{k}^{*}\big((l-k)^{d}),d^{k}\big)}=k+\frac{1}{2}X_{j}+\frac{1}{2}Y_{j}.
\end{equation}
If $l<d$, then $\widetilde{e_{j}}$ acts 
    on $\tau_{-k}^{*}\big((l+k)^{l}),l^{d-l-k}\big)$ as: 
\begin{equation}
    \widetilde{e_{j}}|_{\tau_{-k}^{*}\big((l+k)^{l}),l^{d-l-k}\big)}=-k+\frac{1}{2}X_{j}+\frac{1}{2}Y_{j}.
\end{equation}
\end{lemma}

\begin{proof}
By the definition of $\tilde{e_{j}}$ as in \eqref{eq:v-tilde}:
    \begin{eqnarray*}
        \widetilde{e_{j}}
        &=&e_{j}-\frac{1}{2}\sum_{t<s}(-1)\langle e_{j},(e_{t}-e_{s})^{\vee}\rangle s_{t,s}\\
        &=&e_{j}-\frac{(-1)}{2}\sum_{s>j}s_{j,s}+\frac{(-1)}{2}\sum_{t<j}s_{t,j}\\
        &=&(z-c\sum_{t<j}s_{t,j})-\frac{c}{2}\sum_{s>j}s_{j,s}+\frac{c}{2}\sum_{t<j}s_{t,j}\\
        &=&z-\frac{c}{2}\sum_{t<j}s_{t,j}-\frac{c}{2}\sum_{s>j}s_{j,s}\\
        &=&z-\frac{c}{2}X_{j}-\frac{c}{2}Y_{j}\\
        &=&z+\frac{1}{2}X_{j}+\frac{1}{2}Y_{j}\qquad (\textnormal{we set $c=-1$ in Section \ref{Sec: Hecke alg}}).
    \end{eqnarray*}
    
Let $z=0,k,-k$ respectively, we get the three explicit action on the above (quasi-)rectangular \(W\)-type modules.
\end{proof}

\begin{proposition}\label{prop: computation of d_i}
Let $\lambda=(l^{d})$, $l\geq d$, $\mu=\big((l-k)^{d},d^{k}\big)$ for some $k\geq0$. Let $X=\tau_{0}^{*}(\lambda)$ and $Y=\tau_{k}^{*}(\mu)$. Let 
$f\in\text{Hom}_{S_{n}}(X\otimes\bigwedge^{i}\mathbb{C}^{n},\ Y)$ and $v\in X$. Then
\begin{equation}\label{d_i(f)}
    d_{i}^{\vee}(f)(v\otimes e_{1}\wedge...\wedge e_{i+1})=
    (-1)^{i+1}\Bigg\{\sum_{j=1}^{i} s_{j,i+1}.f\big[(k-X_{i+1})s_{j,i+1}.v\otimes\omega_{\hat{i+1}} \big]-f\big[(k-X_{i+1}).v\otimes\omega_{\hat{i+1}}\big]\Bigg\},
\end{equation}
where $X_{i}$ is the $i$-th Jucys-Murphy element as in Lemma \ref{lem: action of e_j}.
\end{proposition}

\begin{proof}
   For convenience, denote $\omega=e_{1}\wedge...\wedge e_{i+1}$, $\omega_{\hat{j}}=e_{1}\wedge...\wedge \hat{e_{j}}\wedge...\wedge e_{i+1}$.
    We first compute $d_{i}^{\vee}(f)(v\otimes e_{1}\wedge...\wedge e_{i+1})$:
    \begin{eqnarray*}
     d_{i}^{\vee}(f)(v\otimes e_{1}\wedge...\wedge e_{i+1})
    &=&\sum_{j=1}^{i+1} (-1)^{j+1} \bigg\{ \widetilde{e}_j|_Y. f(v \otimes \omega_{\hat{j}}) \ - \ f(\tilde{e}_j|_X. v \otimes \omega_{\hat{j}}) \bigg\}\\
    &=& \sum_{j=1}^{i+1} (-1)^{j+1}f(k\cdot v\otimes\omega_{\hat{j}})\\
    &&+\frac{1}{2}\sum_{j=1}^{i+1} (-1)^{j+1}\bigg\{X_{j}.f(v\otimes\omega_{\hat{j}})-f(X_{j}.v\otimes\omega_{\hat{j}})\bigg\}\\
    &&+\frac{1}{2} \sum_{j=1}^{i+1} (-1)^{j+1}\bigg\{Y_{j}.f(v\otimes\omega_{\hat{j}})-f(Y_{j}.v\otimes\omega_{\hat{j}})\bigg\},
\end{eqnarray*}
where the first equality follows from the definition of the differential in Section \ref{ss koszul resol} and the second equality follows from Lemma \ref{lem: action of e_j}.

For the second term:
\begin{eqnarray*}
    &&\sum_{j=1}^{i+1} (-1)^{j+1}\bigg\{X_{j}.f(v\otimes\omega_{\hat{j}})-f(X_{j}.v\otimes\omega_{\hat{j}})\bigg\}\\
    &&=\sum_{j=2}^{i+1} (-1)^{j+1}\bigg\{X_{j}.f(v\otimes\omega_{\hat{j}})-f(X_{j}.v\otimes\omega_{\hat{j}})\bigg\}\\
    &&=\sum_{j=2}^{i+1} (-1)^{j+1}\bigg\{\sum_{a=1}^{j-1}s_{a,j}.f(v\otimes\omega_{\hat{j}})-f(\sum_{a=1}^{j-1}s_{a,j}.v\otimes\omega_{\hat{j}}) \bigg\} \quad (\mbox{by the definition of $X_i$}) \\
    &&=\sum_{j=2}^{i+1} (-1)^{j+1}\bigg\{\sum_{a=1}^{j-1}f(s_{a,j}.v\otimes s_{a,j}.\omega_{\hat{j}})-f(\sum_{a=1}^{j-1}s_{a,j}.v\otimes\omega_{\hat{j}}) \bigg\}  \quad (\mbox{by the $S_n$-invariance of $f$})\\
    &&=\sum_{j=2}^{i+1} (-1)^{j+1}\sum_{a=1}^{j-1}f\big[s_{a,j}.v\otimes (s_{a,j}.\omega_{\hat{j}}-\omega_{\hat{j}})\big] \quad (\mbox{by the action of $s_{a,j}$ on $\omega_{\hat{j}}$})\\
    &&=\sum_{j=2}^{i+1} (-1)^{j+1}\sum_{a=1}^{j-1}f\big[s_{a,j}.v\otimes ((-1)^{j-a-1}\omega_{\hat{a}}-\omega_{\hat{j}})\big]\\
    &&=\sum_{j=2}^{i+1}\sum_{a=1}^{j-1}(-1)^{2j-a}f(s_{a,j}.v\otimes\omega_{\hat{a}}) - \sum_{j=2}^{i+1}\sum_{a=1}^{j-1}(-1)^{j+1}f(s_{a,j}.v\otimes\omega_{\hat{j}}) \\
    &&=\sum_{j=2}^{i+1}\sum_{a=1}^{j-1}(-1)^{a}f(s_{a,j}.v\otimes\omega_{\hat{a}}) - \sum_{j=2}^{i+1}\sum_{a=1}^{j-1}(-1)^{j+1}f(s_{a,j}.v\otimes\omega_{\hat{j}}) \\
    &&=\sum_{a=2}^{i+1}\sum_{j=1}^{a-1}(-1)^{j}f(s_{a,j}.v\otimes\omega_{\hat{j}}) - \sum_{j=2}^{i+1}\sum_{a=1}^{j-1}(-1)^{j+1}f(s_{a,j}.v\otimes\omega_{\hat{j}}) \quad \mbox{(by swapping the labels $a$ and $j$)}\\
    &&=\sum_{j=1}^{i}\sum_{a=j+1}^{i+1}(-1)^{j}f(s_{a,j}.v\otimes\omega_{\hat{j}}) - \sum_{j=2}^{i+1}\sum_{a=1}^{j-1}(-1)^{j+1}f(s_{a,j}.v\otimes\omega_{\hat{j}}) \quad \mbox{(by interchanging  two summation signs)} \\
    &&=\sum_{j=1}^{i+1}(-1)^{j}\sum_{a=1,a\neq j}^{i+1} f(s_{a,j}.v\otimes\omega_{\hat{j}}).
\end{eqnarray*}

Similarly, for the third term:
\begin{eqnarray*}
    &&\sum_{j=1}^{i+1} (-1)^{j+1}\bigg\{Y_{j}.f(v\otimes\omega_{\hat{j}})-f(Y_{j}.v\otimes\omega_{\hat{j}})\bigg\}\\
    &&=\sum_{j=1}^{i+1} (-1)^{j+1}\sum_{b=j+1}^{n}f\big[s_{b,j}.v\otimes (s_{b,j}.\omega_{\hat{j}}-\omega_{\hat{j}})\big] \quad \mbox{(by the definition of $Y_j$)}\\
    &&=\sum_{j=1}^{i} (-1)^{j+1}\sum_{b=j+1}^{i+1}f\big[s_{b,j}.v\otimes (s_{b,j}.\omega_{\hat{j}}-\omega_{\hat{j}})\big]\\
    &&=\sum_{j=1}^{i} (-1)^{j+1}\sum_{b=j+1}^{i+1}f\big[s_{b,j}.v\otimes (-1)^{b-j-1}\omega_{\hat{b}}\big] - \sum_{j=1}^{i} (-1)^{j+1}\sum_{b=j+1}^{i+1}f\big(s_{b,j}.v\otimes\omega_{\hat{j}}\big) \quad \mbox{(by the $s_{b,j}$-action on $\omega_{\hat{b}}$)}\\
    &&=\sum_{j=1}^{i} \sum_{b=j+1}^{i+1}(-1)^{b}f\big(s_{b,j}.v\otimes \omega_{\hat{b}}\big) - \sum_{j=1}^{i} (-1)^{j+1}\sum_{b=j+1}^{i+1}f\big(s_{b,j}.v\otimes\omega_{\hat{j}}\big)   \\
    &&=\sum_{b=1}^{i} \sum_{j=b+1}^{i+1}(-1)^{j}f\big(s_{b,j}.v\otimes \omega_{\hat{j}}\big) - \sum_{j=1}^{i} (-1)^{j+1}\sum_{b=j+1}^{i+1}f\big(s_{b,j}.v\otimes\omega_{\hat{j}}\big)  \quad \mbox{(by swapping indices $b$ and $j$)}\\
    &&=\sum_{j=2}^{i+1}(-1)^{j} \sum_{b=1}^{j-1}f\big(s_{b,j}.v\otimes \omega_{\hat{j}}\big) - \sum_{j=1}^{i} (-1)^{j+1}\sum_{b=j+1}^{i+1}f\big(s_{b,j}.v\otimes\omega_{\hat{j}}\big)  \quad  \mbox{(by interchanging two summation signs)}\\
    &&=\sum_{j=1}^{i+1}(-1)^{j} \sum_{b=1,b\neq j}^{i+1}f\big(s_{b,j}.v\otimes \omega_{\hat{j}}\big).
\end{eqnarray*}

Substituting these two calculated results back into the first equation, we get:
\begin{eqnarray*}
    d_{i}^{\vee}(f)(v\otimes\omega)
    &=&\sum_{j=1}^{i+1} (-1)^{j+1}f(k\cdot v\otimes\omega_{\hat{j}})\\
    &&+\frac{1}{2}\sum_{j=1}^{i+1}(-1)^{j}\sum_{a=1,a\neq j}^{i+1}f(s_{a,j}.v\otimes\omega_{\hat{j}})\\
    &&+\frac{1}{2}\sum_{j=1}^{i+1}(-1)^{j} \sum_{b=1,b\neq j}^{i+1}f\big(s_{b,j}.v\otimes \omega_{\hat{j}}\big)\\
    &=&\sum_{j=1}^{i+1} (-1)^{j+1}f(k\cdot v\otimes\omega_{\hat{j}})-\sum_{j=1}^{i+1}(-1)^{j+1}\sum_{a=1,a\neq j}^{i+1}f(s_{a,j}.v\otimes\omega_{\hat{j}})\\
    &=&\sum_{j=1}^{i+1} (-1)^{j+1}f\big[(k-\sum_{a=1,a\neq j}^{i+1}s_{a,j}).v\otimes\omega_{\hat{j}} \big]  \quad \mbox{(by combining two terms)} \\
    &=&\sum_{j=1}^{i+1} (-1)^{j+1}f\big[s_{j,i+1}(k-X_{i+1})s_{j,i+1}.v\otimes\omega_{\hat{j}} \big] \quad \mbox{(by $ \sum_{a=1, a\neq j}^{i+1}s_{a,  j} =s_{j,i+1}X_{i+1}s_{j,i+1}$)}\\
    &=&\sum_{j=1}^{i+1} (-1)^{j+1}s_{j,i+1}.f\big[(k-X_{i+1})s_{j,i+1}.v\otimes s_{j,i+1}.\omega_{\hat{j}} \big]   \quad \mbox{(by $S_n$-invariance of $f$)} \\
    &=&\sum_{j=1}^{i} (-1)^{j+1}s_{j,i+1}.f\big[(k-X_{i+1})s_{j,i+1}.v\otimes (-1)^{i-j}\omega_{\hat{i+1}} \big]+(-1)^{i+2}f\big[(k-X_{i+1}).v\otimes\omega_{\hat{i+1}}\big]\\
    & &  \quad \ \quad \ \quad \mbox{(by $s_{j,i+1}$-action on $\omega_{\hat{i+1}}$)} \\
    &=&\sum_{j=1}^{i} (-1)^{i+1}s_{j,i+1}.f\big[(k-X_{i+1})s_{j,i+1}.v\otimes\omega_{\hat{i+1}} \big]+(-1)^{i+2}f\big[(k-X_{i+1}).v\otimes\omega_{\hat{i+1}}\big]\\
    &=&(-1)^{i+1}\Bigg\{\sum_{j=1}^{i} s_{j,i+1}.f\big[(k-X_{i+1})s_{j,i+1}.v\otimes\omega_{\hat{i+1}} \big]-f\big[(k-X_{i+1}).v\otimes\omega_{\hat{i+1}}\big]\Bigg\} .
\end{eqnarray*}
Therefore, we have proved Proposition \ref{prop: computation of d_i}.
\end{proof}

\begin{proposition}\label{prop: computation of d_i for l<d}
    Suppose $l<d$. Let $\lambda=(l^{d})$, and let $\mu=\big((l+k)^{l}),l^{d-l-k}\big)$,  for some $k\geq0$. Let $X=\tau_{0}^{*}(\lambda)$ and $Y=\tau_{-k}^{*}(\mu)$. Let 
    $f\in\text{Hom}_{S_{n}}(X\otimes\bigwedge^{i}\mathbb{C}^{n},\ Y)$ and $v\in X$. Then
    \begin{equation}\label{d_i(f)2}
    d_{i}^{\vee}(f)(v\otimes e_{1}\wedge...\wedge e_{i+1})=
    (-1)^{i}\Bigg\{\sum_{j=1}^{i} s_{j,i+1}.f\big[(k+X_{i+1})s_{j,i+1}.v\otimes\omega_{\hat{i+1}} \big]-f\big[(k+X_{i+1}).v\otimes\omega_{\hat{i+1}}\big]\Bigg\},
   \end{equation}
   where $X_{i}$ is the $i$-th Jucys-Murphy element as in Lemma \ref{lem: action of e_j}.

\end{proposition}
\begin{proof}
    Repeat the calculation in the proof of Proposition \ref{prop: computation of d_i}.
\end{proof}

\subsection{Vanishing differentials}

%\begin{lemma}
%Let \( \phi \vdash i, \nu \vdash i+1\). Let $f \in \mathrm{Hom}_{S_n}(X\otimes \bigwedge^i\mathbb C^n, Y)$. If 
%\[   d_i^{\vee}(f\circ \mathrm{pr}_{i,\nu})\circ \mathrm{pr}_{i+1, \nu'} \neq 0,
%\]
%then \(\mathrm{cont}_T(\boxed{i+1}) =k\).
%\end{lemma}

\begin{lemma} \label{lem vanishing on gt   basis}

Let $f \in \mathrm{Hom}_{S_n}(X\otimes \bigwedge^i\mathbb C^n, Y)$ and let $v_T$ be an element in the Gelfand-Tsetlin basis of $X$. 
\begin{enumerate}
    \item If $l\geq d$, assume that $T$ satisfies the $k$-th column type quasi-symmetric condition. Then 
    \[   f\circ\operatorname{pr}_{i,\phi}\big[(k-X_{i+1}).v_T \otimes e_{1}\wedge...\wedge e_{i}\big] =0.
    \]
    \item If $l<d$, assume that $T$ satisfies the $k$-th row type quasi-symmetric condition. Then 
    \[   f\circ\operatorname{pr}_{i,\phi}\big[(k+X_{i+1}).v_T \otimes e_{1}\wedge...\wedge e_{i}\big] =0 .
     \]
     \end{enumerate}
\end{lemma}

\begin{proof}
 We prove only the first statement here, the other statement is similar. 
 
 We divide the proof into two situations:
            \begin{enumerate}
            
            \item Case 1: \( \mathrm{cont}_{T}(\boxed{i+1})=k\). Then by the scalar action of the Jucys-Murphy element $X_{i+1}$ on the Gelfand–Tsetlin basis as in (\ref{eq:content-eigenvalue}),
            \[  X_{i+1}.v_{T}=kv_{T} . \]
            Hence we have:
             \[f\circ\operatorname{pr}_{i,\phi}\big[(k-X_{i+1}).v_{T} \otimes \omega_{\hat{i+1}}\big]=(k-\operatorname{cont}_{T}(\boxed{i+1}))f\circ\operatorname{pr}_{i,\phi}(v_{T}\otimes\omega_{\hat{i+1}})=0
             \]
             
             \item
            Case 2: \( \mathrm{cont}_{T}(\boxed{i+1})\neq k \).  This means that $\boxed{i+1}$ is not on the diagonal line of  the symmetric part of $T_{\leq i+1}$ (by removing the first $k$ columns). If we remove the box $\boxed{i+1}$ in $T_{\leq i+1}$, then the remaining part $T_{\leq i}$ will not satisfy the $k$-th column type quasi-symmetric condition, because the only way to obtain a symmetric Young diagram by removing a box from a symmetric Young diagram is to remove the box on the top-left to bottom-right diagonal line. 
           Hence by Lemma \ref{projection}, $\operatorname{pr}_{i, \tau}(v_{T}\otimes\omega_{\hat{i+1}})=0$ for any $\tau\vdash i$ which satisfies the $k$-th column type quasi-symmetric condition.
           Therefore we have 
           \[f\circ\operatorname{pr}_{i,\phi}\big[(k-X_{i+1}).v_{T} \otimes \omega_{\hat{i+1}}\big]=0.
             \]
           \end{enumerate}
\end{proof}

\begin{example}
     We use Example \ref{eg: dim of Hom of rec and quasi rec} to give a visual explanation of the key step in the proof of Lemma \ref{lem vanishing on gt basis}: As in Example \ref{eg: dim of Hom of rec and quasi rec}, $k=2$.
            
          Let \[T= \begin{ytableau}
          *(cyan)1  & *(cyan)2  & *(red!40)5  & *(red!40)9  & *(red!40)11& 13  & 14 \\
        *(cyan)3  & *(cyan)4  & *(red!40)7  & *(red!40)12  & *(white)15  &16   &  17 \\
        *(cyan)6  & *(cyan)8  & *(red!40)10  & *(white) 18 & *(white) 19 & 20  & 21
       \end{ytableau}.
       \]
        Then $T_{\leq12}$ is $2$-nd column type quasi-symmetric, $ \mathrm{cont}_{T}(\boxed{12})=2$,  $T_{\leq11}$ is $2$-th column type quasi-symmetric.

         Let 
         \[T'= \begin{ytableau}
          *(cyan)1  & *(cyan)2  & *(red!40)5  & *(red!40)9  & *(red!40)11& 13  & 14 \\
         *(cyan)3  & *(cyan)4  & *(red!40)7  & *(red!40)10  & *(white)15  &16   &  17 \\
         *(cyan)6  & *(cyan)8  & *(red!40)12  & *(white) 18 & *(white) 19 & 20  & 21
          \end{ytableau}.
         \]
        Then $T'_{\leq12}$ is $2$-th column type quasi-symmetric, $ \mathrm{cont}_{T'}(\boxed{12})=0$,  $T‘_{\leq11}$ is not $2$-th column type quasi-symmetric.
\end{example}

\begin{theorem}\label{d_i=o}
Let $X=\tau_{0}^{*}(l^{d})$. If 
$l\geq d$, let 
\[ Y=  \tau_{k}^{*}((l-k)^{d}),d^{k}) ,\]
for some $k \geq 0$. For $l<d$, let 
\[ Y=\tau_{-k}^{*}((l+k)^{l},l^{d-l-k}),\] 
for some \( k\geq0\). In both cases, the differentials in the complex \eqref{eq:ext-complex} vanish, i.e. 
    $d_{i}^{\vee}=0$ for any $i$.
\end{theorem}

%\begin{proof}
%    \textcolor{red}{
%    This is a direct consequence of Proposition \ref{prop: computation of d_i}, Lemma \ref{lem vanishing on gt basis}, Proposition \ref{prop: computation of d_i for l<d} and Lemma \ref{lem vanishing on gt basis for l<d}.
%    }
%    We just prove the $l\geq d$ case and the $l<d$ case is similar. Let $f\in\text{Hom}_{S_{n}}(X\otimes\bigwedge^{i}\mathbb{C}^{n},\ Y)$ and let $v_{T}$ be a Gelfand–Tsetlin basis element of $X$.
    
%    By Proposition \ref{prop: computation of d_i} and Lemma \ref{lem vanishing on gt basis}, 
%    \begin{eqnarray*}
%    d_{i}^{\vee}(f)(v_{T}\otimes e_{1}\wedge...\wedge e_{i+1})&&=
%    \sum_{j=1}^{i}\bigg\{(-1)^{i+1} s_{j,i+1}.f\big[(k-X_{i+1})s_{j,i+1}.v_{T}\otimes\omega_{\hat{i+1}} \big]\\
%    &&\quad+(-1)^{i+2}f\big[(k-X_{i+1}).v_{T}\otimes\omega_{\hat{i+1}}\big]\bigg\}\\
%    &&=0.
%   \end{eqnarray*}
%    Since every $e_{j_{1}}\wedge e_{j_{2}}\wedge...\wedge e_{j_{i+1}}=\sigma.(e_{1}\wedge...\wedge e_{i+1})$ for some $\sigma\in S_{n}$, it's easy to see that $d_{i}^{\vee}(f)(v_{T}\otimes e_{j_{1}}\wedge e_{j_{2}}\wedge...\wedge e_{j_{i+1}})=0$ for any $v_{T}\in X$, then we get $d_{i}^{\vee}(f)=0$ for some $f\in\text{Hom}_{S_{n}}(X\otimes\bigwedge^{i}\mathbb{C}^{n},\ Y)$. 
%    we have $d_{i}^{\vee}=0$.
%\end{proof}

\begin{proof}
       We give a proof only for the $l\geq d$ case; the other case is similar.
       
       By Proposition \ref{prop: decomposition of Hom in section vanish}, we have two decompositions:
       \begin{eqnarray*}
           \mathrm{Hom}_{S_n}(X\otimes \bigwedge^i \mathbb C^n, Y) \cong \bigoplus_{\phi \vdash i} \mathrm{Hom}_{S_i \times S_{n-i}}(\pi_{\phi^t}\otimes \pi_{(\lambda-\phi)^{\circlearrowleft}}, Y),
       \end{eqnarray*}
       \begin{eqnarray*}
           \mathrm{Hom}_{S_n}(X\otimes \bigwedge^{i+1} \mathbb C^n, Y) \cong \bigoplus_{\nu \vdash i+1} \mathrm{Hom}_{S_{i+1} \times S_{n-i-1}}(\pi_{\nu^t}\otimes \pi_{(\lambda-\nu)^{\circlearrowleft}}, Y),
       \end{eqnarray*}
       where $\phi$ and $\nu$ satisfy $k$-th column type quasi-symmetric condition. Thus if $f\in \mathrm{Hom}_{S_n}(X\otimes \bigwedge^{i} \mathbb C^n, Y)$, then we have a decomposition:
       \begin{equation*}
           d_{i}^{\vee}(f)=\sum_{\phi\vdash i}d_{i}^{\vee}(f\circ\operatorname{pr}_{i,\phi})=\sum_{\phi\vdash i}\sum_{\nu\vdash i+1}d_{i}^{\vee}(f\circ\operatorname{pr}_{i,\phi})\circ\operatorname{pr}_{i+1,\nu}.
       \end{equation*}
       
       We firstly show that 
       \[
       d_{i}^{\vee}(f\circ\operatorname{pr}_{i,\phi})(v_T\otimes e_{1}\wedge...\wedge e_{i+1})=0
       \]
       for any 
       $f\in\mathrm{Hom}_{S_n}(X\otimes \bigwedge^i \mathbb C^n, Y)$, any Gelfand–Tsetlin basis element $v_{T}\in X$ and any quasi-symmetric $\phi\vdash i$.

       \begin{enumerate}
           \item 
           If $T_{\leq i+1}$ does not satisfy the $k$-th column type quasi-symmetric condition, then by Lemma \ref{projection}, 
            \begin{equation*}
           \operatorname{pr}_{i+1,\nu}(v_{T}\otimes e_{1}\wedge...\wedge e_{i+1})=0
            \end{equation*} 
            for any $k$-th column type quasi-symmetric $\nu\vdash i+1$. Hence 
           \begin{equation*}
         d_{i}^{\vee}(f\circ\operatorname{pr}_{i,\phi})(v_T\otimes e_{1}\wedge...\wedge e_{i+1})=\sum_{\nu\vdash i+1} d_{i}^{\vee}(f\circ\operatorname{pr}_{i,\phi})\circ\operatorname{pr}_{i+1,\nu}(v_T\otimes e_{1}\wedge...\wedge e_{i+1})=0.
          \end{equation*}

           \item
            From now on, assume that $T_{\leq i+1}$ satisfies the $k$-th column type quasi-symmetric condition, so
            \[
            \operatorname{shape}(T_{\leq i+1})=\operatorname{shape}(\nu),
            \]
          for some $k$-th column type quasi-symmetric $\nu\vdash i+1$. 

             By Proposition \ref{prop: computation of d_i}, we have:
            \begin{eqnarray*}
             d_{i}^{\vee}(f\circ\operatorname{pr}_{i,\phi})(v_T\otimes e_{1}\wedge...\wedge e_{i+1})=
             &&(-1)^{i+1}\sum_{j=1}^{i} s_{j,i+1}.f\circ\operatorname{pr}_{i,\phi}\big[(k-X_{i+1})s_{j,i+1}.v_T\otimes\omega_{\hat{i+1}} \big]\\
             &&+(-1)^{i+2}f\circ\operatorname{pr}_{i,\phi}\big[(k-X_{i+1}).v_T\otimes\omega_{\hat{i+1}}\big].
             \end{eqnarray*}
       
            For $1\leq j\leq i+1$, we know that $s_{j,i+1}$ just changes the first $i+1$ numbers of $v_{T}$ by Theorem \ref{thm: Young's orthogonal form}, so we can express $s_{j,i+1.}.v_{T}$ as: $$s_{j,i+1}.v_{T}=\sum_{t=1}^{N(j)}a_{j,t}v_{T_{t}}\ , $$
            where $\operatorname{shape}(T_{t,\leq i+1})=\operatorname{shape}(T_{\leq i+1})=\operatorname{shape}(\nu)$.

            By Lemma \ref{lem vanishing on gt   basis}, we have
            \[f\circ\operatorname{pr}_{i,\phi}\big[(k-X_{i+1})v_{T_{t}}\otimes\omega_{\hat{i+1}} \big] =0,
            \]
            hence \[
            d_{i}^{\vee}(f\circ\operatorname{pr}_{i,\phi})(v_{T}\otimes e_{1}\wedge...\wedge e_{i})=0.
            \]

    \end{enumerate}
       
      Combining $(1)$ and $(2)$, we obtain that
       \[
       d_{i}^{\vee}(f\circ\operatorname{pr}_{i,\phi})(v_T\otimes e_{1}\wedge...\wedge e_{i+1})=0
       \]
       for any quasi-symmetric $\phi\vdash i$, any $f\in\mathrm{Hom}_{S_n}(X\otimes \bigwedge^i \mathbb C^n, Y)$ and any Gelfand–Tsetlin basis element $v_{T}\in X$ .
       
Since every $e_{j_{1}}\wedge e_{j_{2}}\wedge...\wedge e_{j_{i+1}}=\sigma.(e_{1}\wedge...\wedge e_{i+1})$ for some $\sigma\in S_{n}$, it's easy to see that
\[d_{i}^{\vee}(f\circ\operatorname{pr}_{i,\phi})(v_{T}\otimes e_{j_{1}}\wedge e_{j_{2}}\wedge...\wedge e_{j_{i+1}})=0,
\] 
 hence we have
 \[
 d_{i}^{\vee}(f\circ\operatorname{pr}_{i,\phi})=0.
 \] 
 Therefore \[
 d_{i}^{\vee}(f)=\sum_{\phi\vdash i}d_{i}^{\vee}(f\circ\operatorname{pr}_{i,\phi})=0
 \]
 for any $f\in \mathrm{Hom}_{S_n}(X\otimes \bigwedge^i \mathbb C^n, Y)$, thus $d_{i}^{\vee}=0$.
\end{proof}

\begin{theorem} \label{thm ext in term hom}
    Let $X$ be a Speh $\mathbb H_n$-module  and let $Y$ be a $W$-type $\mathbb{H}_{n}$-module  with the same central character. Then
    \begin{equation}
    \textnormal{Ext}_{\mathbb{H}_{n}}^{i}(X, Y)\cong \textnormal{Hom}_{S_{n}}(X\otimes\bigwedge^{i}\mathbb{C}^{n},Y).
    \end{equation}
\end{theorem}

\begin{proof}
This follows from Theorem \ref{d_i=o}.
\end{proof}

\begin{remark}
    We can use Theorem \ref{thm ext in term hom}, the combinatorial formula for Hom space in Theorem \ref{dim of Hom} and Corollary \ref{cor: dim of Hom for l<d case} to compute $\textnormal{Ext}$ groups for one Speh module and one general \(W\)-type module.
\end{remark}

\section{The first non-vanishing Ext group between two quasi-rectangular type modules}\label{sec first non-vanish Ext for two quasi-rec}

In this section, unless otherwise specified, let \(\lambda =((l-s)^{d} ,d^{s})\) and \((l-t)^{d}),d^{t})\). Let $X=\tau_{s}^{*}\big(\lambda\big)$ and $Y=\tau_{t}^{*}\big(\mu\big)$. According to Proposition \ref{shape for same central character}, $X$ and $Y$ have the same central character. When no confusion can arise, we sometimes use $X$ and $Y$ to represent $\textnormal{Res}_{S_{n}}X$ and $\textnormal{Res}_{S_{n}}Y$ respectively.

\begin{proposition}\label{dim Hom for two Hook type}

    \begin{eqnarray*}
        \textnormal{dim Hom}_{S_{n}}(X\otimes\bigwedge^{i}\mathbb{C}^{n},Y)=\begin{cases}
            0,\quad &i<d|t-s|,\\
            1,\quad& i=d|t-s|,\\
               1,\quad& i=d|t-s|+1,\  \textnormal{and}\ s(l-d-t)=0\\
               2,\quad& i=d|t-s|+1,\  \textnormal{and}\ s(l-d-t)\neq0
        \end{cases}
        \end{eqnarray*}
\end{proposition}

\begin{proof}
By taking dual if necessary, we shall assume that $t\geq s\geq 0$. By the dimension formula in Theorem \ref{dim of Hom},
    \begin{equation}
           \textnormal{dim}\ \textnormal{Hom}_{S_{n}}(X\otimes\bigwedge^{i}\mathbb{C}^{n}, Y)= \sum\limits_{\nu\vdash n-i}\sum\limits_{\xi\vdash i} \textnormal{LR}^{\lambda}_{\nu,\ \xi}\textnormal{LR}^{\mu}_{\nu,\ \xi^{t}}.
        \label{eqn dimension littlewood richardson}
        \end{equation}
        
When $i<d(t-s)$, 
\[ |\nu|=n-i>n-d(t-s)=|\lambda\cap \mu|, \]
hence no $\nu\vdash n-i$ is contained in $\lambda\cap \mu$, and so \( \mathrm{LR}^\lambda_{\nu, \xi}\mathrm{LR}^\mu_{\nu, \xi^t}=0\) for any \( \nu \vdash n-i\). Thus, \eqref{eqn dimension littlewood richardson} implies $\textnormal{dim}\ \textnormal{Hom}_{S_{n}}(X\otimes\bigwedge^{i}\mathbb{C}^{n}, Y)=0$.

When $i=d(t-s)$, 
\begin{eqnarray*}
   \textnormal{dim}\ \textnormal{Hom}_{S_{n}}(X\otimes\bigwedge^{i}\mathbb{C}^{n}, Y)
   &=& \sum\limits_{\xi\vdash i} \textnormal{LR}^{\lambda}_{\lambda\cap \mu,\ \xi}\textnormal{LR}^{\mu}_{\lambda\cap \mu,\ \xi^{t}}\\
   &=& \textnormal{LR}^{\lambda}_{\lambda\cap \mu,\ \xi^*}\textnormal{LR}^{\mu}_{\lambda\cap \mu,\ (\xi^*)^{t}},\quad \textnormal{where } \xi^*=((t-s)^{d})\\
   &=& 1 ,
\end{eqnarray*}
where the first equation follows from the fact that $\lambda\cap \mu$ is the unique Young subdiagram of both $((l-s)^d,d^s)$ and $((l-t)^d, d^t)$ with $d(t-s)$ boxes; the second equation follows from that $\xi^*$ is uniquely determined by $\lambda\cap \mu$, and the third equation follows from the Littlewoord--Richardson rule.

When $i=d(t-s)+1$, if $s(l-d-t)=0$, i.e., one of $\lambda$ and $\mu$ is rectangular, then this case was proved before, so we only consider the case $s(l-d-t)\neq0$. In this case, there exist two diagrams in $\lambda\cap \mu$ with $n-i$ boxes: $\nu_{1}=((l-t+s)^{d-1}),l-t+s-1,d^{s})$ and $\nu_{2}=((l-t+s)^{d},d^{s-1},d-1)$. Then 
\begin{eqnarray*}
   \textnormal{dim}\ \textnormal{Hom}_{S_{n}}(X\otimes\bigwedge^{i}\mathbb{C}^{n}, Y)
   &=& \sum\limits_{\xi\vdash i} \textnormal{LR}^{\lambda}_{\nu_{1},\ \xi}\textnormal{LR}^{\mu}_{\nu_{1},\ \xi^{t}}
   +\sum\limits_{\xi\vdash i} \textnormal{LR}^{\lambda}_{\nu_{2},\ \xi}\textnormal{LR}^{\mu}_{\nu_{2},\ \xi^{t}}\\
   &=& \textnormal{LR}^{\lambda}_{\nu_{1},\ \xi}\textnormal{LR}^{\mu}_{\nu_{1},\ \xi^{t}}+\textnormal{LR}^{\lambda}_{\nu_{2},\ \xi}\textnormal{LR}^{\mu}_{\nu_{2},\ \xi^{t}},\quad \xi=(t-s+1,(t-s)^{d})\\
   &=& 1+1\\
   &=& 2.
\end{eqnarray*}
\end{proof}

\begin{example}\label{example for dim Hom for two hook type}
Let 
\(\lambda=
\begin{ytableau}
*(white) & *(white) & *(white) & *(white) & *(white) & *(white)  \\
*(white) & *(white) & *(white) & *(white) & *(white) & *(white)  \\
*(white) & *(white) & *(white) & *(white) & *(white) & *(white)  \\
*(white) & *(white) & *(white)
\end{ytableau}\)
and $\mu=$
\begin{ytableau}
*(white) & *(white) & *(white) & *(white) & *(white) \\
*(white) & *(white) & *(white) & *(white) & *(white) \\
*(white) & *(white) & *(white) & *(white) & *(white) \\
*(white) & *(white) & *(white) \\
*(white)  & *(white)  & *(white) 
\end{ytableau}.

We explain how to compute the Littlewood--Richardson coefficients to show that $\textnormal{dim Hom}_{S_{21}}(X\otimes\bigwedge^{3}\mathbb{C}^{21},Y)=1$. We first color the Young diagrams \(\lambda\) and \(\mu\) as follows. 
\begin{center}
$\lambda=$
\begin{ytableau}
*(red!40)  & *(blue!35) & *(blue!35) & *(blue!35) & *(blue!35) & *(white) 1\\
*(red!40)  & *(blue!35) & *(blue!35) & *(blue!35) & *(blue!35) & *(white) 2\\
*(red!40)  & *(blue!35) & *(blue!35) & *(blue!35) & *(blue!35) & *(white) 3\\
*(blue!35) & *(blue!35) & *(blue!35)
\end{ytableau}
, $\mu=$
\begin{ytableau}
*(red!40)  & *(red!40) & *(red!40) & *(blue!35) & *(blue!35)  \\
*(blue!35)  & *(blue!35) & *(blue!35) & *(blue!35) & *(blue!35) \\
*(blue!35)  & *(blue!35) & *(blue!35) & *(blue!35) & *(blue!35) \\
*(blue!35) & *(blue!35) & *(blue!35)\\
*(white) 1& *(white) 1& *(white) 1
\end{ytableau}
\end{center}
The colored diagram is the unique $\nu\vdash 18$ which is contained in both \(\lambda\) and \(\mu\). The red diagram in $\lambda$ is the unique $\xi\vdash 3$, and the red diagram in $\mu$ is the unique $\xi^{t}$, such that $\textnormal{LR}^{\lambda}_{\nu,\xi}$ and $\textnormal{LR}^{\mu}_{\nu,\xi^{t}}$ are nonzero. In this case, both $\textnormal{LR}^{X}_{\nu,\xi}$ and $\textnormal{LR}^{Y}_{\nu,\xi^{t}}$ are equal to \(1\), and the corresponding number-fillings are shown above.

We now explain how to show $\textnormal{dim Hom}_{S_{21}}(X\otimes\bigwedge^{4}\mathbb{C}^{21},Y)=2$. We first color the Young diagrams as follows.

\begin{center}
$\lambda=$
\begin{ytableau}
*(red!40)  & *(red!40) & *(yellow!40) & *(yellow!40) & *(yellow!40) & *(white) 1\\
*(red!40)  & *(yellow!40) & *(yellow!40) & *(yellow!40) & *(yellow!40) & *(white) 2\\
*(red!40)  & *(yellow!40) & *(yellow!40) & *(yellow!40) & *(white) 1& *(white) 3\\
*(yellow!40) & *(yellow!40) & *(yellow!40)
\end{ytableau}
, $\mu=$
\begin{ytableau}
*(red!40)  & *(red!40) & *(red!40) & *(yellow!40) & *(yellow!40)  \\
*(red!40)  & *(yellow!40) & *(yellow!40) & *(yellow!40) & *(yellow!40) \\
*(yellow!40)  & *(yellow!40) & *(yellow!40) & *(yellow!40) & *(white) 1\\
*(yellow!40) & *(yellow!40) & *(yellow!40)\\
*(white) 1& *(white) 1& *(white) 2
\end{ytableau}
;
\end{center}

\begin{center}
$\lambda=$
\begin{ytableau}
*(red!40)  & *(red!40) & *(green!40) & *(green!40) & *(green!40) & *(white) 1\\
*(red!40)  & *(green!40) & *(green!40) & *(green!40) & *(green!40) & *(white) 2\\
*(red!40)  & *(green!40) & *(green!40) & *(green!40) & *(green!40) & *(white) 3\\
*(green!40) & *(green!40) & *(white) 1
\end{ytableau}
, $\mu=$
\begin{ytableau}
*(red!40)  & *(red!40) & *(red!40) & *(green!40) & *(green!40)  \\
*(red!40)  & *(green!40) & *(green!40) & *(green!40) & *(green!40) \\
*(green!40)  & *(green!40) & *(green!40) & *(green!40) & *(green!40) \\
*(green!40) & *(green!40) & *(white) 1\\
*(white) 1& *(white) 1& *(white) 2
\end{ytableau}
.
\end{center}
The colored diagrams are the only $\nu_{1}, \nu_{2}\vdash 17$ which are contained in both \(\lambda\) and \(\mu\). The red diagram in $\lambda$ is the unique  $\xi\vdash 4$, and the red diagram in $\mu$ is the unique $\xi^{t}$ such that the Littlewood-Richardson coefficients are nonzero. In this case, their corresponding Littlewood-Richardson coefficients are equal to $1$, and the corresponding number-fillings are shown above.
\end{example}

\begin{theorem}\label{thm: Ext for two quasi-rectangula}
We use the above notation. The minimal $i \in \mathbb Z_{\geq 0}$ satisfying $\textnormal{Ext}_{\mathbb{H}_{n}}^{i}(X, Y)\neq0$ is $i=|t-s|d$. In this case, 
    \[ \textnormal{Ext}_{\mathbb{H}_{n}}^{|t-s|d}(X, Y)\cong\mathbb{C}. \]
\end{theorem}
\begin{proof}
By Proposition \ref{prop dual ext} if necessary, we may and shall assume that $t \geq s \geq 0$.

    Since $\textnormal{dim Hom}_{S_{n}}(X\otimes\bigwedge^{i}\mathbb{C}^{n},Y)=0$ when $i<(t-s)d$, the resolution in \eqref{eq:ext-complex} implies that the minimal \(i \in \mathbb Z_{\geq 0}\) which satisfy \(\mathrm{Ext}^i_{\mathbb H_n}(X,Y)\neq 0\) is at least \((t-s)d\). It suffices to show that 
    \[  \mathrm{Ext}^{(t-s)d}_{\mathbb H_n}(X, Y) \neq 0 .    \]
    
    We now consider the case $i=(t-s)d$.
    Let $f\in\textnormal{Hom}_{S_{n}}(X\otimes\bigwedge^{i}\mathbb{C}^{n},Y)\cong\mathbb{C}$, and let $v_{T}$ be a Gelfand–Tsetlin element of $X$. By Proposition \ref{dim Hom for two Hook type}, it suffices to show that 
    \[  d_i^{\vee}(f\circ \mathrm{pr}_{i, \lambda\cap \mu}) \circ \mathrm{pr}_{i+1, \nu} =0,
    \]
    where \( \nu \) is a subdiagram of both \(\lambda\) and \(\mu\).  The following argument is similar to the proof of Theorem  \ref{d_i=o}; we only give a brief explanation here.
    
    Following the computation in Proposition \ref{prop: computation of d_i}, we have a similar result:
\begin{eqnarray*}
    d_{i}^{\vee}(f\circ\operatorname{pr}_{i,\lambda\cap\mu})(v_{T}\otimes e_{1}\wedge...\wedge e_{i+1})=&&
    (-1)^{i+1}\Bigg\{\sum_{j=1}^{i} s_{j,i+1}.f\circ\operatorname{pr}_{i,\lambda\cap\mu}\big[(t-s-X_{i+1})s_{j,i+1}.v_{T}\otimes\omega_{\hat{i+1}} \big]\Bigg\}\\
    &&+(-1)^{i+2}f\circ\operatorname{pr}_{i,\lambda\cap\mu}\big[(t-s-X_{i+1}).v_{T}\otimes\omega_{\hat{i+1}}\big].
\end{eqnarray*}

Write 
 $$s_{j,i+1}.v_{T}=\sum_{t=1}^{N(j)}a_{j,t}v_{T_{t}}\ , $$
 for some \(a_{j,t} \in \mathbb C\), 
where \(T_{t}\) runs over \(T \in \textnormal{STY}(\lambda)\) such that $\operatorname{shape}(T_{t,\leq i+1})=\operatorname{shape}(T_{\leq i+1})$.

If the content of $\boxed{i+1}$ in $T_{t}$ is $t-s$, then by the scalar action of the Jucys-Murphy element $X_{i+1}$ on Gelfand–Tsetlin basis $T_{t}$, we have:
\begin{eqnarray*}
    f\circ\operatorname{pr}_{i,\lambda\cap\mu}\big[(t-s-X_{i+1})v_{T_{t}}\otimes\omega_{\hat{i+1}} \big]
    =f\circ\operatorname{pr}_{i,\lambda\cap\mu}\big[(t-s-\textnormal{cont}_{T_{t}}(\boxed{i+1})v_{T_{t}})\otimes\omega_{\hat{i+1}}\big]=0.
\end{eqnarray*}

If the content of $\boxed{i+1}$ is not $t-s$, then the sub-diagram consists of the first $i$ numbers is not in the shape of $((t-s)^{d})$, which is the unique satisfying $\xi^{*}\vdash i$ as in the proof of Proposition \ref{dim Hom for two Hook type}. Then by Lemma \ref{projection} (this lemma also holds for this case), 
\begin{eqnarray*}
    f\circ\operatorname{pr}_{i,\lambda\cap\mu}\big[(t-s-X_{i+1})v_{T_{t}}\otimes\omega_{\hat{i+1}} \big]
    =0.
\end{eqnarray*}
Then as in the proof of Theorem \ref{d_i=o}, we get $d_{i=(t-s)d}=0$. As a result, the minimal $i$ which satisfies $\textnormal{Ext}_{\mathbb{H}_{n}}^{i}(X, Y)\neq0$ is $i=(t-s)d$, and in this case, $\textnormal{Ext}_{\mathbb{H}_{n}}^{i}(X, Y)\cong\textnormal{Hom}_{S_{n}}(X\otimes\bigwedge^{i}\mathbb{C}^{n}, Y)\cong\mathbb{C}$.
\end{proof}

\begin{example} \label{ex minimal w type}
We combine Example \ref{example for dim Hom for two hook type} to explain the key part in the proof of Proposition
\ref{thm: Ext for two quasi-rectangula}. In this case, $t-s=1$, and $i=(t-s)d=3$.

Let $T_{1}=$
\begin{ytableau}
*(red!40) 1 & *(white) 4& *(white)... & *(white)... & *(white)... & *(white)...  \\
*(red!40) 2 & *(white)... & *(white)... & *(white)...  & *(white)... & *(white)...  \\
*(red!40)3 & *(white)... & *(white)... & *(white)... & *(white)... & *(white)...  \\
*(white)... & *(white)... & *(white)...
\end{ytableau},
then $\textnormal{cont}_{T_{1}}(4)=1=t-s$, hence $(t-s-X_{4}).v_{T_{1}}=0$, so \begin{equation*}
    f((t-s-X_{4}).v_{T_{1}}\otimes e_{1}\wedge e_{2}\wedge e_{3}))=f(0)=0.
\end{equation*}

Let $T_{2}=$
\begin{ytableau}
*(red!40) 1 & *(red!40) 2& *(white)... & *(white)... & *(white)... & *(white)...  \\
*(red!40) 3 & *(white)... & *(white)...& *(white)...  & *(white)... & *(white)...  \\
*(white) 4 & *(white)... & *(white)... & *(white)... & *(white)... & *(white)...  \\
*(white)... & *(white)... & *(white)...
\end{ytableau},
then $\textnormal{cont}_{T_{2}}(4)=-2\neq t-s$, hence $(t-s-X_{4}).v_{T_{2}}=3\cdot v_{T_{2}}\neq0$. The red part i.e the first $i=3$ numbers is not in the shape of $(1,1,1)$ as in the first colored graph in Example \ref{example for dim Hom for two hook type}, hence 
\begin{equation*}
    f(3\cdot v_{T_{2}}\otimes e_{1}\wedge e_{2}\wedge e_{3})=0.
\end{equation*}
\end{example}

\begin{remark} \label{rm differentials non vanish ex}
As we saw in Example \ref{ex minimal w type}, the complexes involved in the Koszul resolutions of Section \ref{ss koszul resol} for two quasi-rectangular type $W$-modules are more complicated than those in Theorem \ref{thm ext in term hom}. In particular, the quasi-symmetric condition is not sufficient for computing Ext-groups between two quasi-rectangular type modules. More importantly, the differentials in the Koszul resolution in such cases do not necessarily vanish. For example, 
let $\lambda=(2,1)$; and let $X=\tau^*_0(\lambda)$. Then
\[  \dim\mathrm{Hom}_{S_3}(X\otimes \bigwedge^i \mathbb C^3, X)=\left\{ \begin{array}{cc} 1 & i =0, 3 \\  2 & i=1,2 \end{array} \right.
\]
On the other hand,
\[  \mathrm{Ext}^i_{\mathbb H_3}(X, X) = \left\{ \begin{array}{cc} 1 & i=0,1 \\ 0 & \mbox{otherwise} \end{array} \right. 
\]
Thus, the differentials $d_2^{\vee}$ and $d_3^{\vee}$ cannot vanish.
\end{remark}

%%
%{\color{teal} need to combine the following proposition with the previous one}

%\begin{proposition}
%    \textcolor{red}{Let $X=\tau_{-s}^{*}\big((l+s)^{l}),l^{d-l-s}\big)$ and $Y=\tau_{-t}^{*}\big((l+t)^{l}),l^{d-l-t}\big)$, $s\geq t\geq0$ cases, the minimal $i \in \mathbb Z_{\geq 0}$ such that $\textnormal{Ext}_{\mathbb{H}_{n}}^{i}(X, Y)\neq0$ is $i=(t-s)l$. In this case, 
%    \[ \textnormal{Ext}_{\mathbb{H}_{n}}^{(s-t)l}(X, Y)\cong\mathbb{C}. \]}
%\end{proposition}
%\begin{proof}
%    The proof is similar to the proof of Proposition \ref{thm: Ext for two quasi-rectangula}
%\end{proof}

%\begin{remark}
%    \textcolor{red}{explain the difficulty for computing general $Ext^{i}$}
%\end{remark}

\part{Ext-groups between calibrated modules}
\section{Calibrated modules}\label{sec calibrated mod and skew young diagram}

\subsection{Semisimplicity of $W$-type modules}

The following lemma is rather well-known. We provide a brief explanation for the reader's convenience.

\begin{lemma} \label{lem w semisimple}
Let $X$ be a $W$-type $\mathbb H_n$-module. Then for any non-negative integers $n_1,\ldots, n_r$ with $n_1+\ldots+n_r=n$, $X|_{\mathbb H_{n_1}\otimes \ldots \otimes  \mathbb H_{n_r}}$ is a semisimple \(\mathbb H_{n_1}\otimes \ldots \otimes \mathbb H_{n_r}\)-module. 
\end{lemma}

\begin{proof}
Suppose $M=X|_{\mathbb H_{n_1}\otimes \ldots \otimes \mathbb H_{n_r}}$ is not semisimple. There exists an irreducible submodule $L$ of $M$ such that $i:L\hookrightarrow M$ is not split, hence we have
\[
0\rightarrow L\overset{i}{\hookrightarrow} M\rightarrow M/i(L)\rightarrow 0.
\]
Taking the dual, we obtain
\[
0 \rightarrow (M/i(L))^{\bullet}\rightarrow M^{\bullet} \overset{i^{\bullet}}{\twoheadrightarrow}L^{\bullet}\rightarrow 0.
\]
By Lemma \ref{lem dual of bullet}, there exists an isomorphism
\begin{align} \label{eqn dual restrict}
     \Phi:M^{\bullet}\overset{\cong}{\rightarrow}M.  
\end{align}
We define $q: i^{\bullet}\circ\Phi: M\twoheadrightarrow L^{\bullet}$, then we have
\[
L\overset{i}{\hookrightarrow}M\overset{q}{\twoheadrightarrow}L^{\bullet}.
\]
We claim that 
\[q\circ i=0.\]
Otherwise $q\circ i$ is an isomorphism since $L$ and $L^{\bullet}$ are irreducible, then $(q\circ i)^{-1}\circ q: M\twoheadrightarrow L$ is a surjection and $(q\circ i)^{-1}\circ q\circ i=id_{L}$, which contradicts the assumption that $i$ is not split. Hence we have
\[
0\subseteq i(L)\subseteq\operatorname{ker}q\subseteq M.
\]
Since $q$ is surjective, 
\[
M/\operatorname{ker}q\cong L^{\bullet}\cong L.
\]
Therefore, we know that the multiplicity of $L$ in $M$ is greater than $1$ by
Lemma \ref{lem irr bullet self dual}. However, this then contradicts  Proposition \ref{prop multiplicity one weight}.
\end{proof}

\subsection{A multiplicity free result}

\begin{lemma} \label{lem derivative} \cite{grojnowski2001strong}
Let $X$ be any irreducible $\mathbb H_n$-module. Let $y \in \mathbb Z$. Then there exists at most one irreducible $\mathbb H_{n-1}$-module $Y$ such that $\chi_y\boxtimes Y$ appears as a simple submodule of $X|_{\mathbb H_1\otimes \mathbb H_{n-1}}$.
\end{lemma}

\subsection{Calibrated modules}

We now generalize $W$-type modules to a larger family of $\mathbb H_n$-modules.

\begin{definition} \label{def calibrate}
 An irreducible $\mathbb H_n$-module $X$ is said to be \emph{calibrated} if there exists a $W$-type module $Y$ of $\mathbb H_m$ (for some $m \geq n$) and a $W$-type module $Y'$ of $\mathbb H_{m-n}$ such that 
 \[  Y' \boxtimes X
 \]
appears as an irreducible direct summand in $Y|_{\mathbb H_{m-n}\otimes \mathbb H_n}$.
\end{definition}
The notion of calibrated modules originates from \cite{krillo2002gha} while we use a different formulation from  \cite{krillo2002gha}. These modules are sometimes called ladder modules, introduced from different viewpoints \cite{lapidminguezladder} (also see \cite{barabschciubotaruladder}).

\begin{lemma} \label{lem restrict w type}
Let $\lambda$ be a Young diagram with $m$ boxes and let $z \in \mathbb C$. Let $X$ be a $W$-type module of $\mathbb H_{m-n}$, for some $n \leq m$. Then $X \boxtimes Y$ is an $\mathbb H_{m-n}\otimes \mathbb H_n$-submodule of $\tau_z^*(\lambda)$ for some calibrated $\mathbb H_{n}$-module $Y$ if and only if $X \cong \tau^*_z(\mu)$ for some subdiagram $\mu$  of $\lambda$.
\end{lemma}

\begin{proof}
The only if direction follows from that $\tau_z^*(\lambda)$, restricted as a $S_m$-representation, is irreducible and the Littlewood-Richardson rule. 

We now consider the if direction. Indeed, from the Littlewood-Richardson rule again, $\tau^*_z(\lambda)|_{S_{m-n}\times S_{n}}$ contains a subspace $\pi_{\mu}\boxtimes \tau$ for some irreducible $S_{n}$-representation $\tau$. Since the action of $e_1, \ldots, e_{m-n}$ is defined from the action of Jucys-Murphy elements, $\pi_{\mu}$ can be lifted to an action of $\mathbb H_{m-n}$, and so forms an $\mathbb H_{m-n}$-module isomorphic to $\tau_z^*(\mu)$. Now, one applies the $\mathbb H_{n}$-action to obtain a direct summand of the form $\tau^*_z(\mu)\boxtimes Y$ for some calibrated $\mathbb H_{n}$-module $Y$.
\end{proof}

\begin{lemma} \label{lem construct calibrated module}
Let $\lambda -\mu$ be a skew Young diagram. Suppose that $\lambda$ has $m+n$ boxes and $\mu$ has $m$ boxes. Let $z \in \mathbb C$. Let $Y=\tau^{*}_{z}(\lambda)$ and let $Y'=\tau^{*}_{z}(\mu)$. Then there exists a unique calibrated module $X$ of $\mathbb H_{n}$ such that $Y'\boxtimes X$ is a direct summand of 
\[     Y|_{\mathbb H_m \otimes \mathbb H_{n}} .
\]
\end{lemma}

\begin{proof}
We pick a weight $(y_1, \ldots, y_m)$ of $Y'$, and as $\mathbb H_1^{\otimes k}$-module, the corresponding character is of the form:
\[  \chi_{y_1}\boxtimes \ldots \boxtimes \chi_{y_m} . \]
Now it reduces to show that there exists a unique irreducible $\mathbb H_{n}$-module $X$ such that 
\[  \chi_{y_1}\boxtimes \ldots \boxtimes \chi_{y_m} \boxtimes X
\]
is a direct summand of $Y|_{\mathbb H_1^{\otimes m}\otimes \mathbb H_{n}}$.

We proceed with an induction on $m$. When $m=1$, the statement simply follows from Lemmas \ref{lem w semisimple} and \ref{lem derivative}. 

Now suppose $k \geq 2$. From the inductive hypothesis, we have a unique calibrated $\mathbb H_{n+1}$-module $X'$ such that $\chi_{y_1}\boxtimes \ldots \boxtimes \chi_{y_{m-}}\boxtimes X'$ is a direct summand of $Y|_{\mathbb H_{m-1}\otimes \mathbb H_{n+1}}$. Now the lemma follows from Lemma \ref{lem derivative}.
\end{proof}

    For each skew Young diagram $\lambda-\mu$ and $z \in \mathbb C$, we associate with the calibrated module in Lemma \ref{lem construct calibrated module}, denoted by $\beta_z(\lambda-\mu)$.

\begin{remark}
Let \( \lambda\) be a skew Young diagram of \(n\) boxes. Note that for $\nu \vdash n-i$, as \(S_{i}\times S_{n-i}\)-representation, \(\pi_{\lambda} \) decomposes as
\begin{equation} \label{eqn richardson coefficient}  \pi_{\lambda}|_{S_i \times S_{n-i}} = \bigoplus_{\mu \vdash i, \nu \vdash n-i} \mathrm{LR}^{\lambda}_{\mu, \nu} \pi_{\mu}\boxtimes\pi_{\nu}.
\end{equation}

Note that for a fixed \(\mu \vdash i\),  \( \pi_{\mu}\) part in the decomposition \ref{eqn richardson coefficient} can lift to the \(W\)-type \(\tau^*(\mu)\). Consequently, the uniqueness in Lemma \ref{lem construct calibrated module} forces that 
\[  \tau^*(\mu)\boxtimes \beta_z(\lambda-\mu).
\]
contribute all \(S_i \times S_{n-i}\)-representations of the form \(\pi_{\mu} \boxtimes \pi_{\nu}\) for some \(\nu \vdash n-i\) in the decomposition (\ref{eqn richardson coefficient}). This shows that 
 \[ \beta(\lambda-\mu)|_{S_{n-i}} =\bigoplus_{\nu \vdash n-i} \mathrm{LR}^{\lambda}_{\mu, \nu}\pi_{\nu} .\]
\end{remark}

\section{Uniform skew Young diagram}\label{sec uniform skew diagram}

\subsection{Some notions}

\begin{definition}
A skew Young diagram $\lambda-\mu$ is said to be {\it connected} if each row overlaps with the row just immediately it. 
\end{definition}

\begin{example} \label{ex skew young diagram}
Let $\lambda=(7, 4, 3)$ and $\mu=(2,1)$. The skew Young diagram $\lambda-\mu$ takes the form:
\[
\begin{ytableau}
      \none & \none    &        &    &   &     & \\
  \none     &          &        &     & \none  & \none   & \none    \\
           &           &      & \none   & \none  & \none  & \none
\end{ytableau}
\]
The first two boxes of the first row overlap with the last two boxes of the second row. The first two boxes of the second row overlap with the last two boxes of the third row. Thus $\lambda-\mu$ is connected.
\end{example}

\begin{example}
Let $\lambda=(7,4,3)$ and $\mu=(3,3)$. The skew Young diagram $\lambda-\mu$ takes the form:
\[
\begin{ytableau}
      \none & \none    &   \none     &    &   &     & \\
  \none     &   \none &   \none     &     & \none  & \none   & \none    \\
           &           &           & \none   & \none  & \none  & \none
\end{ytableau}
\]
The second row does not overlap with the third row. Thus, $\lambda-\mu$ is not connected. 
\end{example}

\begin{definition}\label{def skew hook} (cf. \cite[Section 8]{vershik2005new})
We say that a skew Young diagram \(\lambda-\mu\) is a \emph{skew hook} if each row overlaps with the row below it in exactly one box. The \emph{length} of a skew hook is the total number of boxes in the diagram.
\end{definition}

\begin{example}
Let $\lambda=(6,4,3,2)$ and $\mu=(3,2,1)$. The skew Young diagram $\lambda-\mu$ takes the form:
\[ \begin{ytableau}
      \none & \none    &   \none     &    &   &      \\
  \none     &   \none &        &     & \none  & \none     \\
   \none        &           &           & \none   & \none  & \none  \\
               &            &    \none    & \none   & \none  & \none  \\
\end{ytableau}
\]
and so is a skew hook.
\end{example}

\begin{example}
Let $\lambda=(6,4,4,4)$ and $\mu=(3,3,3)$. The skew Young diagram $\lambda-\mu$ takes the form:
\[  \begin{ytableau}
         \none & \none    &   \none     &    &   &      \\
  \none     &   \none &  \none      &     & \none  & \none     \\
   \none        &\none      &   \none       &    & \none  & \none  \\
            &       &         &   & \none  & \none  \\
\end{ytableau}
\]
and is a skew hook.
\end{example}

\begin{example}
Let $\lambda=(6,4)$ and $\mu=(2)$. The skew Young diagram $\lambda-\mu$ takes the form:
\[  \begin{ytableau}
         \none & \none    &       &    &   &      \\
                &     &      &     & \none  & \none    
\end{ytableau}
\]
The first row overlaps the second row in two boxes, and hence \(\lambda-\mu\) is not a skew hook.
\end{example}

\begin{definition} \label{def first hook}
\begin{enumerate}
\item The \emph{first skew hook} of a connected skew Young diagram \(\lambda-\mu\) is the largest subdiagram \(\psi\) of \(\lambda-\mu\) such that \(\psi\) contains the entire first row of \(\lambda-\mu\) and is a skew hook.
\item We say that \((\lambda-\mu)-\psi\) is \emph{hook-reduced} from \(\lambda-\mu\).
\end{enumerate}
\end{definition}

We remark that the last two criteria with the maximality uniquely determine the first hook.

\begin{example}
Let $\lambda=(7, 4, 3)$ and $\mu=(2,1)$. The skew Young diagram $\lambda-\mu$ takes the form:
\[
\begin{ytableau}
      \none & \none    &  *(blue!35) &  *(blue!35)   &  *(blue!35)  &   *(blue!35)   & *(blue!35) \\
  \none     &     *(blue!35)      &    *(blue!35)     &     & \none  & \none   & \none    \\
       *(blue!35)     &    *(blue!35)        &      & \none   & \none  & \none  & \none
\end{ytableau}
\] 
The blue-colored boxes form the first hook from $\lambda-\mu$.
\end{example}

\begin{definition}
The skew Young diagram $\lambda-\mu$ is said to be {\it strongly connected} if the following hold:
\begin{enumerate}
\item[(i)] $\lambda-\mu$ is connected; and
\item[(ii)] define $\zeta_1=\lambda-\mu$, and recursively for $i\geq 2$, define $\zeta_i$ to be hook-reduced from $\zeta_{i-1}$, and each $\zeta_i$ is required to be connected.
\end{enumerate}
Let $\eta_1, \eta_2, \ldots, \eta_d$ be the frist skew hooks from $\zeta_1, \ldots,
\zeta_d$ respectively. We shall call $\eta_1, \eta_2, \ldots, \eta_d=\zeta_d$ the {\it sequence of skew hooks} for $\lambda-\mu$, and call $d$ the \emph{hook sequence length} of $\lambda-\mu$.
\end{definition}

\begin{example}
Let $\lambda=(3,3,3,3)$. Then $(2,2,2)$ is hook-reduced from $(3,3,3,3)$ and $(1,1)$ is hook-reduced from $(2,2,2)$. Hence, the sequence of skew hooks for $\lambda$ is $(3,1,1,1)$, $(2,1,1), (1,1)$.
\end{example}

\begin{example}
Let $\lambda=(7,5,4)$ and $\mu=(2,1)$. The skew Young diagram takes the form:
\[
\begin{ytableau}
      \none & \none    &        &    &   &     & \\
  \none     &          &        &     &  & \none   & \none    \\
           &           &      &       & \none  & \none  & \none
\end{ytableau}
\]
The hook-reduced Young diagram from $\lambda - \mu$ takes the form:
\[  \begin{ytableau}
      \none  &  &     \\
             &  &  \none
\end{ytableau}
\]
which is a skew hook. Thus, $\lambda-\mu$ is strongly connected.
\end{example}

\begin{example}
The skew Young diagram in Example \ref{ex skew young diagram} is not strongly connected.
\end{example}

\begin{definition}
A strongly connected skew Young diagram $\lambda-\mu$ is said to be {\it uniform} if for the sequence of skew hooks $\eta_1,\ldots, \eta_d$ of $\lambda-\mu$,  $\eta_{i}$ has precisely two more boxes than $\eta_{i+1}$ for all $i \in \left\{1, \ldots, d-1\right\}$. 
\end{definition}

\begin{example}
Let $\lambda=(4,4,4,4)$ and $\mu=(2,2)$. Then $\lambda-\mu$ is strongly connected and uniform.
\end{example}

\begin{example}
Let $\lambda=(5,4,4)$ and $\mu=(1)$. Then $\lambda-\mu$ takes the form:
\[
\begin{ytableau}
      \none &       &        &    &     \\
            &          &        &     &  \none      \\
           &           &      &       & \none 
\end{ytableau}
\]
and is strongly connected. However, \(\lambda-\mu\) is not uniform since the second skew hook in the sequence of skew hooks for \(\lambda-\mu\) has three, rather than two, less boxes than the first skew hook in the sequence.
\end{example}

\subsection{A property for the $z$-content of a uniform skew Young diagram $\lambda-\mu$}

\begin{definition}
Let $\lambda-\mu$ be a skew Young diagram. For $x \in \mathbb C$, denote by $m(\lambda-\mu, z,x)$, the number of boxes in $\lambda-\mu$ whose $z$-contents are equal to $x$.
\end{definition}

\begin{lemma} \label{lem multiplicity increase}
Let $\lambda-\mu$ be a uniform skew Young diagram. Let $x$ be the smallest integer such that $m(\lambda-\mu, z, x)\neq 0$. Let $k$ be the hook sequence length of $\lambda-\mu$, and let $\ell$ be the length of the first hook. Then
\[  m(\lambda-\mu, z, y) = \left\{ \begin{array}{cc}  y-x+1  & \mbox{ if $x\leq y \leq x+k-1$} \\ k & \mbox{ if $x+k \leq  y \leq x+\ell-k$}   \\  x+\ell -y  & \mbox{ if $x+\ell-k+1\leq y\leq  x+\ell-1$} \\ 0  & \mbox{otherwise} \end{array}  \right.
\]
\end{lemma}
\begin{proof}
Let $\eta_1,\ldots, \eta_k$ be the sequence of skew hooks of $\lambda-\mu$.

The arrangement for the skew hooks in $\lambda-\mu$ forces that the \(z\)-content of the lowest and leftmost (resp. highest and rightmost) box in $\eta_i$ is strictly greater (resp. strictly smaller) than the $z$-content of the leftmost (resp. rightmost) box in $\eta_{i-1}$. 

On the other hand, by the definition of uniformity, for each $i=1,\ldots, k$, $\eta_i$ has exactly two fewer boxes than $\eta_{i-1}$. 

Combining the above two conditions, the \(z\)-content of the lowest and leftmost (resp. highest and rightmost) box of \(\eta_i\) is  the \(z\)-content of the lowest and leftmost box of \(\eta_{i-1}\) plus (resp. minus) one. 

The $z$-content of other boxes in $\eta_i$ can be uniquely determined. In particular, a simple counting then gives the number \(m(\lambda-\mu, z, x)\).
\end{proof}

\begin{example}
Let \(\lambda =(8,8,8,6)\) and \(\mu=(3,2,1)\). Then \(\lambda-\mu\) is a uniform skew Young diagram. Then the \(z\)-contents of \(\lambda-\mu\) are
\[
\begin{ytableau}
       \none    &  \none   & 2 & 3 &  4 & 5 & 6 \\
      \none     &  0    & 1 & 2 &  3 & 4 & 5 \\
     -2      & -1   & 0 & 1 &  2 & 3 & 4   \\
     -3      & -2   & -1 & 0 &  1 & \none & \none \\
     -4     & -3   & -2  & -1 & \none & \none & \none
\end{ytableau}
\]
The sequence of skew hooks of \(\lambda-\mu\) with their \(z\)-contents is given as follow:
\[  \eta_1= \begin{ytableau}
       \none    &  \none   & 2 & 3 &  4 & 5 & 6 \\
      \none     &  0    & 1 & \none & \none & \none & \none \\
     -2      & -1   & \none & \none & \none & \none & \none   \\
     -3      & \none   & \none & \none &  \none & \none & \none \\
     -4     & \none   & \none  & \none & \none & \none & \none
\end{ytableau},  \eta_2=\begin{ytableau}
       \none    &  \none   &  \none & \none &  \none & \none & \none\\
      \none     &  \none   & \none & 2 &  3 & 4 & 5 \\
     \none      & \none   & 0 & 1 &  \none & \none & \none   \\
     \none      & -2   & -1 & \none &  \none & \none & \none \\
     \none     & -3   & \none  & \none & \none & \none & \none
\end{ytableau} , \eta_3=
\begin{ytableau}
       \none    &  \none   & \none  & \none  &  \none  & \none  & \none  \\
      \none     &  \none     & \none  &\none  &  \none & \none  & \none  \\
     \none      & \none   & \none  & \none  &  2 & 3 & 4   \\
     \none       & \none    & \none  & 0 &  1 & \none & \none \\
     \none      & \none    & -2   & -1 & \none & \none & \none
\end{ytableau} .
\]
We observe that the numbers filled in \(\eta_1\) and \(\eta_2\) differ by the numbers in the lowest-leftmost and highest-rightmost boxes in \(\eta_1\), and similarly, the numbers filled in \(\eta_2\) and \(\eta_3\) differ by the numbers in the lowest-leftmost and highest-rightmost boxes in \(\eta_2\).
\end{example}

%To visualize Lemma \ref{lem multiplicity increase}, the $z$-content from each $\eta_i$ looks like:
%\[
%\begin{ytableau}
%    x    &      &        & \none[\cdot]   &  \none[\cdot]   & \none[\cdot]  & &  & x'\\
%     \none     &   x''       &        & \none[\cdot]    &   \none[\cdot]     &  \none[\cdot] &   & x''' & \none \\
%    \none       & \none   &       &  \none[\cdot]    &  \none[\cdot]    &  \none[\cdot] &   &  \none & \none 
%\end{ytableau}
%\]
%where $x'=x+\ell-1$, $x''=x+1$ and $x'''=x+\ell-2$. The first row represents the $z$-contents from the first hook $\eta_0$, the second row represents the $z$-content from the second hook $\eta_1$, and etc. Thus, one sees that the $m(\lambda, z, y)$ first increases, and then becomes steady in the middle, and finally decreases. 

\begin{lemma} \label{lem steady no increase}
Let $\lambda-\mu$ be a uniform skew Young diagram. Let $z \in \mathbb C$. There exists no $x \in \mathbb Z$ such that 
\[  m(\lambda-\mu, z, x-1) =m(\lambda-\mu, z, x) < m (\lambda-\mu, x+1) .
\]
\end{lemma}

\begin{proof}
This follows directly from Lemma \ref{lem multiplicity increase}.
\end{proof}

\section{Ext-non-vanishing conjecture for Arthur type modules}\label{sec ext non-vanishing conj for Arthur}

\subsection{Arthur type modules} \label{ss arthur type modules}

Roughly speaking, Arthur-type modules form a special class of unitary modules that are parabolically induced from Speh modules.

\begin{definition}
For each rectangular Young diagram $\lambda$, there exists a unique $z_0 \in \frac{1}{2}\mathbb Z$ such that 
\[ m(\lambda, z_0, y)=m(\lambda, z_0, -y) \]
for all $y \in \mathbb R$. We shall call these modules $\tau_{z_{0}}^{*}(\lambda)$  \emph{unitary Speh modules}. Here the unitarity coincides with the notion in \cite{barbasch2014unitary}, but we shall not need this.

Let $\mathrm{supp}(\tau^*_{z_0}(\lambda))$ be the set of $y \in \frac{1}{2}\mathbb Z$ such that $m(\lambda, z_0, y)\neq 0$. We emphasize that \(\mathrm{supp}(\tau^*_{z_0}(\lambda))\) is a set, rather than a multiset. 

Let \(X=\tau^*_{z_0}(\lambda)\) and \(X'=\tau^*_{z_0'}(\lambda')\) be two unitary Speh modules. We say that \(X\) and \(X'\) have the same parity if \(z_0-z_0'\in \mathbb Z\). Notice that if two unitary Speh modules \(X\) and \(X'\) have the same parity, then \(\mathrm{supp}(X)\subseteq \mathrm{supp}(X')\) or \(\mathrm{supp}(X')\subset \mathrm{supp}(X)\).

%This set encodes the information of the central character of $\tau^*_{z_0}(\lambda)$, but we do not count the multiplicities of $y$ appearing in the central character of $\tau_{z_0}^*(\lambda)$ for simplicity. 
\end{definition}

The following result is well-known, see for example \cite[Theorem 7.1]{tadic}.
\begin{lemma} \label{lem prop speh modules}
Let $X_1,\ldots, X_k$ be irreducible unitary Speh modules. Then $X_1\times \ldots \times X_k$ is irreducible. Moreover, for any permutation $\sigma \in S_k$,
\[   X_{\sigma(1)} \times \ldots \times X_{\sigma(k)} \cong X_1 \times \ldots \times X_k .
\]
\end{lemma}

\begin{definition}
An irreducible $\mathbb H_n$-module $X$ is said to be of \emph{Arthur type} if $X\cong X_1\times \ldots \times X_k$ for some unitary Speh modules $X_1, \ldots, X_k$.
\end{definition}

\subsection{Ext-non-vanishing conjecture}

%\begin{conjecture}
%Let $X$ and $Y$ be two Arthur-type $\mathbb H_n$-modules. Then $\mathrm{Ext}^i_{\mathbb H_n}(X, Y)\neq 0$ for some $i$ if and only if $X$ and $Y$ have the same central character.
%\end{conjecture}

\begin{lemma} \label{lem decompose }
Let $X$ be a Speh $\mathbb H_n$-module and let $Y$ be a Speh $\mathbb H_m$-module. Let $Z=\tau^*_{z_0}(\lambda)$ be a unitary Speh $\mathbb H_{n+m}$-module for some rectangular Young diagram $\lambda$. Suppose the following conditions hold:
\begin{enumerate}
\item[(i)] $\mathrm{supp}(X)=\mathrm{supp}(Z)$; and
\item[(ii)] $X\times Y$ has the same central character as $Z$.
\end{enumerate}
Then there exist unique $W$-type module $X'$ and Speh module $Y'$ such that 
\begin{enumerate}
\item  $X'$ is associated with a quasi-rectangular shape; and
\item $X$ has the same central character as $X'$; and
\item $Y$ has the same central character as $Y'$; and
\item $\mathrm{Ext}^i_{\mathbb H_{n+m}}(X\times Y, Z) \cong \mathrm{Ext}^i_{\mathbb H_n\otimes \mathbb H_m}(X\boxtimes Y,  X'\boxtimes Y')$.
\end{enumerate}
\end{lemma}

\begin{proof}
By Frobenius reciprocity in Lemma \ref{lem frobenius}, 
\[  \mathrm{Ext}^i_{\mathbb H_{n+m}}(X \times Y, Z) \cong \mathrm{Ext}^i_{\mathbb H_n \otimes \mathbb H_m}(X\boxtimes Y, Z|_{\mathbb H_n \otimes \mathbb H_m}) .
\
\]
On the other hand, by Lemma \ref{lem construct calibrated module},
\[ Z|_{\mathbb H_n \otimes \mathbb H_m} \cong \bigoplus_{\mu} \tau^*_z(\mu) \boxtimes \beta_z(\lambda-\mu) ,
\]
where $\mu$ runs for all Young subdaigrams of $\lambda$. 

By K\"unneth formula, 
\begin{equation}
\mathrm{Ext}^i_{\mathbb H_n \otimes \mathbb H_m}(X\boxtimes Y, Z|_{\mathbb H_n \otimes \mathbb H_m}) \cong \bigoplus_{i_1+i_2=i} \bigoplus_{\mu \vdash n} \mathrm{Ext}^{i_1}_{\mathbb H_n}(X, \tau^*_z(\mu)) \otimes \mathrm{Ext}^{i_2}_{\mathbb H_m}(Y, \beta_z(\lambda-\mu)) . 
\end{equation}
To show the lemma, it suffices to show that there is a unique subdiagram $\mu$ of \(\lambda\) such that $\tau_z^*(\mu)$ has the same central character as $X$.

We first show the existence part. Let \(\mu'\) be the underlying Young diagram associated to \(X\). Suppose the hook length of \(\mu'\) is \(k\). Let \(\mu_0\) be the subdiagram of \(\lambda\) formed by the first \(k\) hooks in the sequence of skew hook of \(\lambda\). Note that the numbers filled in the first skew hook of \(\mu_0\) coincide with the numbers filled in the first skew hook of \(\mu'\), by \(\mathrm{supp}(X)=\mathrm{supp}(Z)\). Since both \(\mu'\) and \(\mu_0\) are uniform and furthermore the hook lengths of \(\mu'\) and \(\mu_0\) are the same, the \(z\)-content of \(\mu'\) and the \(z\)-content of \(\mu_0\) are the same. Thus, \(\tau_z^*(\mu_0)\) and \(X\) have the same central character. Since \(X\times Y\) and \(Z\) have the same central character, \(Y\) and \(\beta_z(\lambda-\mu_0)\) have the same central character.

The uniqueness part is quite similar. Let \(\mu_0\) be the subdiagram of \(\lambda\) such that \(X\) and \(\tau_z^*(\mu_0)\) have the same central character. Then, by Proposition \ref{shape for same central character}, \(\mu_0=(l^p, d^q)\) for some \(l,d,p,q\). With the condition that $\mathrm{supp}(X)=\mathrm{supp}(Z)$, we know that $\mathrm{supp}(X)$ contains the largest and smallest values in $\mathrm{supp}(Z)$, and so $\mu_0$ must contain the entire first row of $\lambda$ and the entire first column of $\lambda$. In other words, $\mu_0$ contains the first hook of $\lambda$. Now, by the uniformity of $X$, one also deduces that the skew hooks in the sequence of skew hooks of $\mu_0$ form a part of the sequence of skew hooks of $\lambda$. This uniquely determines the shape of \(\mu_0\).
\end{proof}

We remark that the central character condition (i) is not significant. One can usually apply the commutativity result in Lemma \ref{lem prop speh modules} to rearrange the terms in the product to satisfy (i).

\begin{example} Let $\lambda$ be the Young diagram $(7,7,7)$ , let $\mu$ be the Young diagram $(3)$ and let $\zeta$ be the Young diagram $(6,6,6,6)$. The three Young diagrams determine three Speh modules, and their $z$-contents take the form respectively:
\[ \begin{ytableau}
       -2    &  -1   & 0  & 1 &  2 & 3 & 4  \\
      -3     &  -2    & -1 & 0 &  1 & 2 & 3  \\
     -4      & -3   &  -2 & -1 &  0 & 1 & 2   
\end{ytableau}, \quad
\begin{ytableau}
    -1 & 0 & 1 
\end{ytableau}
, \quad
\begin{ytableau}
          *(blue!35)-1  &  *(blue!35)0 &  *(blue!35)1 &   *(blue!35)2 &  *(blue!35)3 &  *(blue!35)4 \\
         *(blue!35)-2  &  *(blue!35)-1 &   *(blue!35)0 &   *(blue!35) 1 &  *(blue!35) 2 &  *(blue!35) 3 \\
       *(blue!35) -3  &   *(blue!35) -2 &  *(blue!35) -1 &   *(blue!35) 0 & *(blue!35) 1 &  *(blue!35)  2 \\
        *(blue!35) -4 &   *(blue!35) -3 &  *(blue!35) -2 & -1  & 0 & 1 
\end{ytableau}
\]
Let $\mu=(6,6,6,3)$ and $\mu'=(3)$. Lemma \ref{lem decompose } gives that
\begin{align*}
      & \mathrm{Ext}^i_{\mathbb H_{n+m}}(\tau^*_{-2}(\lambda) \times \tau^*_{-1}(\mu), \tau^*_{-1}(\zeta)) \\
     \cong & \mathrm{Ext}^i_{\mathbb H_n\otimes \mathbb H_m}(\tau^*_{-2}(\lambda)\boxtimes \tau^*_{-1}(\mu), 
      \tau^*_{-1}(\lambda')\boxtimes \tau^*_{-1}(\mu')) .
\end{align*} 
\end{example}

\begin{theorem} \label{thm conj speh}
Conjecture \ref{conj non vanishing Ext} holds when at least one of $X$ or $Y$ is a unitary Speh module.
\end{theorem}

\begin{proof}
We first consider the case where $Y$ is a unitary Speh module. By the definition of an Arthur type module, 
\[   X \cong X_1\times \cdots \times X_k
\]
for some unitary Speh modules $X_1, \cdots, X_k$. Since \(X\) and \(Z\) have the same central character, all \(x_j\) (\(j=1,\ldots k\) have the same parity. 

By Lemma \ref{lem prop speh modules}, we may rearrange $X_i$ such that for any $1\leq i <j \leq k$,
\[   \mathrm{supp}(X_j)  \subset \mathrm{supp}(X_i).
\]
Applying Lemma \ref{lem decompose } multiple times, we obtain
\begin{align*}
      \mathrm{Ext}^i(X_1\times \cdots \times X_k, Y) \cong  \mathrm{Ext}^i(X_1\boxtimes \cdots \boxtimes X_k, Y_1\boxtimes \ldots \boxtimes Y_k)
\end{align*}
for some $W$-type modules $Y_1, \ldots, Y_k$ with the same central character as $X_1, \ldots, X_k$ respectively. By K\"unneth formula, we further have that:
\begin{align} \label{eqn kunneth formula}
\mathrm{Ext}^i(X_1\boxtimes \cdots \boxtimes X_k, Y_1\boxtimes \cdots \boxtimes Y_k) \cong \bigoplus_{i_1+\ldots +i_k=i} \mathrm{Ext}^{i_1}(X_1,Y_1)\otimes \cdots \otimes \mathrm{Ext}^{i_k}(X_k, Y_k) 
\end{align}
By Proposition \ref{thm: Ext for two quasi-rectangula}, for each $j$, there exists $i_j^*$ such that $\mathrm{Ext}^{i_j^*}(X_j, Y_j)\neq 0$. Combining with (\ref{eqn kunneth formula}),
\[ \mathrm{Ext}^i(X_1\boxtimes \cdots \boxtimes X_k, Y_1\boxtimes \cdots \boxtimes Y_k)\neq 0\]
for $i=i_1^*+\ldots+i_k^*$, as desired.

We now consider the case where $X$ is a unitary Speh module. By \cite[Proposition 4.10]{chan2016duality},
\[  \mathrm{Ext}^i(X, Y) \cong \mathrm{Ext}^i(Y^*, X^*) ,
\]
where $*$ is an anti-involution defined as in \cite[Section 4.2]{chan2016duality}. It is a standard fact that the $*$ operation sends a Speh module to a Speh module, and sends an Arthur type module to an Arthur type module. Thus, now the above equality reduces to the previous case, and we are done.
\end{proof}

\section{Elementary transform}\label{sec elementary transform}

\subsection{Addable and removable corner}

%%%%%%%right addable%%
%\begin{definition} \label{def right addable}
%\textcolor{red}{delete right add and left rem}
%We say that $x \in \mathbb C$ is {\it right-addable} to the $z$-content of a skew Young diagram $\lambda-\mu$ if there exists a skew Young diagram $\lambda'-\mu$ such that 
%\begin{enumerate}
%\item $\lambda \subset \lambda'$ and $\lambda'$ has one more box than $\lambda$; and 
%\item the $z$-content of the extra box takes the value $x$.
%\end{enumerate}
% Denote the skew Young diagram $\lambda'-\mu$ by, $A_x^R(\lambda-\mu)$.
%\end{definition}

%\begin{example}
%Let $\lambda =(5,4,2,1)$ and $\mu=(2,2,2,1)$. Then $-4$, $-2$, $0$, $3$ and $5$ are right-addable to the $0$-content of $\lambda-\mu$.
%\end{example}
%%%%%%%%%%%%

\begin{definition}  \label{def left addable} 
 We say that $x \in \mathbb C$ is {\it left-addable} to the $z$-content of a skew Young diagram $\lambda-\mu$ if there exists a corner in $\mu$ whose $z$-content is $x$ and if, after removing this corner from \(\mu\) to obtain a diagram \(\mu'\), the numbers of rows and columns remain unchanged.  Denote  the skew Young diagram $\lambda-\mu'$ by $A^L_x(\lambda-\mu)$. 
%\begin{enumerate}
%\item $\mu'$ is a subdiagram of $\mu$; and 
%\item $\mu$ has one more box than $\mu'$; and 
%\item $x$ appears with precisely one more time in the $c$-content of $\mu$ than the $c$-content in $\mu'$.  
%\end{enumerate}
\end{definition}

\begin{example}
Let $\lambda=(5,4,4)$ and $\mu=(2,1)$. Then $1$ and $-1$ are left addable to the $0$-content $\lambda-\mu$, and
\[   A^L_{1}(\lambda-\mu) = (5,4,4)-(1,1), \quad A^L_{-1}(\lambda-\mu) =(5,4,4)-(2) .\]
\end{example}

\begin{example}
Let \( \lambda=(5,4,4)\) and \(\mu_1=(5,1)\). Then \(4\) is not left addable to the \(0\)-content of \(\lambda-\mu\) although there exists a corner in \(\mu\) whose \(0\)-content is \(4\).
\end{example}

%\begin{lemma}
%Let $\lambda-\mu$ be a skew Young diagram. Suppose, for any $x \in \mathbb Z$, $x$ is not left-addable to the $z$-content of $\lambda-\mu$. Then $\mu=\emptyset$ and so $\lambda-\mu$ is a Young diagram.
%\end{lemma}

The above lemma is immediate from definitions.

\begin{definition} \label{def right removable} 
We say that $x \in \mathbb C$ is \textit{right-removable} from the $z$-content of a skew Young diagram $\lambda-\mu$ if $\lambda$ has a corner that is not in $\mu$ and whose content is $x$. Denote by $R^R_x(\lambda-\mu)$ the skew Young diagram obtained by removing the corner from $\lambda-\mu$.

%if there exists a skew Young diagram $\lambda'-\mu$ such that
%\begin{enumerate}
%\item $\lambda'$ is a subdiagram of $\lambda$; and
%\item $\lambda$ has one more box than $\lambda'$; and 
%\item $x$ appears with precisely one more time in the $c$-content of $\lambda$ than the $c$-content in $\lambda'$.
%\end{enumerate}
\end{definition}

\begin{example}
Let $\lambda=(5,4,4, 3)$ and $\mu=(1)$. Then the values $5, 2, 0$ are right-removable from the $1$-content of $\lambda-\mu$. Then 
\[  R^R_5(\lambda-\mu)=(4,4,4,3)-(1), \quad R^R_3(\lambda-\mu)=(5,4,3,3) , \quad R^R_0(\lambda-\mu)=(5,4,4,2)=(1) .\]
\end{example}

%\begin{definition} \label{def left removable}
%We say that $x \in \mathbb C$ is {\it left-removable} from the $z$-content of a skew Young diagram $\lambda-\mu$ if there exists a skew Young diagram $\lambda-\mu'$ such that 
%\begin{enumerate}
%\item $\mu$ is a subdiagram of $\mu'$,  and $\mu'$  has one more box than $%\mu$; 
%\item the extra box in $\mu'$ has the $z$-content $x$.
%\end{enumerate}
%Denote the skew Young diagram  $\lambda-\mu'$ by $R^L_x(\lambda-\mu)$.
%\end{definition}

%\begin{example}
%Let $\lambda=(5,4,4,3)$ and $\mu=(1)$. Then $2$ and $0$ are left-removable from the $1$-content of $\lambda-\mu$. Then 
%\[  R^L_2(\lambda-\mu)=(5,4,4,3)-(2), \quad R^L_0(\lambda-\mu)=(5,4,4,3)-(1,1) . \]
% \end{example}

\begin{lemma}
The notation $A^R_x(\lambda-\mu), A^L_x(\lambda-\mu), R^R_x(\lambda-\mu), R^L_x(\lambda-\mu)$ is well-defined, that is, the skew Young diagrams in Definition  \ref{def left addable} and \ref{def right removable} obtained are unique. 
\end{lemma}

\begin{proof}
This follows from the observations that the corners of the $z$-content of $\lambda$ take different values, and similarly the corners of the $z$-content of $\lambda$ take different values.
\end{proof}

\begin{lemma} \label{lem addable property}
Let $\lambda-\mu$ be a skew Young diagram. If $y \in \mathbb C$ is left-addable to the $z$-content of $\lambda -\mu$, then 
\[  m(\lambda-\mu, z, y) \leq m(\lambda-\mu, z, y-1) .
\]
\end{lemma}

\begin{proof}
The condition that $y \in \mathbb C$ is left-addable to the $z$-content of $\lambda-\mu$ guarantees that if a row has a box with $z$-content $y-1$, then the row also has a box with $z$-content $y$. The lemma then follows.
\end{proof}

Similarly, we can prove:

%\begin{lemma} \label{lem removable property}
%Let $\lambda-\mu$ be a skew Young diagram. If $y \in \mathbb Z$ is right-removable to the $z$-content of $\lambda-\mu$, then 
%\[    m(\lambda-\mu, z, y-1) \leq m(\lambda-\mu, z, y) .\]
%\end{lemma}

\begin{lemma} \label{lem addable control}
In the setting of Lemma \ref{lem multiplicity increase}, for $y \in \mathbb C$ satisfying $x+l- k+1 \leq y \leq x+l$, $y$ is not  left-addable to the $z$-content of $\lambda-\mu$. 
\end{lemma}

\begin{proof}
Let $y \in \mathbb C$ satisfy $x+l-k+1 \leq y \leq x+l-1$. Suppose that $y$ is left-addable to the $z$-content of $\lambda-\mu$. Then, there exists a row which has a box with the $z$-content $y+1$, but does not have a box with the $z$-content $y$. 

On the other hand, there is at most one row that has a box with the $z$-content equal to $y$, but does not have a box with the $z$-content equal to $y+1$.

These imply that
\[    m(\lambda-\mu, z, y)\leq m(\lambda-\mu, z, y+1) .\]
This contradicts Lemma \ref{lem multiplicity increase}.
\end{proof}

We now show that the inequality in Lemma \ref{lem addable property} becomes an equality when $\lambda-\mu$ is uniform.

\begin{lemma} \label{cor addable control}
Suppose that $y \in \mathbb C$ is left-addable to the $z$-content of a uniform skew Young diagram $\lambda-\mu$. Then, \[   m(\lambda-\mu, z, y-1) =   m(\lambda-\mu, z, y) .  \]
\end{lemma}

\begin{proof}
By Lemma \ref{lem addable property}, we only have to show that
\[   m(\lambda-\mu, z, y-1) \leq m(\lambda-\mu, z, y) .
\]

It follows from Lemma \ref{lem addable control} that $y <l-k$ in the notation of Lemma  \ref{lem multiplicity increase}. Then the corollary also follows from Lemma \ref{lem addable control}.
\end{proof}

%\begin{lemma} \label{lem remove control}
%Let $\lambda-\mu$ be a skew Young diagram. Suppose $y$ is right removable from the $z$-content of $\lambda-\mu$. Then 
%\[   m(\lambda-\mu, z, y-1) \leq m(\lambda-\mu, z, y) .\]\end{lemma}

%\begin{proof}
% This follows from the observation that if a row has a box with the $z$-content $y-1$, then the rwo also has a box with the $z$-content $y$. Thus,
%\[   m(\lambda-\mu, z, y-1) \leq m(\lambda-\mu, z, y) .\]\end{proof}

%\begin{corollary}
%Let $\lambda-\mu$ be a uniform skew Young diagram. Suppose that $y \in \mathbb Z$ is left addable from the $z$-content of $\lambda-\mu$. Then, 
%\[  m(\lambda-\mu, z, y-1) = m(\lambda-\mu, z, y) .\]
%\end{corollary}

%\begin{proof}
%Note that in Lemma \ref{lem remove control}, if the strict inequality holds, then the row that has the rightmost box with the content $y$ has only one box. However, this contradicts \end{proof}

\begin{proposition} \label{prop addable removable}
Let $\lambda-\mu$ be a uniform skew Young diagram. If $y \in \mathbb C$ is left addable to the $z$-content of $\lambda -\mu$, then $y$ is right removable from  the $z$-content of $\lambda -\mu$. 
\end{proposition}

\begin{proof}
%By lemma \ref{lem addable property}, we have:
%\[   m(\lambda-\mu, z, y) \leq m(\lambda-\mu, z, y-1)  \] and any row with a $z$-content $y$ also has the $z$-content $y-1$. However,

By Lemma \ref{cor addable control} and $y$ being left addable, 
\begin{align}  \label{eqn add remove 1}
m(\lambda-\mu, z, y) = m(\lambda-\mu, z, y-1)  .
\end{align}

Suppose $y$ is not right removable from the $z$-content of $\lambda-\mu$. 
\begin{enumerate}
 \item[(a)] There exists a row whose right end point has the $z$-content $y$. As $y$ is not right removable, there must also be a row whose leftmost box has the $z$-content $y-1$, and thus the row does not have a box with the $z$-content $y$. Combining with the first paragraph, this implies that
 \begin{align*}
 m(\lambda-\mu, z, y) < m(\lambda-\mu, z, y-1) ,
 \end{align*}
but this is also impossible from the above discussion.
 \item[(b)] There does not exist a row whose right end point has the $z$-content $y$. Then any row with the $z$-content $y$ also has the $z$-content $y+1$. On the other hand, there exists a row  with the $z$-content $y+1$, but not with the $z$-content $y$ (the row that contributes the addability of $y$). Thus, we have that:
 \begin{align} \label{eqn add remove 2}
 m(\lambda-\mu, z, y) < m(\lambda-\mu, z, y+1) . 
 \end{align}
Combining  (\ref{eqn add remove 1}) and (\ref{eqn add remove 2}), we have:
\[  m(\lambda-\mu, z, y-1) = m(\lambda-\mu, z, y) < m(\lambda-\mu, z, y+1) .
\]
This contradicts Lemma \ref{lem steady no increase}.
\end{enumerate}
As the above two cases exhaust all possibilities, $y$ is right removable from the $z$-content of $\lambda-\mu$.
\end{proof}

\subsection{Elementary transform}

\begin{definition} \label{def transform 1}
Let $\lambda-\mu$ be a skew Young diagram. Suppose $y \in \mathbb C$ is right removable from the $z$-content of $\lambda-\mu$. By Proposition \ref{prop addable removable}, $y$ is also left addable to the $z$-content of $\lambda-\mu$. Then $y$ is still right removable from the $z$-content of $A^L_y(\lambda-\mu)$. Thus $R^R_y\circ A^L_y(\lambda-\mu)$ is well-defined. Note that \(A\)

We say that $R^R_y\circ A^L_y(\lambda-\mu)$ is the \emph{elementary transform} of $\lambda-\mu$ with respect to $y$. 
\end{definition}

\begin{lemma} \label{lem conn uniform transform}
We use the notation in Definition \ref{def transform 1}. Then $ R^R_y \circ A^L_y(\lambda-\mu)$ is still uniform.
\end{lemma}

\begin{proof}
Note that the multiplicities of each $y \in \mathbb C$ in the $z$-contents of $\lambda-\mu$ and $R^R_y\circ A^L_y(\lambda-\mu)$ are the same. Then the connectedness of $\lambda-\mu$ implies the connectedness of $\lambda'-\mu'$, and the uniformity of $\lambda-\mu$ also implies the uniformity of $R^R_y\circ A^L_y(\lambda-\mu)$.
\end{proof}

%\begin{definition} \label{def sequence of transform}
%For two uniform connected skew Young diagrams $\lambda-\mu$ and $\lambda-\mu'$, we write 
%\[  \lambda - \mu \Dashv \lambda'-\mu'\]
%if there exists a sequence of $y_1,\ldots, y_k \in \mathbb Z$ such that $\lambda'-\mu'$ is obtained from $\lambda-\mu$ by a sequence of elementary transforms with respect to $y_1,\ldots y_k$, that is,
%\[ \lambda'-\mu' =T_{y_r}\circ \cdots \circ T_{y_1}(\lambda-\mu) .\]
%\end{definition}

%\begin{lemma} \label{lem transform inverse}
%Let $\lambda-\mu$ be a uniform connected skew Young diagram. 
%\begin{enumerate} 
%\item Suppose $y$ is left-addable to the $z$-content of $\lambda-\mu$. Let $\lambda'-\mu'$ be in Definition \ref{def transform 2} i.e.
%\[  \lambda-\mu \vdash^y \lambda'-\mu'\]
%Then $\lambda'-\mu' \dashv^y \lambda-\mu$.
%\item Suppose $y$ is right-removable to the $z$-content of $\lambda-\mu$. Let $\lambda'-\mu'$ be in Definition \ref{def transform 1} i.e.
%\[   \lambda-\mu  \dashv^y \lambda'-\mu' .\]
%Then $\lambda'-\mu' \vdash^y \lambda-\mu$.
%\end{enumerate} \end{lemma}

\subsection{Transform to a Young diagram}

\begin{lemma} \label{lem young diag connected}
Any rectangular Young diagram is uniform.
\end{lemma}

The above lemma follows from a simple counting on the boxes in the skew hooks in the sequence of skew hooks of a rectangular Young diagram. We omit a proof.

\begin{corollary} \label{cor skew yound to young diagram}
Let $\lambda-\mu$ be a uniform skew Young diagram. Then there exists a Young diagram $\gamma$ such that 
\[  \gamma = R^R_{y_k}\circ A^L_{y_k}\circ \ldots  \circ R^R_{y_1} \circ  A^L_{y_1}(\lambda-\mu) .
\]
for some $y_1, \ldots, y_k \in \mathbb C$.
\end{corollary}

\begin{proof}
If $\lambda-\mu$ is a Young diagram, that is $\mu=\emptyset$, then there is nothing to prove. We now assume that $\lambda-\mu$ is not a rectangular Young diagram. Then there exists some $y \in \mathbb Z$ such that $y$ is left-addable to the $z$-content of $\lambda-\mu$. Then one obtains a skew Young diagram $R^R_y\circ A_y^L(\lambda-\mu)$.

By Lemma \ref{lem conn uniform transform}, one repeats the process, and it terminates when one obtains a Young diagram $\gamma$. This shows the corollary.
\end{proof}

\section{Uniform calibrated $\mathbb H_n$-modules}\label{sec uniform calibrated modules}

\begin{definition}
 A calibrated $\mathbb H_n$-module $\tau^*_z(\lambda-\mu)$ is {\it uniform} (resp. \emph{strongly connected})  if $\lambda-\mu$ is uniform (resp. strongly connected).
\end{definition}

The important feature of a uniform calibrated \(\mathbb H_n\)-module \(\tau^*_z(\lambda-\mu)\) is as follows.

\begin{proposition}
A calibrated $\mathbb H_n$-module $\omega$ is uniform if and only if $\omega$ has the same central character as some Speh module of $\mathbb H_n$.
\end{proposition}

\begin{proof}
The proof of the "only if" direction is a simple comparison of central character. We only provide an example to illustrate the main idea. We consider the skew Young diagram with the $0$-content: 
\[
\begin{ytableau}
       \none    &  \none   & 2 & 3 &  4 & 5 & 6 \\
      \none     &  0    & 1 & 2 &  3 & 4 & 5 \\
     -2      & -1   & 0 & 1 &  2 & 3 & 4   \\
     -3      & -2   & -1 & 0 &  1 & \none & \none \\
     -4     & -3   & -2  & -1 & \none & \none & \none
\end{ytableau}
\]
The sequence of skew hooks for $\lambda-\mu$ is:
\[ \begin{ytableau}
       \none    &  \none   & 2 & 3 &  4 & 5 & 6 \\
      \none     &  0    & 1 & \none & \none & \none & \none \\
     -2      & -1   & \none & \none & \none & \none & \none   \\
     -3      & \none   & \none & \none &  \none & \none & \none \\
     -4     & \none   & \none  & \none & \none & \none & \none
\end{ytableau}
\begin{ytableau}
       \none    &  \none   &  \none & \none &  \none & \none & \none\\
      \none     &  \none   & \none & 2 &  3 & 4 & 5 \\
     \none      & \none   & 0 & 1 &  \none & \none & \none   \\
     \none      & -2   & -1 & \none &  \none & \none & \none \\
     \none     & -3   & \none  & \none & \none & \none & \none
\end{ytableau}
\begin{ytableau}
       \none    &  \none   & \none  & \none  &  \none  & \none  & \none  \\
      \none     &  \none     & \none  &\none  &  \none & \none  & \none  \\
     \none      & \none   & \none  & \none  &  2 & 3 & 4   \\
     \none       & \none    & \none  & 0 &  1 & \none & \none \\
     \none      & \none    & -2   & -1 & \none & \none & \none
\end{ytableau}
\]
One can now unfold the skew hooks to obtain the corresponding Speh module determined by the following Young diagram with the $-2$-content by Proposition \ref{shape for same central character}:
\[ \begin{ytableau}
       -2    &  -1   & 0  & 1 &  2 & 3 & 4 & 5 & 6 \\
      -3     &  -2    & -1 & 0 &  1 & 2 & 3 & 4 & 5 \\
     -4      & -3   &  -2 & -1 &  0 & 1 & 2 & 3 & 4  
\end{ytableau}
\]
We remark that the number of rows is determined by the hook length of $\lambda-\mu$ (i.e. the number of terms in the sequence of the skew hook).

For the "if" direction, we have observed in Lemma \ref{lem young diag connected} that a rectangular Young diagram is uniform. The connectedness of the rectangular Young diagram then implies that $\lambda-\mu$ is also connected, and its uniformity implies that $\lambda-\mu$ is also uniform. 
\end{proof}

\subsection{Connections to removable points}

\begin{lemma}
Let $\lambda$ be a Young diagram with $n$ boxes and let $z \in \mathbb C$. Suppose $y \in \mathbb C$ is right-removable from the $z$-content of $\lambda$. Let $\mu=R^R_y(\lambda)$, that is the Young diagram obtained by removing the box with the right-removable value $y$ from $\lambda$. Then $\tau^*_z(\mu) \boxtimes \chi_y$ is a direct summand of $\tau^*_z(\lambda)|_{\mathbb H_{n-1}\otimes \mathbb H_1}$.
\end{lemma}

\begin{proof}
By the Littlewood-Richardson rule, $\pi_{\mu}$ is a submodule of $\pi_{\lambda}|_{S_{n-1}\times S_1}$. Since $\tau^*_z(\lambda)$ and $\tau^*_y(\mu)$ are constructed from Jucys-Murphys elements, $\tau^*_z(\mu)$ is an $\mathbb H_{n-1}$-submodule of $\tau^*_z(\lambda)$. By extending $\tau^*_z(\mu)\boxtimes\operatorname{triv}$ to an $\mathbb H_{n-1}\otimes \mathbb H_1$-module, and comparing their central characters, the semisimplicity of $\tau^*_z(\lambda)|_{\mathbb H_{n-1}\otimes \mathbb H_1}$ from  Lemma \ref{lem w semisimple} implies that $\tau^*_z(\mu)\boxtimes \chi_{y}$ is a direct summand of $\tau^*_z(\lambda)|_{\mathbb H_{n-1}\otimes \mathbb H_1}$.
\end{proof}

\begin{corollary} \label{cor removable restriction}
Let $\lambda-\mu$ be a skew Young diagram of $n$ boxes.
%and let $z \in \mathbb C$.  Let $y \in \mathbb C$ be right removable from the $z$-content of $\lambda-\mu$. Then $\beta_z(R_y^R(\lambda-\mu))\boxtimes \chi_y$ is a direct summand in $\beta_z(\lambda-\mu)|_{\mathbb H_{n-1}\boxtimes \mathbb H_1}$.
\begin{enumerate}
 \item Let $y \in \mathbb C$ be right removable from the $z$-content of $\lambda-\mu$. Then $\beta_z(R_y^R(\lambda-\mu))\boxtimes \chi_y$ is a direct summand in $\beta_z(\lambda-\mu)|_{\mathbb H_{n-1}\boxtimes \mathbb H_1}$.
\item Let $y \in \mathbb C$ be left removable from the $z$-content of $\lambda-\mu$. Then $\chi_y \boxtimes \beta_z(R^L_y(\lambda-\mu))$ is a direct summand in $\beta_z(\lambda-\mu)|_{\mathbb H_1\boxtimes \mathbb H_{n-1}}$.
\end{enumerate}
\end{corollary}

\begin{proof}
Suppose that $\mu$ has $m$ boxes. Hence, $\lambda$ has $m+n$ boxes. By Lemma \ref{lem construct calibrated module}, and the right removability of the value \(y\),
\[  \tau^*_z(\lambda)|_{\mathbb H_{m+n-1} \otimes \mathbb H_1} 
\]
has a direct summand of $\beta_z(\lambda')\boxtimes \chi_y$ for some subdiagram \(\lambda'\) of \(\lambda\) such that \(\mu \subset \lambda'\). Since $\mu \subset \lambda'$, we then obtain that
\[  \tau^*_z(\mu) \boxtimes \beta_z(\lambda'-\mu) \boxtimes \chi_y
\]
is a direct summand of $\beta_z(\lambda'-\mu)|_{\mathbb H_m\otimes \mathbb H_{n-1}\otimes \mathbb H_1}$. Moreover, it is the unique summand of the form $\tau^*_z(\mu) \boxtimes X \boxtimes \chi_y$ for some irreducible $\mathbb H_{n-1}$-module \(X\).

It follows from definitions that $\tau^*_{z}(\mu) \boxtimes \beta_z(\lambda-\mu)$ is a direct summand of $\tau^*_z(\lambda)|_{\mathbb H_{m}\otimes \mathbb H_n}$. By the uniqueness in Lemma \ref{lem construct calibrated module},  
\[  \beta_z(\lambda-\mu)|_{\mathbb H_{n-1}\otimes \mathbb H_1}
\]
contains a direct summand of the form $X \boxtimes \chi_y$ for some irreducible $\mathbb H_{m+n-1}$-module $X$ so that $\beta_z(\lambda)|_{\mathbb H_m \otimes \mathbb H_{n-1}\otimes \mathbb H_1}$ contains the direct summand of the form $\tau^*_z(\mu)\boxtimes X' \boxtimes \chi_y$ for some irreducible $\mathbb H_{n-1}$-module $X'$.

Combining with the above uniqueness, we have that:
\[  X \cong X' \cong \beta_z(R^R_y(\lambda-\mu)) .  
\]
This proves (1) of the corollary. 

A proof for (2) is quite similar to that for (1) by using restrictions in stages.
\end{proof}

\section{Extension between Speh modules and calibrated modules}\label{sec ext between speh and calibrated}

%\subsection{Calibrated/ladder representation}

%\begin{proposition} \label{prop calibrated summand from W}
%Let $\pi$ be a ladder $\mathbb H_n$-module. Then there exists a $W$-lifted module $\omega$ of $\mathbb H_n$ for some $n$ such that $\omega'\otimes \pi$ is a direct summand in
%\[   \omega|_{\mathbb H_{m-n} \otimes \mathbb H_m} .
%\]
%\end{proposition}

%\begin{remark}
%If we fix the value $m$ such that the $W$-lifted $\mathbb H_n$-module $\omega$ exists, then we have the uniqueness of such $\omega$.
%\end{remark}

\subsection{Transfer theorem}

We first recall the following commutativity of parabolic induction. The lemma is quite standard; see, for example, \cite[Lemma 9.1]{chanhomobranch}, and we omit a proof. 

\begin{lemma}  \label{lem commute speh}
Let $\lambda$ be a rectangular Young diagram. Let $z\in \mathbb C$. Let $y \in \mathbb C$ such that $m(\lambda, z, y)\neq 0$. Then 
\[      \tau^*_z(\lambda) \times \chi_y \cong \chi_y \times \tau_z^*(\lambda) . \]
\end{lemma}

\begin{theorem} \label{thm transform one step calibrate}
Let $\lambda-\mu$ be a skew Young diagram with $n$ boxes. Let $y\in \mathbb C$ be left-addable to the $z$-content of $\lambda-\mu$. Let $X$ be a unitary Speh module of $\mathbb H_n$. Then 
\[  \dim \mathrm{Ext}^i_{\mathbb H_n}(X, \beta_z(\lambda-\mu)) = \dim \mathrm{Ext}^i_{\mathbb H_n}(X, \beta_z(R^R_y\circ A^L_y(\lambda-\mu))).
\]
\end{theorem}

\begin{proof}
Let $\omega$ be the rectangular Young diagram and let $z \in \mathbb C$ such that $X\cong \tau_z^*(\omega)$. By the definition of addability of $y$, $m(\omega, z, y)\neq 0$.

By Frobenius reciprocity and Corollary \ref{cor removable restriction}, for all $i$,
\[  \mathrm{Ext}^i_{\mathbb H_{n+1}}(\chi_y \times X, \beta_z(A^L_y(\lambda-\mu))) \cong  \mathrm{Ext}^i_{\mathbb H_{1}\otimes \mathbb H_{n}}(\chi_y \boxtimes X, \beta_z(A^L_y(\lambda-\mu)|_{\mathbb H_{1}\boxtimes \mathbb H_n}) . \]
On the other hand, by Corollary \ref{cor removable restriction}, we have that $\chi_y\boxtimes \beta_z(\lambda-\mu)$ is the unique direct summand in
\[  \beta_z(A_y^L(\lambda-\mu))|_{\mathbb H_1\otimes \mathbb H_{n}}
\]
of the form $\chi_y \boxtimes X$ for some $\mathbb H_{n}$-module $X$. Thus, by Proposition \ref{prop cc vanish},
\[  \mathrm{Ext}^i_{\mathbb H_{n+1}}( \chi_y \times X, \beta_z(A^L_y(\lambda-\mu))) \cong \mathrm{Ext}^i_{\mathbb H_1 \otimes \mathbb H_{n}}(\chi_y\boxtimes X, \chi_y\boxtimes \beta_z(\lambda-\mu)|_{\mathbb{H}_1\boxtimes\mathbb{H}_n}) .
\]
By the Ku\"nneth formula,
\[  \mathrm{Ext}^i_{\mathbb H_1\otimes \mathbb H_n}(\chi_y \boxtimes X, \chi_y\boxtimes\beta_z(\lambda-\mu))\cong \bigoplus_{i_1+i_2=i} \mathrm{Ext}^{i_1}_{\mathbb H_1}(\chi_y, \chi_y)\otimes\mathrm{Ext}^{i_2}_{\mathbb H_{n}}(X, \beta_z(\lambda-\mu))  .
\]
Combining with 
\[ \mathrm{Ext}^i_{\mathbb H_1}(\chi_y, \chi_y)\cong \begin{cases} \mathbb C & \mbox{ if \(i=0,1\)} \\ 0 & \mbox{ otherwise } \end{cases} ,
\]
we obtain
\begin{align} \label{eqn ext dimension}
\dim \mathrm{Ext}^{i}_{\mathbb H_{n}}(X, \beta_z(\lambda-\mu))=\sum_{t=0}^{i}(-1)^{t}\dim \mathrm{Ext}^{i-t}_{\mathbb H_{n+1}}(\chi_y \times X,  \beta_z(A^L_y(\lambda-\mu))).
\end{align}

Similarly, 
\[  \mathrm{Ext}^i_{\mathbb H_{n+1}}(X\times \chi_y, \beta_z(A^L_y(\lambda-\mu))) \cong \mathrm{Ext}^i_{\mathbb H_{n} \otimes \mathbb H_1}( X \boxtimes \chi_y,  \beta_z(R^R_y\circ A^L_y(\lambda-\mu))\boxtimes \chi_y ) .
\]
By a similar argument to the above one, we deduce 
\begin{align} \label{eqn ext dimension 2}
\dim \mathrm{Ext}^i_{\mathbb H_n}(X, \beta_z(R^R_y\circ A^L_y(\lambda-\mu))) =\sum_{t=0}^{i}(-1)^{t}\dim \mathrm{Ext}^{i-t}_{\mathbb H_{n+1}}( X \times \chi_y,  \beta_z(A^L_y(\lambda-\mu))). 
\end{align}
Combining (\ref{eqn ext dimension}) and (\ref{eqn ext dimension 2}) with Lemma \ref{lem commute speh}, we obtain the lemma.
\end{proof}

\begin{corollary} \label{cor transfer thm}
Let $X$ be a unitary Speh $\mathbb H_n$-module. Let $Y$ be a calibrated $\mathbb H_n$-module with the same central character as $X$. Then there exists a unique $W$-type module $Z$ such that 
\begin{enumerate}
\item[(i)] for all $i$,
\[  \mathrm{Ext}^i_{\mathbb H_n}(X, Y) \cong \mathrm{Ext}^i_{\mathbb H_n}(X, Z) ;
\]
\item[(ii)] the Young diagram associated to $Z$ has the same number of rows as the skew-Young diagram associated to $Y$.
\end{enumerate}

\end{corollary}

\begin{proof}
The first condition follows from repeated applications of Theorem \ref{thm transform one step calibrate} and Corollary \ref{cor skew yound to young diagram}. 

We now explain the second condition. It \(X\) is a unitary Speh \(\mathbb H_1\)-module, then the corollary is straightforward. Thus, we may and shall assume that \(n \geq 2\). By induction, it suffices to see that for each left-addable value \(y \in \mathbb C\) to \(\lambda-\mu\), \(R^R_y(\lambda-\mu)\) has the same number of rows as \(\lambda-\mu\). To this end, suppose the rightmost box $\Box$ in a row of a skew Young diagram $\lambda-\mu$ is removed so that the resulting Young diagram $\lambda'-\mu$ is still a Young diagram. If the resulting Young diagram has fewer rows, then the corresponding row must have consisted only of the box $\Box$, and so it must be in the last row of \(\lambda-\mu\). However, the connectedness condition forces that the row above also contains a box in the same column as \(\Box\). The $z$-content of the $\Box$ is the least value among all numbers in the $z$-content of $\lambda'-\mu$, and hence cannot be left-addable to the  $z$-content of $\lambda'-\mu$. Thus, \(R^R_y(\lambda-\mu)\) has the same number of rows as \(\lambda-\mu\). It follows from Definition \ref{def left addable} that the number of rows fo \(A^L_y\circ R^R_y(\lambda-\mu)\) is the same as that of \(R^R_y(\lambda-\mu)\). This shows the second condition. 

The uniqueness is determined by the central character of $Z$ and the number of rows of the Young diagram associated to $Z$. 
\end{proof}

%%%%%%%%%%%

%\subsection{Proof of transfer theorem}

%\begin{proof}
%Let $\sigma$ and $\pi$ be as in Theorem \ref{thm transfer thm}. By Corollary \ref{cor ladder transform}, there exists a lifted $W$-type multisegment $\widetilde{\mathfrak m}$ such that 
%\[   \mathfrak m \vdash^{a_1} \mathfrak m_1 \vdash^{a_2} \ldots \vdash^{a_k} \mathfrak m_k =\widetilde{\mathfrak m}  \]
%for some $a_1, a_2, \ldots, a_k$ lying the the central character of $\sigma$.

%In each of the above steps, we only add the nodes without removing the nodes. As a result, we have another segment, denoted by $\overline{\mathfrak m}$.

%Then
%\[  \mathrm{Ext}^i(\sigma \times \mathrm{St}(a_1)\times \ldots \times \mathrm{St}(a_k), \mathrm{St}(\overline{\mathfrak m}) ) \cong \mathrm{Ext}^i(\sigma \boxtimes 
%\mathrm{St}([a_1])\boxtimes \ldots \boxtimes \mathrm{St}([a_k]), \pi \boxtimes \mathrm{St}([a_1])\boxtimes \ldots \boxtimes \mathrm{St}([a_k])\]
%and
%\[  \mathrm{Ext}^i(\mathrm{St}([a_1])\times \ldots \times \mathrm{St}([a_k]) \times \sigma, \mathrm{St}(\overline{\mathfrak m})) \cong \mathrm{Ext}^i(\mathrm{St}([a_1]) \boxtimes \ldots \boxtimes \mathrm{St}([a_k]) \boxtimes \sigma, \mathrm{St}([a_1]) \boxtimes \ldots \boxtimes \mathrm{St}([a_k]) \boxtimes\mathrm{St}(\widetilde{\mathfrak m})) .\]

%By using the K\"unneth formula, we reocver and inductively have:
%\[ \mathrm{Ext}^i(\sigma, \pi) \cong \mathrm{Ext}^i(\sigma, \mathrm{St}(\widetilde{\mathfrak m})) .\]
%\end{proof}

\subsection{A combinatorial description of a result of Orlik}

We explain a combinatorial description of a result \cite{orlik2005extensions} for type $A_n$. For each skew hook $\lambda-\mu$ with $n$ boxes, choose $z$ to be the number such that $1$ is the least value among all $z$-contents of boxes in $\lambda-\mu$. Then the $z$-contents of $\lambda-\mu$ are the numbers $1, \ldots, n$ with multiplicity one. Suppose that $\lambda-\mu$ has $k$ rows. Then we divide these $n$ numbers into $k$-groups according to the $z$-contents of boxes in each row. If the numbers $i$ and $i+1$ are not in the same row, then we say that $i$ and $i+1$ are divided. 

 We denote the corresponding division by $\mathcal D_{\lambda-\mu}$. Denote by \(|\mathcal D_{\lambda-\mu}|\) the number \(n\) minus the number of parts in the division. Let $\lambda_1=(1^n)$, and the division $\mathcal D_{(1^n)}$ takes the form:
\[  1\ |\ 2\ | \ldots |\ n
\]
i.e. each group on the division contains only one number and hence \(|\mathcal D_{(1^n)}|=0\). Let $\lambda_2=(n)$. Then $\mathcal D_{(n)}$ corresponds to the division with only one part  and hence \(|\mathcal D_{(n)}|=n-1\). For \(i=1,\ldots, n-1\), the pairs \(i\) and \(i+1\) are not divided.  Let \(\lambda =(3,3,2)\). Then \(\mathcal D_{(6,4,2)-(3,1)}\) takes the form:
\[    |1\ 2  | 3\ 4\ 5 | 6\ 7\ 8| .
\]
Then \(|\mathcal D_{(6,4,2)-(3,1)}|=5\). Moreover, the numbers \(2\) and \(3\) are divided, and the numbers \(5\) and \(6\) are divided.

For two divisions $\mathcal D_1$ and $\mathcal D_2$, we define the union of the divisions $\mathcal D_1$ and $\mathcal D_2$, denoted by $\mathcal D_1 \cup \mathcal D_2$, such that the numbers $i$ and $i+1$ are divided in $\mathcal D_1\cup \mathcal D_2$ if $i$ and $i+1$ are divided in both $\mathcal D_1$ and $\mathcal D_2$.

For two divisions $\mathcal D_1$ and $\mathcal D_2$, we define the intersection of the divisions $\mathcal D_1$ and $\mathcal D_2$, denoted by $\mathcal D_1 \cap \mathcal D_2$, such that the numbers $i$ and $i+1$ are divided if and only if $i$ and $i+1$ are divided in at least one of $\mathcal D_1$ and $\mathcal D_2$.

\begin{theorem} \cite{orlik2005extensions} \label{thm combin generalized steinberg}
Let $X=\tau^*_z(\lambda-\mu)$ and let $X'=\tau^*_{z'}(\lambda'-\mu')$ for some skew hooks $\lambda-\mu$ and $\lambda'-\mu'$. Let $\mathcal D=\mathcal D_{\lambda-\mu}$ and $\mathcal D'=\mathcal D_{\lambda'-\mu'}$. Suppose $X$ and $X'$ have the same central character. Then 
\[ \mathrm{Ext}^i(X, X') \cong \begin{cases} \mathbb C & \mbox{ if $i=|\mathcal D\cup \mathcal D'|-|\mathcal D\cap \mathcal D'|, |\mathcal D \cup \mathcal D'|-|\mathcal D\cap \mathcal D'|+1$} \\
0 & \mbox{ otherwise } \end{cases}.
\]
\end{theorem}

\begin{proof}
This theorem is a reformulation of \cite[Theorem 1]{orlik2005extensions}, and we only sketch the main idea.  Let \(\Delta=\left\{e_1-e_2, e_2-e_3, \ldots, e_{n-1}-e_n\right\}\) be the set of simple roots associated to the root system of type \(A_{n-1}\). Given a division \(\mathcal D_{\lambda-\mu}\), we associate with a subset \(I\) of \(\Delta\) such that \(e_i-e_{i+1} \in I\) if and only if \(i\) and  \(i+1\) are not divided in \(\mathcal D\). The subset \(I\) determines the corresponding parameter of \(\tau^*_z(\lambda-\mu)\) in the sense of Orlik \cite{orlik2005extensions}. Note that, given two divisions \(\mathcal D\) and \(\mathcal D'\), if \(I\) and \(I'\) are the corresponding subsets of \(\Delta\) described above, then \(|I\cup I'|=|\mathcal D\cup \mathcal D'|\) and \(|I \cap I'|=|\mathcal D\cap \mathcal D'|\). The theorem is then translated from \cite[Theorem 1]{orlik2005extensions}.
\end{proof}

\begin{example} \label{ex one hook case}
Let $\lambda=(n)$. Let $\gamma$ be a skew hook with $n$ boxes and with $k$ rows. Let $\gamma'$ be the hook with $n$ boxes and with $k$ rows. Let $z, z' \in \frac{1}{2}\mathbb Z$ such that $\tau^*_z(\gamma)$ and $\tau^*_{z'}(\gamma')$ respectively have the same central character as $\tau_0^*(\lambda)$. It follows from \cite{orlik2005extensions} that
\[  \mathrm{Ext}^i_{\mathbb H_n}(\tau^*_0(\lambda), \tau^*_z(\gamma)) \cong \mathrm{Ext}^i_{\mathbb H_n}(\tau^*_0(\lambda), \tau^*_{z'}(\gamma')) \cong \begin{cases} \mathbb C & \mbox{ if $i=k-1, k$} \\ 0 & \mbox{otherwise} \end{cases} .
 \]
 This illustrates an example of Corollary \ref{cor transfer thm}.
\end{example}

\subsection{Conjectural extensions}

For two skew hooks $\lambda-\mu$ and $\lambda'-\mu'$ with \(n\) boxes, we fix the numbers $z$ and $z'$ such that $\tau^*_z(\lambda-\mu)$ and $\tau^*_{z'}(\lambda'-\mu')$ have the same central character. Define $i^*(\lambda-\mu, \lambda'-\mu')$ 
to be the smallest integer $i^*$ such that
\[   \mathrm{Ext}^{i^*}_{\mathbb H_n}(\tau^*_z(\lambda-\mu), \tau^*_{z'}(\lambda'-\mu'))
\]
is non-zero. By Theorem \ref{thm combin generalized steinberg}, such a number exists.

\begin{conjecture} \label{conj smallest ladder non vanishing}
Let $X$ and $X'$ be two uniform calibrated $\mathbb H_n$-modules with the same central character. Let $X=\tau^*_z(\lambda-\mu)$ and let $X'=\tau^*_{z'}(\lambda'-\mu')$ for some $z, z' \in \mathbb C$ and some uniform skew Young diagrams $\lambda-\mu$ and $\lambda'-\mu'$. Let $\eta_1,\ldots, \eta_d$ be the sequence of skew hooks for $\lambda-\mu$ and let $\eta'_1, \ldots, \eta_d'$ be the sequence of skew hooks for $\lambda'-\mu'$. Then the smallest integer $i^*$ such that $\mathrm{Ext}^{i^*}_{\mathbb H_n}(X, X')\neq 0$ is 
\[  i^*(\eta_1, \eta_1')+\ldots + \ldots + i^*(\eta_d, \eta_d') .
\]
\end{conjecture}

Conjecture \ref{conj smallest ladder non vanishing} is based on the validity on several examples provided below.

\begin{proposition} \label{prp conj speh w type}
Conjecture \ref{conj smallest ladder non vanishing} holds when $X$ is a Speh $\mathbb H_n$-module and $X'$ is a $W$-type $\mathbb H_n$-module.
\end{proposition}

\begin{proof}
Let $X=\tau^*_z((l^d))$ for some $l$ and $d$.

We only explain the argument for $l\geq d$. By Corollary \ref{cor: shape for same central character}, the underlying Young diagram of \(X'\) is \(((l-k)^{d}, d^k)\) for some $k$. Then the sequence of skew hooks of $(l^d)$ is
\[  \eta_1=(l, 1^{d-1}), \eta_2=(l-1, 1^{d-2}), \ldots, \eta_d=(l-d+1)  .
\]
 The sequence of the skew hooks for \( ((l-k)^d, d^k) \) are 
\[  \eta'_1=(l-k, 1^{d+k-1}), \eta'_2=(l-k-1, 1^{d+k-2}), \ldots, \eta'_d=(l-k-d+1, 1^{k}) .
\]
By Theorem \ref{thm combin generalized steinberg}, for $j=1,\ldots, d$
\[  i^*(\eta_j, \eta_j')=k .
\]
Hence, the sum
\[  i^*(\eta_1, \eta_1')+\ldots + \ldots + i^*(\eta_d, \eta_d')  
\]
is equal to $kd$. This coincides with Theorem \ref{thm: Ext for two quasi-rectangula} and hence verifies the conjecture in this case.
\end{proof}

\begin{proposition} \label{prop conj smallest ladder 1}
    Conjecture \ref{conj smallest ladder non vanishing} holds when $X$ and $X'$ are quasi-rectangular $W$-type modules.
\end{proposition}

\begin{proof}
This follows from Theorem \ref{thm: Ext for two quasi-rectangula} and an argument similar to Proposition \ref{prp conj speh w type}.
\end{proof}

%\begin{lemma}
%Let $\lambda$ be a rectangular Young diagram of $n$ boxes and let $\lambda'-\mu'$ be a uniform and connected Young diagram of $n$ boxes. Suppose furthermore that, for some $z, z' \in \mathbb C$, the $z$-content of $\lambda$ coincides with the $z'$-content of $\lambda'-\mu'$. Then the number of boxes in the first column of $\lambda'-\mu'$ is at least the number of boxes in the first column of $\lambda$.
%\end{lemma}

%\begin{proof}
%Let $z_0$ be the least value in the $z$-content of $\lambda$ (and $z'$-content of $\lambda'-\mu'$).

%Let $k$ (resp. $k'$) be the number of boxes in the first column of $\lambda$ (resp. $\lambda'-\mu'$). Then we must have that 
%\[  m(\lambda, z_0+k-2)< m(\lambda, z_0+k-1)=m(\lambda'-\mu', z_0+k) \]
%and similarly, 
%\[  m(\lambda,      \]
%\end{proof}

\begin{proposition} \label{prop conj two row speh}
Conjecture \ref{conj smallest ladder non vanishing} holds when $X$ is a Speh $\mathbb H_n$-module whose underlying Young diagram has precisely two rows, and $X'$ is a uniform calibrated $\mathbb H_n$-module.
\end{proposition}

\begin{proof}
Let $\lambda$ be the underlying rectangular Young diagram of $X$. We shall assume that the number of columns of \(\lambda\) is greater than $2$. Otherwise, the statement follows from Proposition \ref{prp conj speh w type}.

Let $\eta_1, \eta_2$ be the sequence of skew hooks for $\lambda$. Then 
\[ \eta_1=(\frac{n}{2},1 ), \quad \eta_2=(\frac{n}{2} ) . \]

Let $\eta_1', \eta_2'$ be the sequence of skew-hooks for $\lambda'-\mu'$. Note that the uniform property forces the first column of $\eta_1'$ to have at least $2$ boxes. By Theorem \ref{thm combin generalized steinberg}, $i^*(\eta_1, \eta_1')$ is equal to the number of rows minus $2$. Note that $\eta_2$ has only one row, and so Theorem \ref{thm combin generalized steinberg} implies that $i^*(\eta_2, \eta_2')$ is equal to the number of rows of \(\eta_1'\) minus $1$. 

On the other hand, Corollary \ref{cor transfer thm} provides a $W$-type module $Z$ such that 
\begin{eqnarray} \label{eqn ext isomorphism}
\mathrm{Ext}^i_{\mathbb H_n}(X, Y) \cong \mathrm{Ext}^i_{\mathbb H_n}(X, Z) . 
\end{eqnarray}
Let $\widetilde{\lambda}$ be the rectangular Young diagram associated with $Z$. Let $\widetilde{\eta}_1, \widetilde{\eta}_2$ be the sequence of skew-hooks for $\widetilde{\lambda}$. The elementary transform preserves the number of rows for the skew hooks involved in the sequence $\left\{ \widetilde{\eta}_1, \widetilde{\eta}_2\right\}$ of skew hooks. We now apply the same argument as above to compute \(i^*(\eta_1, \eta_1')\) and \(i^*(\eta_2, \eta_2')\) and hence, 
\[  i^*(\eta_1, \widetilde{\eta}_1) =i^*(\eta_1, \eta_1'), \quad i^*(\eta_2, \widetilde{\eta}_2)=i^*(\eta_2, \eta'_2) .
\]

Now, combining this with (\ref{eqn ext isomorphism}) and Proposition \ref{prp conj speh w type}, we prove the proposition.
\end{proof}

We expect that a similar argument will generalize Proposition \ref{prop conj two row speh} to the Speh module $X$ corresponding to a rectangular Young diagram of multiple rows. However, it is likely to require some other inputs or approaches for the general case of Conjecture \ref{conj non vanishing Ext}.

%\subsection{Switching process}

%Let $\mathfrak m$ be a lifted $W$-type multisegment. By {\color{blue} add}, the mutlisegment takes the form
%\[ \left\{ [a, a+k], [a-1, a-1+k], \ldots, [a-r,a-r+k], [a-r-1, a-1], \ldots, [a-r-2, a-2], \ldots , [a-r-s, a-s] \right\} . \]

%We define the Speh-completion $\overline{\mathfrak m}^{Sp}$ of $\mathfrak m$+1-k to be the multisegment 
%\[ \left\{ [a, a+k], [a-1, a-1+k], \ldots, [a-r-s, a-r-s+k]  \right\}\]
%The complemented Speh multisegment, denoted by $\overline{\mathfrak m}^{Sp}$, is the multisegment:
%\[  \left\{ [a, a+k'], [a-1, a-1+k'], \ldots, [a-s+1,a-s+1+k'] \right\},\]
%where $k'=k-(r+1)$

\bibliographystyle{alpha}
\bibliography{reference}

@article{hayashi2015decomposition,
  title={A decomposition rule for certain tensor product representations of the symmetric groups},
  author={Hayashi, Takahiro},
  journal={Journal of Algebra},
  volume={434},
  pages={46--64},
  year={2015},
  publisher={Elsevier}
}

@book{fulton1997young,
  title={Young tableaux: with applications to representation theory and geometry},
  author={Fulton, William},
  series= {London Mathematical Society Student Texts},
  number= {35},
  year={1997},
  publisher={Cambridge University Press}
}

@article{chan2016duality,
  title={Duality for Ext-groups and extensions of discrete series for graded Hecke algebras},
  author={Chan, Kei Yuen},
  journal={Advances in Mathematics},
  volume={294},
  pages={410--453},
  year={2016},
  publisher={Elsevier}
}

@incollection{barbasch2014unitary,
  title={Unitary Hecke algebra modules with nonzero Dirac cohomology},
  author={Barbasch, Dan and Ciubotaru, Dan},
  booktitle={Symmetry: Representation Theory and Its Applications: In Honor of Nolan R. Wallach},
  pages={1--20},
  year={2014},
  publisher={Springer}
}

@article{babasical1999Ktype,
  title={Classification of One K-Type Representations},
  author={Barbasch, Dan and  Moy, Allen},
  journal={Transactions of the American Mathematical Society},
  volume={294},
  pages={4145-4261},
  year={1999},
  publisher={E}
}

@article{chan2019bernstein,
  title={Bernstein--Zelevinsky derivatives: a Hecke algebra approach},
  author={Chan, Kei Yuen and Savin, Gordan},
  journal={International Mathematics Research Notices},
  volume={2019},
  number={3},
  pages={731--760},
  year={2019},
  publisher={Oxford University Press}
}

@article{grojnowski2001strong,
  title={Strong multiplicity one theorems for affine Hecke algebras of type A},
  author={Grojnowski, Ian and Vazirani, Monica},
  journal={Transformation groups},
  volume={6},
  number={2},
  pages={143--155},
  year={2001},
  publisher={Springer}
}

@article{krillo2002gha,
  title={Representations of graded Hecke algebras},
  author={Cathy Kriloff and Arun Ram},
  journal={Represent. Theory},
  volume={6}, 
  pages={31--69},
  year={2002}
}

@article{okounkov1996new,
  title={A new approach to representation theory of symmetric groups},
  author={Okounkov, Andrei and Vershik, Anatoly},
  journal={Selecta Mathematica New Series},
  volume={2},
  number={4},
  pages={581--606},
  year={1996},
  publisher={Basel; Boston: Birkhauser, c1993-1994.}
}

@book{ceccherini2010representation,
  title={Representation theory of the symmetric groups: the Okounkov-Vershik approach, character formulas, and partition algebras},
  author={Ceccherini-Silberstein, Tullio and Scarabotti, Fabio and Tolli, Filippo},
  volume={121},
  year={2010},
  publisher={Cambridge University Press}
}

@article{vershik2005new,
  title={A new approach to the representation theory of the symmetric groups. II},
  author={Vershik, Anatoli M and Okounkov, A Yu},
  journal={Journal of Mathematical Sciences},
  volume={131},
  number={2},
  pages={5471--5494},
  year={2005},
  publisher={Springer}
}

@article{lusztig1989affine,
  title={Affine Hecke algebras and their graded version},
  author={Lusztig, George},
  journal={Journal of the American Mathematical Society},
  volume={2},
  number={3},
  pages={599--635},
  year={1989},
  publisher={JSTOR}
}

@article{qadri2025extbranch,
title={Non-tempered Ext Branching Laws for the $p$-adic General Linear Group},
author={ Qadri, Mohammed Saad},
journal={ International Mathematics Research Notices },
volume={2025},
number={1},
year={2025},
}

@article{orlik2005extensions,
  title={On extensions of generalized Steinberg representations},
  author={Orlik, Sascha},
  journal={Journal of Algebra},
  volume={293},
  number={2},
  pages={611--630},
  year={2005},
  publisher={Elsevier}
}

@article{opdamsolleveldhomo09,
author = {Opdam, Eric and Solleveld, Maarten},
title= {Homological algebra for affine Hecke algebras},
journal={Adv. Math.},
volume={220},
pages={1549-1601},
year={2009}
}

@article{opdamsolleveld10,
author ={Opdam, Eric and Solleveld, Maarten},
title= {Discrete series characters for affine Hecke algebras and their formal degrees},
journal={Acta Math.},
volume={205},
number={1},
pages={105-187},
year={2010}
}

@article{opdamsolleveldExttemp,
author={Opdam, E. and Solleveld, M.},
title={Extensions of tempered representations},
journal={Geom. Funct. Anal.},
volume={23}, 
pages={664–714},
year={2013},
}

@article{casselmanvanishing, 
author={Casselman, W.},
title={On a $p$-adic vanishing theorem of Garland},  
journal={Bull. Amer. Math. Soc.},
volume={80},
number={5}, 
pages={1001-1004},
year={1974}
}

@article{chanqarciinductiveExt,
author={Chan, K.Y. and Qadri, M.S.},
title={An inductive Ext non-vanishing theorem for the $p$-adic general linear group},
journal={arXiv:2601.20635},
year={2026}
}

@incollection{PrasadExt, 
author={Prasad, D.},
title={Homological aspects of branching laws},
series={Proceedings of Symposia in Pure Mathematics},
booktitle={The Langlands Program---III. Periods, $L$-functions, and the Relative Trace Formula.},
volune={112},
number={3},
year={2025},
publisher={Amer. Math. Soc., Providence, RI}
}

@article{barbaschciubotrapa,
author={Barbasch, D. and Ciubotaru, D. and Trapa, P.}, 
title={Dirac cohomology for graded affine Hecke algebras}, 
journal={Acta Math.},
volume={202}, 
number={2},
year={2012},  
pages={197-227}
}

@article{chanhomobranch,
author={Chan, K.Y.},
title={Homological branching law for $(\mathrm{GL}_{n+1}(F),\mathrm{GL}_n(F))$: projectivity and indecomposability},
journal={Invent. Math.},
year={2021},
pages={299–345}
}

@incollection{barabschciubotaruladder, 
author={Barbasch, D. and Ciubotaru, D.},
year={2015},
title={Ladder representations of $\mathrm{GL}_n(\mathbb Q_p)$},
booktitle={In: Nevins, M., Trapa, P. (eds) Representations of Reductive Groups. Progress in Mathematics},
volume={312},
publisher={Birkhäuser, Cham}
}

@article{lapidminguezladder,
author={Lapid, E. and M\'inguez, A.},
title={On a determinantal formula of Tadi\'c},
journal={American Journal of Mathematics},
volume={136}, 
number={1},
year={2014}, 
pages={111-142}
}

@article{tadic, 
author={Tadi\'c, M.},
title={Classification of unitary representations in irreducible representations of general linear group (non-
Archimedean case)},
journal={Ann. Sci. Éc. Norm. Sup\'er.},
volume={19},
number={4},
year={1986},
pages={335-382}
}

@article{barbasch2013hermitian,
  title={Hermitian forms for affine Hecke algebras},
  author={Barbasch, Dan and Ciubotaru, Dan},
  journal={arXiv preprint arXiv:1312.3316},
  year={2013}
}

@article{chen2026extension,
  title={The Extension groups of one $K$-type representations of graded Hecke algebra},
  author={Chen, Quanyong},
  journal={arXiv preprint arXiv:2611.04535},
  year={2026}
}

@article{chanhuang,
author={Chan, K.Y. and Huang, S.},
title={Discrete seires for the graded Hecke algebra of type \(H_4\)},
journal={Journal of Algebra},
volume={712},
year={2027},
pages={326-383}
}
\end{document}